\documentclass[rpreprint,3p,10pt,sort&compress, colorlinks,linkcolor=blue, citecolor=blue,fleqn]{elsarticlex}

\usepackage[utf8]{inputenc}

\usepackage{amssymb}
\usepackage{amsmath,mathrsfs}
\usepackage{color}
\usepackage[dvipsnames]{xcolor}
\usepackage{placeins}

\font\LARGEbsf=cmssdc10 scaled 2100      

\usepackage{tikz}
\usetikzlibrary{spy}

\usepackage{siunitx}

\usepackage{stfloats}
\usepackage{siunitx}

\usepackage{bm}

\usepackage{graphicx}

\usepackage{caption}
\usepackage{subcaption}

\usepackage{multicol}

\usepackage{multirow}

\usepackage{mhchem}

\usepackage[toc,page]{appendix}

\usepackage{bm}

\usepackage[displaymath, mathlines, pagewise]{lineno}

\usepackage{float}

\usepackage[makeroom]{cancel}
\usepackage{pgfplots}
\pgfplotsset{compat=1.18, width=7cm, height=5.5cm}
\usepackage{gensymb}
\usepackage{siunitx}

\usepackage[colorlinks,linkcolor=blue,anchorcolor=green,citecolor=blue]{hyperref}

\usepackage[numbers]{natbib}
\begin{document}
	
	\begin{frontmatter} 
		\title{A multi-field decomposed model order reduction approach for thermo-mechanically coupled gradient-extended damage simulations}
		
		\author[els,*]{Qinghua Zhang}
		\ead{qinghua.zhang@ifam.rwth-aachen.de}
		
		\author[els]{Stephan Ritzert}
		\author[els]{Jian Zhang}
		\author[els]{Jannick Kehls}
		\author[els,rvt]{Stefanie Reese}
		\author[els]{Tim Brepols}
		\cortext[*]{Corresponding author}
		
		\address[els]{Institute of Applied Mechanics, RWTH Aachen University, Mies-van-der-Rohe-Str. 1, 52074 Aachen, Germany}
		\address[rvt]{University of Siegen, Adolf-Reichwein-Straße 2a, 57076 Siegen, Germany}
		
		\begin{abstract}
			Numerical simulations are crucial for comprehending how engineering structures behave under extreme conditions, particularly when dealing with thermo-mechanically coupled issues compounded by damage-induced material softening. However, such simulations often entail substantial computational expenses. To mitigate this, the focus has shifted towards employing model order reduction (MOR) techniques, which hold promise for accelerating computations. Yet, applying MOR to highly nonlinear, multi-physical problems influenced by material softening remains a relatively new area of research, with numerous unanswered questions. Addressing this gap, this study proposes and investigates a novel multi-field decomposed MOR technique, rooted in a snapshot-based Proper Orthogonal Decomposition-Galerkin (POD-G) projection approach. Utilizing a recently developed thermo-mechanically coupled gradient-extended damage-plasticity model as a case study, this work demonstrates that splitting snapshot vectors into distinct physical fields (displacements, damage, temperature) and projecting them onto separate lower-dimensional subspaces can yield more precise and stable outcomes compared to conventional methods. Through a series of numerical benchmark tests, our multi-field decomposed MOR technique demonstrates its capacity to significantly reduce computational expenses in simulations involving severe damage, while maintaining a high level of accuracy.
		\end{abstract}

		\begin{keyword}
			Multi-field decomposition; Model order reduction; Proper Orthogonal Decomposition-Galerkin; Gradient-extended damage; Multiphysics; Damage softening; Finite strain.
		\end{keyword}
		
	\end{frontmatter}
	
	\section{Introduction}
	\label{Introduction}
	Predicting and analyzing damage initiation and propagation in ductile metallic materials under thermal influence is of significant interest in research and industry. Many models in the literature concentrate solely on fracture \cite{wu2017unified,amor2009regularized,miehe2010phase,fredddi2017numerical} and damage \cite{forest2009micromorphic,nguyen2020gradient,junker2022efficient} mechanisms, neglecting the influence of temperature. As a result, simulations derived from these models typically suffice only for specific cases that deviate significantly from actual manufacturing conditions. Moreover, multiphysics damage/fracture modeling \cite{ruan2023thermo,lamm2024gradient,felder2022thermo} considering temperature-dependent plastic flow is not only a challenging task \cite{dittmann2020phase,felder2022thermo}, but also usually leads to computationally expensive simulations \cite{zhuang2022phase,dinachandra2022adaptive,zhuang2023transverse}. This is mainly due to an increased number of degrees of freedom (in the case of gradient-extended damage/phase-field of fracture) and the coupling between damage, plasticity, and temperature. Consequently, applying efficient model order reduction (MOR) techniques to accelerate multiphysics models is an interesting and promising area of research, especially for industries focusing on real-time simulations.
	
	In solid mechanics, various MOR techniques have been employed to reduce the simulation cost by projecting linear and nonlinear equations onto a subspace. Notable methods include Proper Orthogonal Decomposition-Galerkin (POD-G) projection approach \cite{liang2002proper,kerfriden2012local,radermacher2013comparison}, Discrete Empirical Interpolation Method (DEIM) \cite{tiso2013discrete,radermacher2016pod}, Energy-Conserving Sampling and Weighting Method (ECSWM) \cite{farhat2015structure} as well as Hierarchical Tensor Approximation (HTA) \cite{ballani2013black,kastian2020two}, etc. Although the works mentioned above clearly show that MOR can be used very successfully in the field of solid mechanics, most studies in the literature deal with constitutive models without softening behavior, such as damage. The latter is known to cause additional complications, posing severe challenges to the application of MOR approaches, as demonstrated e.g. by \citet{kerfriden2011bridging} and \citet{kastian2023discrete}.
	
	The phenomenon can be attributed to the fact that the material's stiffness is degraded by a local damage-induced degradation function after damage initiation. This subsequently leads to a severe drop in the force-displacement curve. In certain loading scenarios, even snap-back phenomena can occur in the softening region, as noted by \citet{peerlings1996gradient} and \citet{saji2024new}, even though this snap-back phenomenon is hard to observe in experiments, see e.g. \citet{rezaei2017prediction}. Due to the complicated behavior exhibited in the damage-induced softening region, a huge parametric space is required to accurately reproduce the mechanical response under multi-axial loading conditions, as explained by \citet{selvaraj2024adaptive}.
	
	To the author's knowledge, \citet{kerfriden2011bridging} were the first to systematically address the aforementioned challenge by employing a full-order model solver with a Krylov algorithm as a corrector whenever the standard MOR solver fails to adequately reproduce damage evolution in the softening region. Then, the corrected nodal solutions are collected as enhanced snapshots to enrich the projection matrix. In light of the large computations required by the full-order model, this corrective technique requires significant computational costs. 
	Later, \citet{kerfriden2013partitioned} utilized the domain decomposed MOR method to reduce simulations involving a comparatively simple local damage model. The relative error between the solutions obtained by using the full-order and the reduced model reached around $10^{-1} \sim 10^{-4}$ \cite{kerfriden2013partitioned}. However, whether the MOR approach can reconstruct the force-displacement curve in the softening region is not clearly given in their contribution. Interestingly enough, the domain partitioned and hyper-reduction-based empirical cubature method (ECM) by \citet{hernandez2017dimensional} was employed by \citet{rocha2020accelerating} to reduce a fracture simulation. The latter authors proposed a possibility for reconstructing the force-displacement response in the softening region. However, both the full-order and reduced-order models exhibited strong oscillations. Another interesting investigation by \citet{oliver2017reduced} incorporates the material state decomposition method in conjunction with domain decomposition techniques to reduce a multi-scale fracture model simulation. In this context, MOR is applied at the integration point level and not at the global level, wherein microscopic strain snapshots are collected separately in the elastic, elastoplastic, and softening regions for the representative volume element (RVE). Using this method,  the domain is decomposed into a cohesive zone and an undamaged zone during snapshot collection.
	
	Recently, \citet{selvaraj2024adaptive} employed an adaptive and variable MOR approach for a cohesive zone fracture model employing explicit time integration. However, their approach predominantly relies on the domain decomposition method, in which full-order and reduced-order models are applied to the cohesive zone and the undamaged domain, respectively. This strategy can avoid directly solving a damage zone via MOR techniques. However, as the force-displacement curve enters the softening region, substantial fluctuations are observed, which renders the damage evolution unstable. Addressing this concern, \citet{mishra2024enhanced} have refined the methodology by incorporating Enhanced-Transformation Field Analysis (E-TFA), leading to improved results. Inspired by \citet{radermacher2016pod}, \citet{kastian2023discrete} utilized a Proper Orthogonal Decomposition-Galerkin (POD-G) projection approach in conjunction with the DEIM to approximate the non-linear stiffness matrix, thereby reducing a gradient-extended damage model. Nevertheless, it must be stated that the relative error between the reduced and full-order simulations remained rather high in the softening region. This emphasizes again the persistent challenges in effectively applying MOR in damage simulations.
	
	In light of the aforementioned observations, the primary objective of this study is to introduce a novel multi-field decomposed MOR technique and apply it exemplarily to a thermo-mechanically coupled gradient-extended damage-plasticity model. The aim is to numerically investigate its stability and effectiveness in reducing computational costs, particularly by addressing the complexities in the damage softening region. Hence, the paper is structured as follows: \hyperref[Modelling]{Section. \textit{2}} recapitulates the derivations of a thermo-mechanically coupled multi-physical gradient-extended damage model, proposed by \citet{felder2022thermo}. Afterwards, \hyperref[MOR theory]{Section. \textit{3}} introduces the novel multi-field decomposed MOR method, which involves collecting and decomposing converged nodal snapshots into displacement, damage, and temperature fields, and integrating them into the thermo-mechanically coupled multi-physical gradient-extended damage model. \hyperref[Model validation]{Section. \textit{4}} is concerned with two academic numerical tests, conducted to evaluate the performance of the proposed multi-field MOR approach. In \hyperref[Results and discussion]{Section. \textit{5}}, a comparison between the novel multi-field MOR and classical MOR is provided to explain why the former should be preferred for reducing computational costs within the damage-softening region. Furthermore, the advantages and limitations of the multi-field decomposed MOR approach are discussed in detail. Finally, the concluding remarks are given in \hyperref[Conclusion]{Section. \textit{6}}.
	
	
	\section{Material model}
	\label{Modelling}	
	
	\subsection{Kinematics and general Helmholtz free energy function}
	Following the well-known multiplicative decomposition method proposed by \citet{lu1975decomposition}, the deformation gradient $\bm{F}$ is decomposed into mechanical $\bm{F}_m$ and isothermal $\bm{F}_{\theta}$  parts \cite{stojanovic1964finite}. $\bm{F}_m$ can be further multiplicatively decomposed into elastic $\bm{F}_{e}$ and plastic $\bm{F}_{p}$ parts:  
	\begin{equation}
		\bm{F} = {\bm{F}_m} \, {\bm{F}_\theta } = \underbrace {{\bm{F}_e} \, {\bm{F}_p}}_{{\bm{F}_m}} \, {\bm{F}_\theta}.
		\label{eq1}
	\end{equation}
	\noindent Assuming isotropic thermal expansion of the material, $\bm{F}_{\theta}$ is expressed using an exponential function $\vartheta$ based on temperature $\theta$, see  \citet{stojanovic1964finite}:
	\begin{equation}
		{\bm{F}_\theta } = \vartheta \left( \theta  \right) \, \bm{I},
		\quad
		\vartheta \left( \theta  \right) =  \exp {\left( \alpha \, \left( {\theta  - {\theta _0}} \right) \right)},
		\label{eq2}
	\end{equation}
	
	\noindent where $\theta _0$ and $\bm{I}$ represent the reference temperature and the second-order identity tensor, respectively. $\alpha$ denotes the temperature-independent thermal expansion coefficient. With the kinematic equations at hand, the elastic and plastic right Cauchy-Green tensors are expressed as: 
	\begin{equation}
		{\bm{C}_e} = \bm{F}_e^{\rm{T}} \, {\bm{F}_e} = \frac{1}{{{\vartheta ^2}}} \, \bm{F}_p^{ - \rm{T}} \, \bm{C} \, \bm{F}_p^{ - 1}, \quad {\bm{C}_p} = \bm{F}_p^{\rm{T}} \, {\bm{F}_p}.
		\label{eq3}
	\end{equation}
	
	Subsequently, the Helmholtz energy density function $ \psi $ is formulated in terms of the elastic ($\bm{C}_{e}$) and plastic ($\bm{C}_{p}$) right Cauchy-Green tensors:
	\begin{equation}
		\psi  = \underbrace {{f_d}(D)\left( {{\psi _e} \, ({\bm{C}_e}, \, \theta ) + {\psi _p} \, ({\bm{C}_p}, \, {{\xi} _p}, \, \theta )} \right)}_{{\rm{Elastoplastic \, part}}}  \, + \underbrace {{\psi _d} \, ({\xi _d}, \, \theta )}_{{\rm{Damage \, hardening \, part}}} + \underbrace {{\psi _{\bar d}} \, (D, \, \bar D, \, \nabla \bar D, \, \theta )}_{{\rm{Micromorphic \, extension}}} + \underbrace {{\psi _\theta } \, (\theta )}_{{\rm{Caloric \, part}}}.
		\label{eq4}
	\end{equation}
	
	\noindent Following the hypothesis of energy equivalence by \citet{cordebois1982damage}, the quadratic degradation function $f_{d}(D)=(1-D)^2$ is employed to degrade the elastoplastic energy whenever the local damage variable $D$ evolves. $\psi _e$ represents the elastic energy density which depends on both the elastic right Cauchy-Green deformation tensor $\bm{C}_e$ and temperature $\theta$. $\psi _p$ denotes the plastic energy density which is formulated in terms of the plastic right Cauchy-Green tensor $\bm{C}_p$, the accumulated plastic strain $\xi _p$ as well as $\theta$, respectively. Furthermore, a damage-hardening energy density $\psi _d$ is postulated as a function of the damage-hardening variable $\xi _d$ and $\theta$. 
	
	As for the micromorphic extension part, $\psi _{\bar d}$ achieves a coupling between the local damage variable $D$ and non-local (micromorphic) damage variable $\bar {D}$ and introduces the influence of the gradient $\nabla {\bar D}$ (related to the reference configuration) into the formulation. Additionally, the caloric energy density $\psi_{\theta}$ does not need to be further specified for some specific temperature-dependent materials, like polymers \cite{ames2009thermo,zhang2022exploring,zhang2021molecular}, since a constant heat capacity is assumed. For reasons of brevity, no further details are provided here. Instead, the interested reader is referred to \citet{felder2022thermo}. Finally, the energy density is formulated in terms of invariants of $\bm{C}_e$. These invariants can then equivalently be reformulated in terms of $\bm{C}$ and $\bm{C}_p$, see \hyperref[eq3]{Eq.\eqref{eq3}}. The total energy density function becomes $\psi \left( \bm{C}, \,  \bm{\zeta}_{{\rm {int}}},  \, \bar{D}, \,  \nabla {\bar D}, \, \theta \right)$, where the internal variables are defined as  ${\bm{\zeta}_{{\rm {int}} }} = [{\bm{C}_p}, \, { \xi _p}, \, D, \, {\xi _d}]$. 
	
	\subsection{Specific choice for the Helmholtz free energy}
	\begin{itemize}
		\item Elastoplasticity
	\end{itemize}
	
	A compressible Neo-Hookean-type material model is chosen to formulate the elastic energy density function $\psi_e$ in terms of the elastic right Cauchy-Green tensor $\bm{C}_e$ as:
	\begin{equation}
		{\psi _e} = \frac{\mu }{2} \, \left( {{\rm{tr}}({\bm{C}_e}) - 3 - \ln \left(\det \left({\bm{C}_e} \right) \right)} \right) + \frac{\Lambda }{4} \, \left( {\det \left({\bm{C}_e} \right) - 1 - \ln \left( \det \left ({\bm{C}_e} \right) \right)} \right).
		\label{eq5}
	\end{equation}
	
	\noindent The plastic energy density function $\psi_{p}$ is expressed as follows:
	\begin{equation}
		{\psi _p} = \underbrace {\frac{a}{2} \, \left( {{\rm{tr}}({\bm{C}_p}) - 3 - \ln \left( {\det ({\bm{C}_p})} \right)} \right)}_{\rm{linear \,kinematic \, hardening}} + \underbrace { e_{p} \, \left( {{{\xi} _p} + \frac{{\exp ( - f_{p}\,{{\xi} _p}) - 1}}{f_p}} \right)}_{\text{nonlinear isotropic hardening}} + \underbrace {\frac{1}{2} \, P \,{\xi} _p^2,}_{\text{linear isotropic hardening}}
		\label{eq6}
	\end{equation}
	
	\noindent where the Lam\'e constants $\mu$ and $\Lambda$ in \hyperref[eq5]{Eq.\eqref{eq5}}, as well as the plastic material parameters $a$, $e_p$, $f_p$, and $P$ in \hyperref[eq6]{Eq.\eqref{eq6}} are temperature dependent.
	\begin{itemize}
		\item Damage hardening
	\end{itemize}
	
	The damage-hardening energy density function is chosen analogously to the nonlinear isotropic hardening of the plastic part, where $e_d$ and $f_d$ are constant material parameters related to damage hardening:
	\begin{equation}
		{\psi _d} = {e_d} \, \left( {{{\xi _d}} + \frac{{\exp ( - {f_d}\,{{\xi _d}}) - 1}}{{{f_d}}}} \right)
		\label{eq7}
	\end{equation}
	
	\begin{itemize}
		\item Micromorphic extension
	\end{itemize}
	\begin{equation}
		{\psi _{\bar d}} = \underbrace {\frac{{H}}{2} \, {{\left( {D - \bar D} \right)}^2}}_{{\text{local - nonlocal coupling}}} + \frac{{A}}{2} \, \nabla \bar D \cdot \nabla \bar D
		\label{eq8}
	\end{equation}
	
	The miromorphic extension energy density function $\psi _{\bar{d}}$ realizes a strong coupling between the local damage variable $D$ and the non-local damage variable $\bar {D}$ by means of a constant penalty parameter $H$. $A$ denotes a constant length scale parameter that incorporates the influence of the gradient of $\bar{D}$ in the formulation and suitably regularizes the latter.
	\subsection{Thermodynamical consistency}
	To obtain thermodynamically consistent constitutive relations, an extended form of the Clausius-Duhem inequality is employed:
	\begin{equation}
		\bm{S} \, {:} \, \frac{1}{2} \, \dot {\bm{C}} + \underbrace {{a_i} \, \dot {\bar D} + {\bm{b}_i} \cdot \nabla \bar D}_{\rm{micromorphic \, extension}} - \dot \psi  - \eta \, \dot \theta  - \frac{1}{\theta } \, \bm{q} \cdot \nabla \theta  \ge 0 .
		\label{eq9}
	\end{equation}
	
	\noindent For any arbitrary process, inequality \eqref{eq9} must hold such that the second law of thermodynamics is fulfilled, see \citet{gurtin2010mechanics}. Inserting the time derivative of the Helmholtz free energy density function $\psi$ (see \hyperref[eq4]{Eq.\eqref{eq4}}) into inequality \eqref{eq9}, leads after some rearrangements to inequality \eqref{eq10}. For details, the interested reader is referred to the \hyperref[AppendixB]{Appendix B}. 
	\begin{equation}
		\begin{gathered}
			{\bm{S}} \, {\text{:}} \, \frac{1}{2} \, {\dot{\bm {C}}} + \left( {{a_i} - \frac{{\partial {\psi _{\bar d}}}}{{\partial \bar D}}} \right)\dot {\bar D} + \left( {{{\bm{b}}_i} - \frac{{\partial {\psi _{\bar d}}}}{{\partial \nabla \bar D}}} \right)\nabla \dot {\bar D}   - {f_d}(D)\left[ {\frac{{\partial {\psi _e}}}{{\partial {{\bm{C}}_e}}} \, {\text{:}} \, {{{\dot{\bm {C}}}}_e} + \frac{{\partial {\psi _p}}}{{\partial {{\bm{C}}_p}}} \, {\text{:}} \, {{{\dot{\bm {C}}}}_p} + \frac{{\partial {\psi _p}}}{{\partial {\xi _p}}} \, {{\dot \xi }_p}} \right] \hfill \\
			- \left[ {\frac{{\partial {f_d}(D)}}{{\partial D}}\left( {{\psi _e} + {\psi _p}} \right) + \frac{{\partial {\psi _{\bar d}}}}{{\partial D}}} \right]\dot {D} - \frac{{\partial {\psi _d}}}{{\partial {\xi _d}}} \, {{\dot \xi }_d}  - \left( {\frac{{\partial \psi }}{{\partial \theta }} + \eta } \right)\dot \theta - \frac{1}{\theta } \, {\bm{q}} \cdot \nabla \theta  \geqslant 0 \hfill \\ 
		\end{gathered} 
		\label{eq10}
	\end{equation}
	
	\noindent Following the classical Coleman-Noll procedure and assuming zero dissipation for the miromorphic damage variable $\bar{D}$ and its gradient $\nabla \bar{D}$, the second Piola-Kirchhoff stress $\bm{S}$, as well as the micromorphic forces ($a_i$, $\bm{b}_i$) are obtained. Furthermore, the Mandel stress $\bm{M}$ and the back stress tensor $\bm{\chi}$ in the intermediate configuration are defined as:
	\begin{equation}
		{{\bm{S}} = 2 \, {f_d}(D) \, \frac{1}{{{\vartheta ^2}}} \, {\bm{F}}_p^{ - 1} \, \frac{{\partial {\psi _e}}}{{\partial {{\bm{C}}_e}}} \, {\bm{F}}_p^{ - {\text{T}}}}, \quad
		{a_i} = \frac{{\partial {\psi _{\bar d}}}}{{\partial \bar D}}, \quad \bm{b}_{i} = \frac{{\partial {\psi _{\bar d}}}}{{\partial \nabla \bar D}}, \quad
		\bm{M} :={2 \, f_{d}(D) \, {\frac{{\partial {\psi _e}}}{{\partial {{\bm{C}}_e}}}}  \,  \bm{C}_{e}  }, \quad
		\bm{\chi} :={2 \, f_{d}(D) \, {\bm{F}}_p^{\text{T}} \, \frac{{\partial {\psi _p}}}{{\partial {{\bm{C}}_p}}} \, {{\bm{F}}_p}}.
		\label{eq11}
	\end{equation}
	\noindent The entropy $\eta$ and heat flux $\bm{q}$ are determined as follows:
	\begin{equation}
		\eta  =  - \frac{{\partial \psi }}{{\partial \theta }}, \quad \bm{q} =  - {\bm{K}_{\theta}} \, \nabla \theta.
		\label{eq12}
	\end{equation}
	
	\noindent The conjugated driving forces to damage ${Y}$, isotropic hardening $q_p$, and damage hardening $q_d$ are defined, respectively, as:
	\begin{equation}
		{Y}: =  - \left( {\frac{{\partial {f_d}(D)}}{{\partial D}} \, \left( {{\psi _e} + {\psi _p}} \right) + \frac{{\partial {\psi _{\bar d}}}}{{\partial D}}} \right), \quad
		{q_p} := {f_d}(D) \, \frac{{\partial {\psi _p}}}{{\partial {\xi _p}}}, \quad
		{q_d} := \frac{{\partial {\psi _d}}}{{\partial {\xi _d}}}.
		\label{eq13}
	\end{equation}
	
	\noindent With this, the remaining part of inequality \eqref{eq10} becomes:
	\begin{equation}
		\left( {{\bm{M}} - {\bm{\chi }}} \right):{{\bm{D}}_p} - {q_p} \, {{\dot \xi }_p} - {q_d} \, {{\dot \xi }_d} + Y \, \dot D \geqslant 0,
		\label{eq14}
	\end{equation}
	\noindent where $\bm{D}_{p}=\frac{1}{2} \, \bm{F}^{{\rm{-T}}}_{p} \, \dot {\bm{C}}_{p} \, \bm{F}^{{\rm{-1}}}_{p}$ represents the symmetric part of the plastic velocity gradient. In \hyperref[eq12]{Eq.\eqref{eq12}}, $\bm{K}_{\theta}$ denotes the positive semi-definite conductivity tensor. Following \citet{dittmann2020phase}, $\bm{K}_{\theta}$ is assumed as the following function:
	\begin{equation}
		\bm{K}_{\theta} = \left( f_{d}(D) \, K_{0} + \left( 1-f_{d}(D) \right) \, K_{c} \right) \, \bm{C}^{-1}.
		\label{eq15}
	\end{equation} 
	\noindent $K_0$ and $K_c$ denote the heat conduction parameters of the undamaged state ($f_{d}(D)=1$) and broken state ($f_{d}(D)=0$), respectively. As such, the heat conduction becomes influenced by the damage state in the material. In line with \citet{dittmann2020phase}, the plastic material parameters are degraded by a thermal softening function $f_{\theta} (\theta)$ with a softening parameter $\omega$, see \citet{simo1992associative}:
	\begin{equation}
		f_{\theta}(\theta) =  1 - \omega \, \left(\theta-\theta _{0}\right).
		\label{eq16}
	\end{equation} 
	
	\subsection{Two-surface approach: plastic yield and damage loading functions}
	Next, following the \textit{two-surface} approach proposed by \citet{brepols2017gradient,brepols2018micromorphic,brepols2020gradient},  plasticity and damage are modeled as distinct but coupled phenomena. For plasticity, a von Mises-type criterion is employed. i.e. 
	\begin{equation}
		{\Phi _p} = \sqrt {\frac{3}{2}} \, \| {{\bm{\tilde M}'} - {\bm{\tilde \chi }'}} \| - \left( {{\sigma _{y0}} + {{\tilde q}_p}} \right).
		\label{eq17}
	\end{equation} 
	\noindent $\sigma _{y0}$ is the initial yield stress of the material. The operator ${\left(  \bullet  \right)^\prime } = \left(  \bullet  \right) - \tfrac{1}{3}{\text{tr}}\left(  \bullet  \right) \, \bm{I}$ in \hyperref[eq17]{Eq.\eqref{eq17}} represents the deviatoric part of a second order tensor. It can be observed from \hyperref[eq17]{Eq.\eqref{eq17}} that the plastic yield criterion is formulated in terms of effective quantities $\left( {\tilde  \bullet } \right) = \tfrac{(\bullet)}{{{f_d}\left( D \right)}}$. For damage, the damage loading function is defined as:
	\begin{equation}
		{\Phi _d} = Y - \left( {{Y_0} + {q_d}} \right).
		\label{eq18}
	\end{equation}
	
	\noindent $Y_0$ represents the damage threshold that governs the damage initiation. 
	\begin{itemize}
		\item Plastic quantities
	\end{itemize}
	\begin{equation}
		{{{\bm{D}}}_p} = {\dot \lambda _p} \, \frac{{\partial {\Phi _p}}}{{\partial {\bm{M}}}} = {\dot{\lambda} _p} \, \sqrt {\frac{3}{2}} \, \frac{1}{{{f_d}\left( D \right)}}\frac{{{\bm{\tilde M}'} - {\bm{\tilde \chi}'}}}{{\| {{\bm{\tilde M}'} - {\bm{\tilde \chi}'}} \|}} 
		\label{eq19}
	\end{equation} 
	\begin{equation}
		{\dot \xi _p} = -{\dot \lambda _p} \, \frac{{\partial {\Phi _p}}}{{\partial {q_p}}} = {\dot \lambda _p} \, \frac{1}{{{f_d}\left( D \right)}}
		\label{eq20}
	\end{equation}
	\noindent Here, $\dot {\lambda}_p$ and $\dot {\xi}_p$ are the plastic multiplier and the accumulated plastic strain rate, respectively. 
	\begin{itemize}
		\item Damage quantities
	\end{itemize}
	\begin{equation}
		\dot D = {{\dot \lambda }_d} \, \frac{{\partial {\Phi _d}}}{{\partial Y}} = {{\dot \lambda }_d}, \quad {{\dot \xi }_d} = -{{\dot \lambda }_d} \, \frac{{\partial {\Phi _d}}}{{\partial {q_d}}} = {{\dot \lambda }_d}
		\label{eq21}
	\end{equation}
	\noindent Here, $\dot{\lambda}_d$ is the damage multiplier. Due to the \textit{two-surface} nature of the model, there are two sets of Karush-Kuhn-Tucker (KKT) conditions that are formulated as:
	\begin{equation}
		\begin{gathered}
			{\Phi _p} \leqslant 0, \quad {{\dot \lambda }_p \geq 0}, \quad {\Phi _p} \, {{\dot \lambda }_p} = 0, \hfill \\
			{\Phi _d} \leqslant 0, \quad {{\dot \lambda }_d \geq 0}, \quad {\Phi _d} \, {{\dot \lambda }_d} = 0. \hfill \\ 
		\end{gathered}
		\label{eq22}
	\end{equation}
	
	\subsection{Balance equations}
	The basic constitutive equations and derivations are based on the general micromorphic approach developed by \citet{forest2009micromorphic} in order to overcome the mesh dependency issues known for conventional continuum damage models. Assuming an isothermal and quasi-static setting at finite strains, the relevant partial differential equations (PDEs) of the formulation are given as:
	\label{MABE}
	\begin{equation}
		\begin{gathered}
			{\underline{\text{Balance of linear momentum}}} \hfill \\
			\begin{array}{*{20}{c}}
				{{\text{Div}}\left( {{\bm{FS}}} \right) + \bar {\bm{f}}}&{ = \bm{0}}&{{\text{in }}\Omega } \\ 
				{\bm{FS} \cdot \bm{n}}&{ = \bar{ \bm t}}&{{\text{on }}\partial {\Omega _t}} \\ 
				\bm{u}&{ = \tilde{\bm {u}}}&{{\text{on }}\partial {\Omega _u}} 
			\end{array} \hfill \\ 
		\end{gathered}
		\label{eq26}
	\end{equation}
	\begin{equation}
		\begin{gathered}
			{\underline{\text{Micromorphic balance}}} \hfill \\
			\begin{array}{*{20}{c}}
				{{\text{Div}}\left( {\bm{b}_i} - {\bm{b}_e} \right) - {a_i} + {a_e}}&{ = {0}}&{{\text{in }}\Omega } \\ 
				{\left({\bm{b}_{i}-\bm{b}_e}\right) \cdot \bm{n}}&{ = {a_c}}&{{\text{on }}\partial {\Omega _c}} \\ 
				{\bar D}&{ = \tilde {\bar D}}&{{\text{on }}\partial {\Omega _{\bar{D}}}} 
			\end{array} \hfill \\ 
		\end{gathered}
		\label{eq27}
	\end{equation}
	\begin{equation}
		\begin{gathered}
			{\underline{\text{Energy balance}}} \hfill \\
			\begin{array}{*{20}{c}}
				{ \underbrace {- \dot e + \bm{S}:\frac{1}{2} \, \dot{\bm{C}} + {a_i} \, \dot {\bar D} + \bm{b_i} \cdot \nabla \dot {\bar D} }_{r_{\rm{int}}} + {r_{{\text{ext}}}}} - {\text{Div}}\left( \bm{q} \right)&{ = 0}&{{\text{in }}\Omega } \\ 
				{\bm{q} \cdot \bm{n}}&{ =  - \bar{q}}&{{\text{on }}\partial {\Omega _q}} \\ 
				\theta &{ = \tilde \theta }&{{\text{on }}\partial {\Omega _\theta }} 
			\end{array} \hfill \\ 
		\end{gathered}
		\label{eq28}
	\end{equation}
	
	\noindent Here, $\Omega$ denotes the body of consideration in the reference configuration. \hyperref[eq27]{Eq.\eqref{eq27}} describes the micromorphic balance in terms of the internal body forces $a_i$, $\bm{b}_i$ together with external micromorphic forces $a_e$, $\bm{b}_e$ \cite{felder2022thermo}. In this study, $a_e$ and $\bm{b}_e$ are assumed to be zero. $\bar{\bm{f}}$ and $\bm{n}$ denote the general body force vector as well as the unit outward normal vector, respectively. In \hyperref[eq28]{Eq.\eqref{eq28}}, the internal energy $e$ and corresponding time derivative $\dot e$ are obtained by a Legendre transformation $e=\psi + \theta \, \eta$ and $\dot {e}=\dot {\psi} + \dot{\theta} \, \eta + \theta \, \dot{\eta}$, respectively. The Dirichlet boundaries are denoted as $\partial \Omega _u$, $\partial \Omega _{\bar{D}}$, and $\partial \Omega _{\theta}$ with prescribed displacement $\tilde{\bm{u}}$, non-local damage $\tilde{\bar{D}}$, and temperature $\tilde{\theta}$, respectively. Analogously, Neumann boundary conditions are prescribed on $\partial \Omega _t$, $\partial \Omega _{c}$, and $\partial \Omega _{q}$ with the macroscopic traction force $\bar{\bm{t}}$, the micromorphic traction force $a_c$, and the heat flux $\bar{q}$, respectively. Then, the total internal heat generation is expressed as $r_{\rm{int}}={r_e}+{r_p}+{r_d}$ with the elastic part $r_{e}$, the plastic part $r_{p}$, and the damaged part $r_{d}$. For details, the interested reader is referred to \hyperref[AppendixA]{Appendix A}. 
	
	\subsection{Weak form of the problem}
	Having the strong forms given in \hyperref[eq26]{Eq.\eqref{eq26}}, \hyperref[eq27]{Eq.\eqref{eq27}}, and \hyperref[eq28]{Eq.\eqref{eq28}} at hand, the weak form of the formulation is derived by multiplication with appropriate test functions $\delta \bm{u}$, $\delta \bar D$, and $\delta \theta$. The residual energies of the element for displacement, damage, and temperature fields are defined as $\mathit{g} _u^e$, $\mathit{g} _d^e$, and $\mathit{g} _{\theta}^e$, respectively. 
	
	\begin{itemize}
		\item Displacement field
	\end{itemize}
	\begin{equation}
		\mathit{g} _u^e \left( {\bm{u}}, \, \bar{D}, \, \theta ,\delta {\bm{u}} \right) = \int\limits_{{\Omega _e}} {{\bm{S}} \, \, {:} \, \delta {\bm{E}} \, dV}  - \int\limits_{{\Omega _e}} {{\bar{\bm{f}}} \cdot \delta {\bm{u}} \, dV}  - \int\limits_{\partial {\Omega _e}} {{\bar{\bm{t}}} \cdot \delta {\bm{u}} \, dA}
		\label{eq29}
	\end{equation}
	
	\begin{itemize}
		\item Damage field
	\end{itemize}
	\begin{equation}
		\mathit{g} _d^e \left( {\bm{u}}, \, \bar{D}, \, \theta , \, \delta {\bar{D}} \right) = \int\limits_{{\Omega _e}} {H \, \left( {D - \bar D} \right) \, \delta \bar D \, dV}  - \int\limits_{{\Omega _e}} {A \, \nabla \bar D \cdot \nabla (\delta \bar {D})  \, dV}
		\label{eq30}
	\end{equation}
	
	\newpage
	\begin{itemize}
		\item Temperature field
	\end{itemize}
	\begin{equation}
		\mathit{g} _\theta ^e ( {\bm{u}}, \, \bar D, \, \theta , \, \dot \theta , \, \delta {\theta} ) = \int\limits_{{\Omega _e}} {c \, \dot \theta \, \delta \theta \, dV}  - \int\limits_{{\Omega _e}} {{\bm{q}} \cdot \nabla  (\delta \theta) \, dV}  - \int\limits_{{\Omega _e}} {{r_{{\text{int}}}} \, \delta \theta \, dV}  - \int\limits_{\partial \Omega _e} {\bar q \, \delta \theta \, dA}
		\label{eq31}
	\end{equation}
	
	\noindent Here, $\partial \Omega _e$ represents the relevant Dirichlet ($\partial \Omega_{e, \theta}$) and Neumann ($\partial \Omega_{e, t}$) boundaries of an element $e$, with $\partial \Omega_{e}=\partial \Omega_{e, \theta} \cup \partial \Omega_{e, t}$. For brevity, a unified boundary $\partial \Omega _e$ is employed here. $\delta \bm{E}$ denotes the test function of the Green-Lagrange strain tensor $\delta \bm{E} = \frac{1}{2} \, [\bm{F}^{{\rm{T}}} \, \nabla (\delta \bm{u}) + \nabla ^{\rm{T}} (\delta (\bm{u})) \, \bm{F}]$. The element residual force $\bm{R}_e$ is derived from the virtual energy $\mathit{g}^{e}=\mathit{g} _u^{e}+\mathit{g} _d^{e}+\mathit{g} _{\theta}^{e} $ with respect to all degrees of freedom ${\bm{U}_e} = \left[ {{\bm{u}}, \, {\bar{D}}, \, {\theta} } \right]$, named displacement $\bm{u}$, non-local damage ${\bar{D}}$, and temperature  ${\theta}$:
	\begin{equation}
		\bm{R}_e = \frac{{\partial {\mathit{g}^{e} \, (\bm{u},\, \bar{D}, \, \theta)}}}{{\partial \bm{U}_e}}, \quad \bm{K}^{e} = \frac{{\partial \bm{R}_e \, (\bm{u}, \, \bar{D}, \, \theta)}}{{\partial \bm{U}_e}}.
		\label{eq32}
	\end{equation}
	
	\noindent Specifically, after taking the partial derivative of the internal energy $\mathit{g}^{e}_{\rm{int}}$ and the external energy $\mathit{g}^{e}_{\rm{ext}}$ concerning the nodal degrees of freedom, the residual force is expanded as $\bm{R}_{e}=\bm{F}^{e}_{\rm{int}}-\bm{F}^{e}_{\rm{ext}}$. $\bm{F}^{e}_{\rm{int}}$ and $\bm{F}^{e}_{\rm{ext}}$ denote the internal and external force vectors of the element, respectively. The detailed parts of the element stiffness matrix $\bm{K}^{e}$ are expressed as:
	\begin{equation}
		{\bm{K}}_{uu}^e = \frac{{{\partial ^2}{\mathit{g} ^e}}}{{\partial {{\bm{u}}^2}}}, \, {\bm{K}}_{u\bar D}^e = \frac{{{\partial ^2}{\mathit{g} ^e}}}{{\partial {\bm{u}} \partial {{\bar{D}}} }}, 
		\, {\bm{K}}_{u\theta }^e = \frac{{{\partial ^2}{\mathit{g} ^e}}}{{\partial {\bm{u}}\partial {{\theta}}  }}, \, {{K}}_{\bar D\bar D}^e = \frac{{{\partial ^2}{\mathit{g} ^e}}}{{\partial {{{{\bar{D}}} }^2}}},  \,  {{K}}_{\bar D\theta }^e = \frac{{{\partial ^2}{\mathit{g} ^e}}}{{\partial {{\bar{D}}}  \partial {{\theta}}  }}, \, {{K}}_{\theta \theta }^e = \frac{{{\partial ^2}{\mathit{g} ^e}}}{{\partial {{{\theta}}  ^2}}}.
		\label{eq33}
	\end{equation}
	\noindent $\bm{K}^{e}_{\bar{D}u}$, $\bm{K}^{e}_{\theta u}$, and ${K}^{e}_{\theta \bar{D}}$ are obtained analogously. The derivatives are computed by an automated differential toolbox: \textit{AceGen} \cite{korelc2002multi,korelc2016automation}. Subsequently, assembling individual element residual vectors $\bm{R}_e$, stiffness matrices $\bm{K}^e$, and increment vectors of unknowns $\Delta \bm{U}_e$, the global terms are obtained:
	\begin{equation}
		{{\bm{R} = \raise3pt
				\hbox{$\hbox{\scriptsize $N_e$}\atop{\hbox{\LARGEbsf A}\atop {\scriptstyle e=1}}$} {\bm{R}_e}}, \quad {\bm{K} = \raise3pt
				\hbox{$\hbox{\scriptsize $N_e$}\atop{\hbox{\LARGEbsf A}\atop {\scriptstyle e=1}}$} {\bm{K}^{e}}}
			, \quad {\Delta \bm{U} = \raise3pt
				\hbox{$\hbox{\scriptsize $N_e$}\atop{\hbox{\LARGEbsf A}\atop {\scriptstyle e=1}}$} {\Delta \bm{U}_e}}
		},   
		\label{eq34}
	\end{equation}
	\noindent where operator $\hbox{$\hbox{\scriptsize $N_e$}\atop{\hbox{\LARGEbsf A}\atop {\scriptstyle e=1}}$} \left( \bullet \right)$ in \hyperref[eq34]{Eq.\eqref{eq34}} denotes the 
	well-known assembly operator applied to all $N_e$ elements in the given domain. Finally, the global stiffness matrix $\bm{K}$, the residual force vector $\bm{R}$, and the nodal solution increment vector $\Delta \bm{U}$ are expressed as:
	\begin{equation}
		\underbrace {\left[ {\begin{array}{*{20}{c}}
					{{{\bm{K}}_{uu}}}&{{{\bm{K}}_{u\bar D}}}&{{{\bm{K}}_{u\theta }}} \\ 
					{{{\bm{K}}_{\bar Du}}}&{{{\bm{K}}_{\bar D\bar D}}}&{{{\bm{K}}_{\bar D\theta }}} \\ 
					{{{\bm{K}}_{\theta u}}}&{{{\bm{K}}_{\theta \bar D}}}&{{{\bm{K}}_{\theta \theta}}} 
			\end{array}} \right]}_{\bm{K}}\underbrace {\left\{ {\begin{array}{*{20}{c}}
					{\Delta {\bm{u}}} \\ 
					{\Delta \bar{\bm{D}}} \\ 
					{\Delta {\bm{\theta }}} 
			\end{array}} \right\}}_{\Delta {\bm{U}}} =  - \underbrace {\left\{ {\begin{array}{*{20}{c}}
					{{{\bm{R}}_u}} \\ 
					{{{\bm{R}}_{\bar D}}} \\ 
					{{{\bm{R}}_\theta }} 
			\end{array}} \right\}}_{\bm{R}}.
		\label{eq35}
	\end{equation}
	
	\section{Model order reduction}
	\label{MOR theory}
	Having $\bm{R}$ (see \hyperref[eq34]{Eq.\eqref{eq34}}) at hand, the classical approach is to solve it iteratively using the Newton-Raphson method. However, this can be computationally expensive due to the high dimension of $\bm{R}$. For this reason, MOR is employed to accelerate the computations.
	\subsection{Proper orthogonal decomposition} 
\begin{figure}[!ht]
		\centering
		\includegraphics[width=16cm]{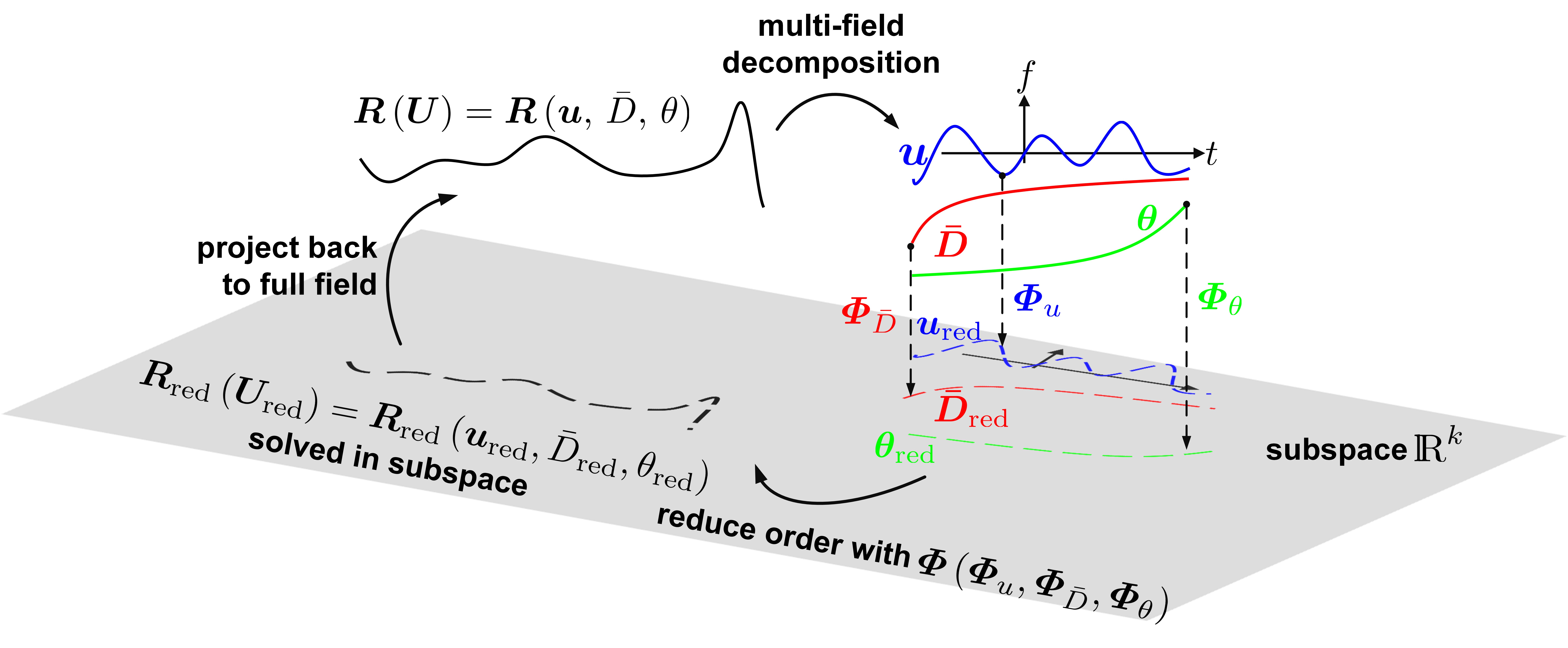}
		\caption{Schematic illustration of the multi-field decomposed MOR process and the corresponding projection. The solid and dashed lines represent the functions of the full-order model (full space) and the reduced-order model (reduced subspace), respectively. The process starts from the full residual vector $\bm{R}$ decomposition and finishes at backward projection. Then, the full nodal solution $\bm{U}$ is decomposed into a displacement field ($\bm{u}$, blue), non-local damage field (${\bar{D}}$, red), and temperature field ($\theta$, green). $\bm{u}_{\text{red}}$, $\bar{D}_{\text{red}}$, and $\theta_{\text{red}}$ denote the reduced displacement, non-local damage, and temperature, respectively. $\bm{\Phi}$ denotes the assembled projection matrix which consists of individual projection matrices $\bm{\Phi}_{u}$, $\bm{\Phi}_{\bar{D}}$, and $\bm{\Phi}_{\theta}$ for the displacement, non-local damage, and temperature fields, respectively. $\bm{R}_{\text{red}}$ is the reduced residual vector in terms of the reduced quantities.}
		\label{fig0}
	\end{figure}

	Proper orthogonal decomposition (POD \cite{breuer1991use,chatterjee2000introduction,holmes2012turbulence}) based model order reduction has been widely used for many problems, such as bio-mechanical problems \cite{radermacher2016pod,geelen2023operator}, heat transfer problems \cite{georgaka2020hybrid,georgaka2018parametric}, topology optimization \cite{nguyen2020three}, machine learning  \cite{chen2021physics,zhang2022hidenn}, etc. The main idea is to seek a lower-dimensional orthonormal projection matrix $\bm{\Phi } = \left[ {{{\bm{S}}_1},{{\bm{S}}_2}, \ldots, {{\bm{S}}_k}} \right] \in \text{I\!R}^{n \times k}$ which projects the global residual vector $\bm{R}$ and stiffness matrix $\bm{K}$ (see \hyperref[eq35]{Eq.\eqref{eq35}}) onto subspaces $\text{I\!R}^{k}$ and $\text{I\!R}^{k \times k}$, respectively, where $\text{I\!R}^{k} \subset \text{I\!R}^n$, $k \ll n$, see \hyperref[fig0]{Fig.1}. 
	
	This can be done via POD-based linear projection. In this study, nodal solutions are collected as a data matrix ${\bm{U}_{\rm{snapshot}}} = [ 
	\bm{U}^{1}_1, \ldots , \bm{U}^{n_{\rm{step}}}_i ] \in {\text{I\!R}^{n \times n_{\rm{step}}}}$ by performing $n_{\rm{step}}$ simulations with the full-order model in advance. $n$  represents the total number of degrees of freedom (DOF). $\bm{U}^{n_{\rm{step}}}_{i}$ denotes the nodal data of the $i^{th}$ iteration at time step $n_{\rm{step}}$. A minimized $\bm{\Phi}$ is sought which can project the full-order model nodal solution vector $\bm{U}^{j}_{i}$ onto the lower-dimensional subspace. Therefore, the minimization process is formulated as follows:  
	\begin{equation}
		\mathop {{\text{Arg min}}}\limits_{\{ {{\bm{S}}_1}, \, \cdots , \, {{\bm{S}}_j}\} } {\sum\limits_{j = 1}^{{n_{{\text{step}}}}} {\| {{\bm{U}}_i^j - \underbrace {\sum\limits_{i = 1}^k {( {{\bm{S}}_j^{\text{T}} \, {\bm{U}}_i^j} ) \, {{\bm{S}}_j}} }_{{\bm{U}}_{i{\text{, project}}}^j}} \|^2 _2} } \quad \xrightarrow{{{\text{subject to}}}} \quad{\bm{S}}_a^{\text{T}} \, {{\bm{S}}_b} = \delta _{ab},
		\label{eq36}
	\end{equation}
	\noindent where $\delta _{ab}$ denotes the Kronecker-Delta, with ($\delta _{ab}=1$, if $a=b$; $\delta _{ab}=0$, if $a \neq b$). $\bm{S}_j$ denotes the corresponding projection vector for the snapshot at step $j$. From \hyperref[eq36]{Eq.\eqref{eq36}}, a projection matrix can be sought in order to obtain the minimal difference between the nodal solution snapshot $\bm{U}^{j}_{i}$ and its projection term $\bm{U}^{j}_{i, \, \rm{project}}$. This is simplified to a minimization of the norm square problem. Hence, the error between full-order and reduced solution vectors is given as the sum of squared singular values $\sigma _j^2$ from the truncating mode number $k+1$ to the end $r={\rm{rank}} \, (\bm{U}_{\rm{snapshot}})$.
	\begin{equation}
		{\sum\limits_{j = 1}^{{n_{{\text{step}}}}} {\| {{\bm{U}}_i^j - {\sum\limits_{i = 1}^k {( {{\bm{S}}_j^{\text{T}} \, {\bm{U}}_i^j} ) \, {{\bm{S}}_j}} } } \|^2 _2} } = {\sum_{j=k+1}^{r} {\sigma _{j}^{2}}}
		\label{eq37}
	\end{equation}
	
	Based on \citet{volkwein2013proper}, \citet{golub2013matrix} as well as \citet{alter2000singular}, the solution of \hyperref[eq37]{Eq.\eqref{eq37}} can be obtained from the singular values of the singular value decomposition ($\mathbf{SVD}$) operator over the snapshots matrix $\bm{U}_{\rm{snapshot}}$. The singular values of ${\bm{U}_{\rm{snapshot}}}$ are denoted and ordered as ${\sigma _1} \geqslant {\sigma _2} \cdots {\sigma _k} \cdots {\sigma _r} \geqslant 0$. Here, the index $k$ represents a truncated column number (number of modes) of the projection matrix. Therefore, the decomposition of the snapshot matrix is expressed as follows:
	\begin{equation}
		{\bm{S} \, \bm{\Sigma}} \,{{\bm{V}}^{\text{T}}} = {\bm{U} _{{\text{snapshot}}}} .
		\label{eq38}
	\end{equation}
	\noindent Here, ${{\bm{S}}} = \left[ {{{\bm{S}}_1}, \ldots {{\bm{S}}_k}, \ldots ,{{\bm{S}}_r}} \right] \in {\text{I\!R}^{n \times r}}$ is a left orthonormal singular basis matrix. Matrices $\bm{S}$ and $\bm{V}$ are required to fulfill the basic orthogonality condition ${\bm{S}}^{\text{T}} \,  {{\bm{S}}} = {{\bm{V}}^{\text{T}}} \, {\bm{V}} = {\bm{I}}$ on the subspace. ${\bm{\Sigma }}$ denotes the diagonal matrix of singular values in decreasing order, expressed as ${\bm{\Sigma }} = {\text{diag}}\left( {{\sigma _1}, \ldots {\sigma _k}, \ldots,{\sigma _r}} \right) \in {\text{I\!R}^{r \times r}}$. $\bm{V}^\text{T}$ is the right orthonormal singular basis matrix, expressed as ${\bm{V}}^{\text{T}} = \left[ {{{\bm{V}}_1},{{\bm{V}}_2}, \ldots ,{{\bm{V}}_{n_{\rm{step}}}}} \right] ^{\text{T}} \in {\text{I\!R}^{r \times {n_{\rm{step}}}}}$. By having $\bm{S}$ at hand, the left orthonormal singular basis matrix is truncated at a user-selected column (mode number) $k$ in order to obtain the projection matrix ${\bm{\Phi }} = \left[ {{{\bm{S}}_1},{{\bm{S}}_2}, \ldots, {{\bm{S}}_k}} \right] \in {\text{I\!R}^{n \times k}}$, where $k \ll n$. 
	The converged full-order nodal solution vector $\bm{U}$ and its increment $\Delta {\bm{U}}$ are approximated by projecting the reduced solution vector $\bm{U}_{\rm{red}}$ and its increment $\Delta \bm{U}_{\rm{red}}$ to the full system via: 
	\begin{equation}
		{\bm{U}} \approx {\bm{\Phi }} \, {{\bm{U}}_{{\text{red}}}}, \quad \Delta {\bm{U}} \approx {\bm{\Phi }} \, \Delta {{\bm{U}}_{{\text{red}}}}.
		\label{eq39}
	\end{equation}
	
	\subsection{Novel multi-field decomposed model order reduction}
	The main idea of the novel multi-field decomposed MOR approach for thermo-mechanically coupled gradient-extended damage is to decompose the nodal solution vector into three individual fields, see \hyperref[eq40]{Eq.\eqref{eq40}}, namely the displacement field $\bm{\hat{u}}_{u}$, damage field $\bm{\hat{u}}_{\bar{D}}$, and temperature field $\bm{\hat{u}}_{\theta}$.
	\begin{figure}[!ht]
		\centering
		\includegraphics[width=14cm]{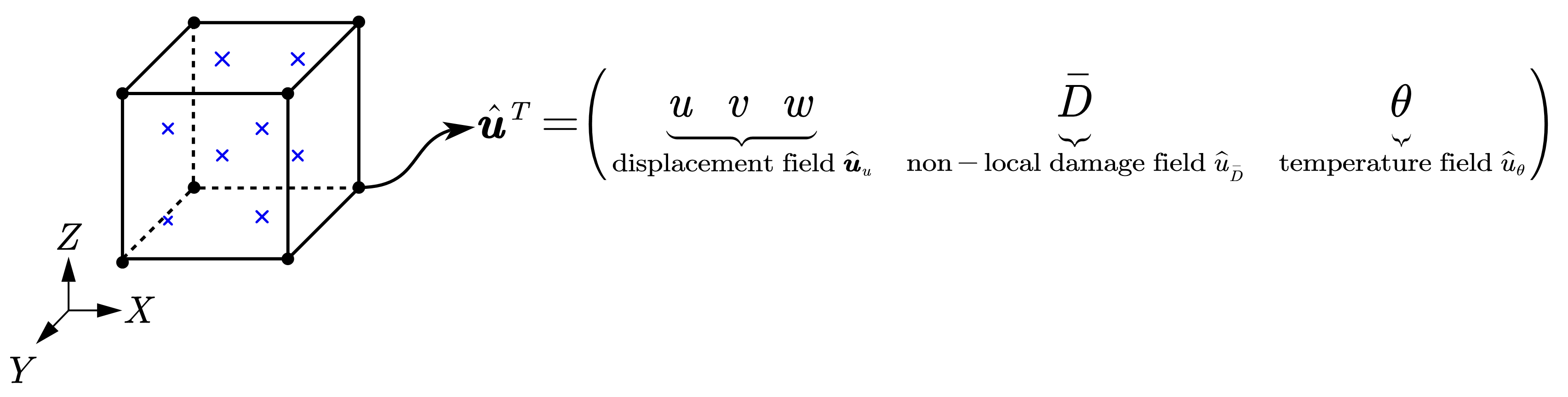}
		\caption{Schematic illustration of the nodal degrees of freedom of a linear 8-node hexahedral element. The black points and blue cross markers are the nodal points and the Gaussian integration points, respectively.}
		\label{fig1}
	\end{figure}
	
	From the dimensional perspective, the displacement field $\hat{\bm{u}}_u$ could be further decomposed into three fields i.e. ${u}$, ${v}$, and ${w}$ (the displacements in $X$, $Y$, and $Z$ direction). However, such decomposition is deemed unnecessary and will not be pursued here. Nonetheless, developing a reduced-order model for cases in which damage
	occurs poses significantly greater difficulties, since the local damage $D$  leads to a strongly non-linear behavior of the material model, see \hyperref[eq4]{Eq.\eqref{eq4}}, \hyperref[eq11]{Eq.\eqref{eq11}}, and \hyperref[eq15]{Eq.\eqref{eq15}}. The monolithically increasing stress tensors (stress-like quantities) are degraded by the degradation function. For the reduced-order model, it is a big challenge to capture the damage-induced highly non-linear softening behavior in the reduced subspace $\text{I\!R}^{k}$. 
	As also indicated by \citet{kerfriden2011bridging} and \citet{kastian2023discrete}, this is probably one of the reasons why MOR for problems involving material damage has hardly been dealt with in the literature so far. Instead, most works concentrate on  hyper-elasticity \cite{radermacher2016pod,bhattacharjee2016nonlinear}, plasticity \cite{lange2024monolithic}, viscoelasticity \cite{radermacher2016pod}, and viscoplasticity  \cite{ghavamian2017pod}, etc. To address the above-mentioned challenge, a novel multi-field decomposed MOR approach is proposed in the present work by splitting the nodal solution vector into displacement, non-local damage, and temperature parts as follows: 
	\begin{equation}
		{\hat {\bm{u}}} = \left( {\begin{array}{*{20}{c}}
				u \\ 
				v \\ 
				w \\ 
				{\bar D} \\ 
				\theta  
		\end{array}} \right) = \underbrace {\left( {\begin{array}{*{20}{c}}
					u \\ 
					v \\ 
					w \\ 
					0 \\ 
					0 
			\end{array}} \right)}_{{{\hat{\bm{u}}}_u}} + \underbrace {\left( {\begin{array}{*{20}{c}}
					0 \\ 
					0 \\ 
					0 \\ 
					{\bar D} \\ 
					0 
			\end{array}} \right)}_{{\hat{\bm{u}}_{\bar D}}} + \underbrace {\left( {\begin{array}{*{20}{c}}
					0 \\ 
					0 \\ 
					0 \\ 
					0 \\ 
					\theta  
			\end{array}} \right)}_{{\hat{\bm{u}}_\theta }}.
		\label{eq40}
	\end{equation}
	Hence, the global solution vector is decomposed into a displacement field ${{\bm{U}}_u} = {\left[ {{\bm{\hat u}}_u^1, \, {\bm{\hat u}}_u^2, \, \ldots  ,{\bm{\hat u}}_u^{\nu}} \right]^{\text{T}}}$, non-local damage field ${{\bm{U}}_{\bar D}} = {\left[ {{\bm{\hat u}}_{\bar D}^1, \, {\bm{\hat u}}_{\bar D}^2, \, \ldots  , \, {\bm{\hat u}}_{\bar D}^{\nu}} \right]^{\text{T}}}$, and temperature field ${{\bm{U}}_\theta } = {\left[ {{\bm{\hat u}}_\theta ^1, \, {\bm{\hat u}}_\theta ^2, \,  \ldots  , \, {\bm{\hat u}}_\theta ^{\nu}} \right]^{\text{T}}}$, respectively. $\nu$ represents the total number of nodes in the meshed specimen. There are two reasons for decomposing the nodal solution vector $\hat{\bm{u}}$ according to \hyperref[eq40]{Eq.\eqref{eq40}}. First, the split in \hyperref[eq40]{Eq.\eqref{eq40}} is reasonable because it effectively decouples the different field quantities and thus makes individual projections for them possible. Second, filling the vectors with zeros is beneficial from the implementational point of view because the dimensions of the vectors and matrices for the different field quantities remain the same in the MOR approach.
	Subsequently, the decomposed snapshots are defined as ${\bm{U}}_{{\text{snapshot}}}^u$ (displacement snapshots), ${\bm{U}}_{{\text{snapshot}}}^{\bar D}$ (damage snapshots), and ${\bm{U}}_{{\text{snapshot}}}^\theta$ (temperature snapshots) as follows:
	\begin{equation}
		{\bm{U}}_{{\text{snapshot}}}^u = \left[ {{\bm{U}}_u^1, \ldots {\bm{U}}_u^k, \ldots, {\bm{U}}_u^r} \right], \quad {\bm{U}}_{{\text{snapshot}}}^{\bar D} = \left[ {{\bm{U}}_{\bar D}^1, \ldots {\bm{U}}_{\bar D}^l, \ldots, {\bm{U}}_{\bar D}^r} \right], \quad {\bm{U}}_{{\text{snapshot}}}^\theta  = \left[ {{\bm{U}}_\theta ^1, \ldots {\bm{U}}_\theta ^m, \ldots, {\bm{U}}_\theta ^r} \right].
		\label{eq41}
	\end{equation}
	\noindent The indices $k$, $l$, and $m$ are individually chosen truncation numbers, i.e. the number of displacement modes, damage modes, and temperature modes, respectively, defining the dimensions of the projection matrices separately for each field. By employing the singular value decomposition $\mathbf{SVD} (\bullet)$ for the different snapshot matrices, one obtains
	\begin{equation}
		\left\{ {{{\bm{S}}_u}, \, {{\bm{\Sigma }}_u}, \, {\bm{V}}_u^{\text{T}}} \right\} = {\mathbf{SVD}}\left( {{\bm{U}}_{{\text{snapshot}}}^u} \right), \,
		\left\{ {{{\bm{S}}_{\bar D}}, \, {{\bm{\Sigma }}_{\bar D}}, \, {\bm{V}}_{\bar D}^{\text{T}}} \right\} = {\mathbf{SVD}}\left( {{\bm{U}}_{{\text{snapshot}}}^{\bar D}} \right), \,
		\left\{ {{{\bm{S}}_\theta }, \, {{\bm{\Sigma }}_\theta }, \, {\bm{V}}_\theta ^{\text{T}}} \right\} = {\mathbf{SVD}}\left( {{\bm{U}}_{{\text{snapshot}}}^\theta } \right).
		\label{eq42}
	\end{equation}
	Afterwards, the projection matrices of the displacement field $\bm{\Phi}_{u, \, k}$, non-local damage field $\bm{\Phi}_{\bar{D}, \, l}$, and temperature field $\bm{\Phi}_{\theta, \, m}$ are computed by truncating $\bm{S}_u$, $\bm{S}_{\bar {D}}$, and $\bm{S}_{\theta}$ at $k$, $l$, and $m$, respectively:
	\begin{equation}
		{{\bm{S}}_u} = \left[ {\begin{array}{*{20}{c}}
				{\underbrace {{\bm{S}}_u^1, \cdots {\bm{S}}_u^k}_{{{\mathbf{\Phi }}_{u, \, k}}},}& \cdots , &{{\bm{S}}_u^r} 
		\end{array}} \right], \quad {{\bm{S}}_{\bar D}} = \left[ {\begin{array}{*{20}{c}}
				{\underbrace {{\bm{S}}_{\bar D}^1, \cdots {\bm{S}}_{\bar D}^l}_{{{\mathbf{\Phi }}_{\bar D, \, l}}},}& \cdots , &{{\bm{S}}_{\bar D}^r} 
		\end{array}} \right], \quad {{\bm{S}}_\theta } = \left[ {\begin{array}{*{20}{c}}
				{\underbrace {{\bm{S}}_\theta ^1, \cdots {\bm{S}}_\theta ^m}_{{{\mathbf{\Phi }}_{\theta , \, m}}},}& \cdots , &{{\bm{S}}_\theta ^r} 
		\end{array}} \right].
		\label{eq43}
	\end{equation}
	
	\noindent Finally, the total projection matrix ${\mathbf{\Phi }} \in {\text{I\!R}^{n \times (k + l + m)}}$ is obtained with the mixed mode numbers $k$, $l$, and $m$ for the three fields:
	\begin{equation}
		{\mathbf{\Phi }} = \left[ {\begin{array}{*{20}{c}}
				{\underbrace {{{\mathbf{\Phi }}_{u, \, k}}}_{{\text{displacement field}}},}&{\underbrace {{{\mathbf{\Phi }}_{\bar D, \, l}}}_{{\text{non-local damage field}}},}&{\underbrace {{{\mathbf{\Phi }}_{\theta , \, m}}}_{{\text{temperature field}}}} 
		\end{array}} \right]
		\label{eq44}
	\end{equation}
	
	\noindent To have a meaningful reduced-order model, it should hold $(k + l + m) \ll n$. 
	For simulations, the nonlinear discrete residual equation $\bm{R}=\bm{0}$ is solved by means of Newton-Raphson method, see \hyperref[eq35]{Eq.\eqref{eq35}}. Therefore, the reduced residual vector $\bm{R}_{\rm{red}}$,  the reduced stiffness matrix $\bm{K}_{\rm{red}}$, and the reduced solution increment $\Delta {{\bm{U}}_{{\text{red}}}}$ are computed as: 
	\begin{equation}
		{{\bm{R}}_{{\text{red}}}}{\text{ = }}\left[ {\begin{array}{*{20}{c}}
				{\bm{\Phi }}_{u, \, k}^{\text{T}} \, {{\bm{R}}_u} \\ [0.1cm]
				{\bm{\Phi }}_{\bar D, \, l}^{\text{T}} \, {{\bm{R}}_{\bar D}} \\[0.1cm] 
				{\bm{\Phi }}_{\theta , \, m}^{\text{T}} \, {{\bm{R}}_\theta } 
		\end{array}} \right], \quad 
		{{\bm{K}}_{{\text{red}}}}{\text{ = }}\left[ {\begin{array}{*{20}{c}}
				{{\bm{K}}_{uu}^{{\text{red}}}}&{{\bm{K}}_{u\bar D}^{{\text{red}}}}&{{\bm{K}}_{u\theta }^{{\text{red}}}} \\ [0.1cm]
				{{\bm{K}}_{\bar Du}^{{\text{red}}}}&{{\bm{K}}_{\bar D\bar D}^{{\text{red}}}}&{{\bm{K}}_{\bar D\theta }^{{\text{red}}}} \\[0.1cm] 
				{{\bm{K}}_{\theta u}^{{\text{red}}}}&{{\bm{K}}_{\theta \bar D}^{{\text{red}}}}&{{\bm{K}}_{\theta \theta }^{{\text{red}}}} 
		\end{array}} \right], \quad
		\Delta {{\bm{U}}_{{\text{red}}}}{\text{ = }}\left[ {\begin{array}{*{20}{c}}
				{\bm{\Phi }}_{u, \, k}^{\text{T}} \, \Delta {{\bm{u}}} \\ [0.1cm]
				{\bm{\Phi }}_{\bar D, \, l}^{\text{T}} \, \Delta {\bar{\bm D}} \\ [0.1cm]
				{\bm{\Phi }}_{\theta , \, m}^{\text{T}} \, \Delta {\bm{\theta }} 
		\end{array}} \right].
		\label{eq45}
	\end{equation}
	The individual terms of the reduced stiffness matrices are given as follows:
	\begin{equation}
		\begin{aligned}
			{\bm{K}}_{uu}^{{\text{red}}} & = {\bm{\Phi }}_{u, \, k}^{\text{T}} \, {\bm{K}} \, {{\bm{\Phi }}_{u, \, k}}, \quad  {\bm{K}}_{u\bar D}^{{\text{red}}} = {\bm{\Phi }}_{u, \, k}^{\text{T}} \, {\bm{K}} \, {{\bm{\Phi }}_{\bar D, \, l}}, \quad {\bm{K}}_{u\theta }^{{\text{red}}} = {\bm{\Phi }}_{u, \, k}^{\text{T}} \, {\bm{K}} \, {{\bm{\Phi }}_{\theta , \, m}}, \\[0.1cm]
			{\bm{K}}_{\bar Du}^{{\text{red}}} & = {\bm{\Phi }}_{\bar D, \, l}^{\text{T}} \, {\bm{K}} \, {{\bm{\Phi }}_{u, \, k}}, \quad  {\bm{K}}_{\bar D\bar D}^{{\text{red}}} = {\bm{\Phi }}_{\bar D, \, l}^{\text{T}} \, {\bm{K}} \, {{\bm{\Phi }}_{\bar D, \, l}}, \quad {\bm{K}}_{\bar D\theta }^{{\text{red}}} = {\bm{\Phi }}_{\bar D, \, l}^{\text{T}} \, {\bm{K}} \, {{\bm{\Phi }}_{\theta , \, m}}, \\[0.1cm]
			{\bm{K}}_{\theta u}^{{\text{red}}} & = {\bm{\Phi }}_{\theta , \, m}^{\text{T}} \, {\bm{K}} \, {{\bm{\Phi }}_{u, \, k}}, \quad {\bm{K}}_{\theta \bar D}^{{\text{red}}} = {\bm{\Phi }}_{\theta , \, m}^{\text{T}} \, {\bm{K}} \, {{\bm{\Phi }}_{\bar D, \, l}}, \quad {\bm{K}}_{\theta \theta }^{{\text{red}}} = {\bm{\Phi }}_{\theta , \, m}^{\text{T}} \, {\bm{K}} \, {{\bm{\Phi }}_{\theta , \, m}}.
		\end{aligned}
	\end{equation}
	
	\noindent The main task of the finite element method (FEM) is to compute a suitable solution vector $\bm{U}=[\bm{u}, \, \bm{\bar{D}}, \, \bm{\theta}]$ by iteratively solving the quasi-static discrete nonlinear residual equation $\bm{R} \, (\bm{U})=\bm{0}$:
	\begin{equation}
		\begin{gathered}
			{\bm{R}} \, ({{\bm{U}}_{i + 1}^j}) \approx {\bm{R}} ( {{\bm{U}}_i^j}) +  {\bm{K}} ( {{\bm{U}}_i^j})  \, \Delta {\bm{U}}_{i + 1}^j = \bm{0}.
		\end{gathered}
		\label{eq47-1}
	\end{equation}
	The reduced equations are expressed as:  	
	\begin{equation}
		\begin{gathered}
			{{\bm{\Phi }}^{\text{T}}} \,{\bm{R}} \, ({{\bm{U}}_{i + 1}^j}) \approx {{\bm{\Phi }}^{\text{T}}} \, {\bm{R}} ( {{\bm{U}}_i^j}) + {{\bm{\Phi }}^{\text{T}}} \, {\bm{K}} ( {{\bm{U}}_i^j}) \, \bm{\Phi} \, \Delta {{\bm{U}}_{{\text{red}}}}_{i+1}^j = \bm{0} , \hfill \\
			{\Delta {{\bm{U}}_{{\text{red}}}}_{i+1}^j} =  - {\bm{K}_{\rm{red}}^{-1}} \, {\bm{R}}_{\rm{red}} \, ({{\bm{U}}_{i }^j}), \hfill \\
			{\bm{U}}_{i + 1}^j = {\bm{U}}_i^j + \bm{\Phi} \, \Delta {{\bm{U}}_{{\text{red}}}}_{i+1}^j, \hfill \\
			\left\| {{\bm{R}_{\rm{red}}} ( {{\bm{U}}_{i + 1}^j})} \right\| \leqslant tol,  \hfill \\
			i \leftarrow i + 1 . \hfill \\ 
		\end{gathered}
		\label{eq47}
	\end{equation}
	
	\noindent The primary computational burden associated with the Newton-Raphson method in the context of the full-order model comes from the inverse of the high-dimensional stiffness matrix, denoted as $\bm{K}^{-1}$. By contrast, utilizing the inverse of the reduced stiffness matrix $\bm{K}^{-1}_{\text{red}}$ leads to a noticeable reduction in computational complexity due to its lower dimension. Once the reduced residual $\bm{R}_{\text{red}}$ reaches a predefined tolerance ($tol$), the converged increment in the reduced-order solution $\Delta \bm{U}_{\text{red}}$ is deemed found. Subsequently, the reduced solution must be projected back to the full-order system in order to calculate solution vector $\bm{U}$ for the converged state:
	\begin{equation}
		\Delta \bm{U}_{\rm{full}} \approx \Delta \bm{U}^{j}_{i+1} = \bm{\Phi} \, \Delta {\bm{U}_{\rm{red}}}^{j}_{i+1}.
		\label{eq48}
	\end{equation}
	\noindent It is important to note that whenever the reduced residual converges or does not converge, the project and backward project are employed throughout the Newton iteration. 
	
	To better compare the accuracy and performance of the reduced-order simulations with the full-order model simulations, the relative error in displacement $\varepsilon _u$ and the CPU time ratio $\tau$ are defined and utilized in this study as follows:  
	\begin{equation}
		\varepsilon _{u}  = \frac{1}{N}\sum\limits_{j = 1}^N {\frac{{| {u_{p, \, \rm{red}}^j - u_{p, \, \rm{full}}^j} |}}{{| {u_{p, \, \rm{full}}^j} |}}}, \quad 
		\tau  = \sum\limits_{j = 1}^N {\frac{{T_{\rm{red}}^j}}{{{T^{j}_{\rm{full}}}}}}.
		\label{eq49}
	\end{equation}
	
	\noindent Here, $N$ represents the actual number of time steps. $u^{j}_{p, \, \rm{red}}$ and $u^{j}_{p,\, \rm{full}}$ denote the displacements obtained using the reduced and full-order model at a selected nodal point $p$ at time $j$. $T^{j}_{\rm{red}}$ and $T^{j}_{\rm{full}}$ are the solution time costs of the reduced and the full-order model at time $j$, respectively. The smaller the CPU time ratio $\tau$, the better the performance of the reduced-order model becomes. In the context of non-local damage scenarios, it is essential to recognize that a relative error measure for non-local damage variables cannot be used. To avoid a division by zero for points in the structure that remain undamaged, the absolute damage error $\varepsilon_{\bar{D}}$ is employed to evaluate the accuracy of the non-local damage variable in this study. 
	\begin{equation}
		\varepsilon _{\bar{D}}  = \frac{1}{N}\sum\limits_{j = 1}^N {{| {\bar{D}_{p,\, \rm{red}}^j - \bar{D}_{p, \, \rm{full}}^j} |}}
		\label{eq50}
	\end{equation}
	Here, $\bar{D}^{j}_{p, \, \rm{red}}$ and $\bar{D}^{j}_{p, \, \rm{full}}$ represent the non-local damage variable obtained using the reduced-order and full-order model at a selected nodal point $p$ at time $j$, respectively.
	
	\section{Numerical examples}
	\label{Model validation}
	In the following, two three-dimensional numerical examples are considered to assess the performance and properties of the newly developed multi-field decomposed MOR approach. For this, the latter is implemented into \textit{Python}. The material model and element are implemented by using the academic finite element software \textit{FEAP} \cite{taylor2014feap}. The visualization of all the contour plots is realized by the open-source software \textit{ParaView} \cite{ahrens200536}.
	
	\subsection{Test 1: Asymmetrically notched specimen} \label{ANS}
	
	\begin{figure}[!ht]
		\centering
		\includegraphics[width=9cm]{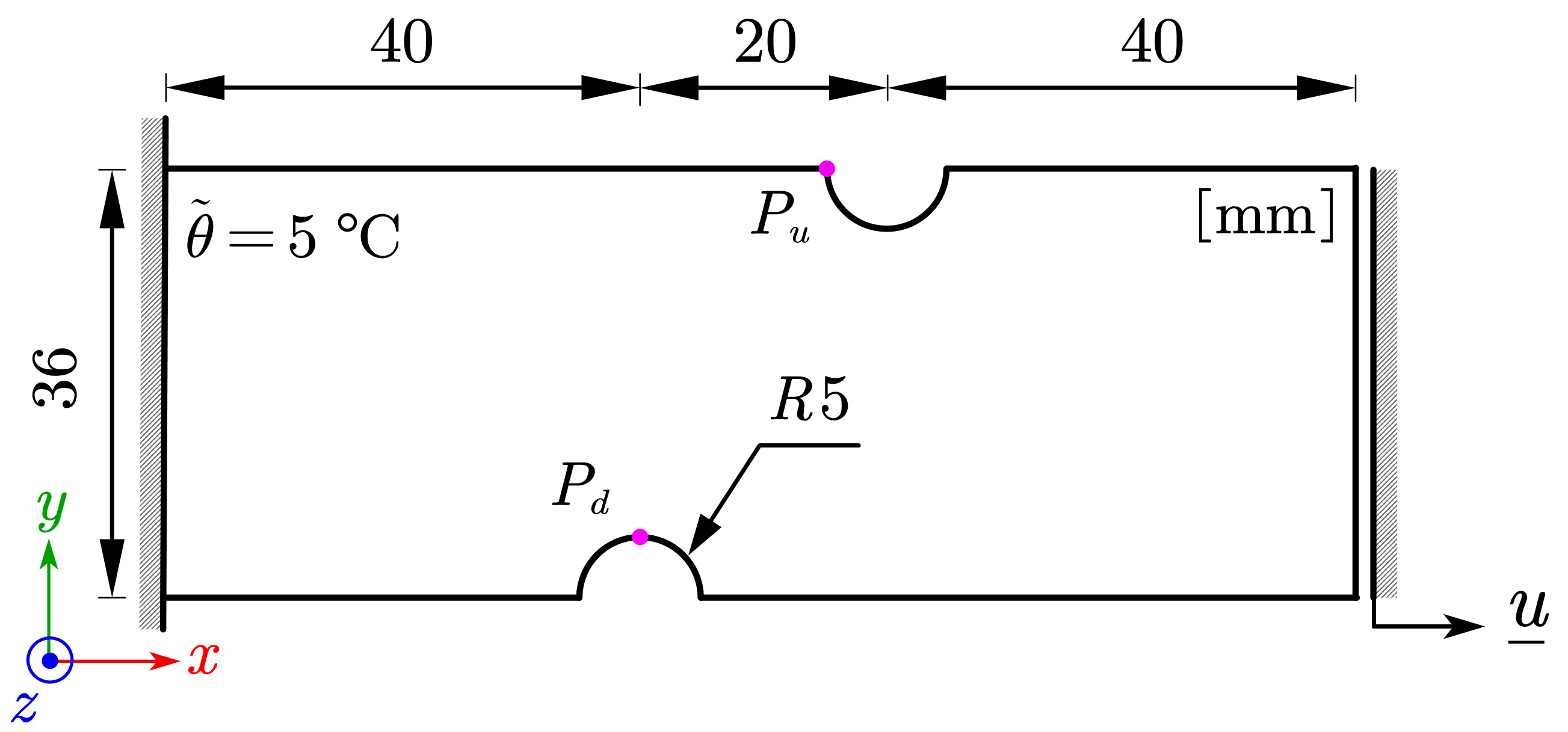}
		\caption{Schematic illustration of the geometry and boundary value problem of Test 1. The left edge is completely fixed and a horizontal displacement of $\underline{u}=3 \, \rm{mm}$ is prescribed on the right edge. The thickness of the specimen is 1 $\rm{mm}$. A constant temperature of $\tilde{\theta} = \SI{5}{\celsius}$ is assumed. Points $P_u$ and $P_d$ are selected to compare the errors of the displacement and non-local damage variables, respectively. The thickness direction is unconstrained.}
		\label{fig2}
	\end{figure}
	
	In this section, an asymmetrically notched specimen (based on an example from \citet{ambati2016phase}) is chosen to investigate the performance of the multi-field decomposed MOR approach, see \hyperref[fig2]{Fig.3}. The thickness of the specimen is meshed with one element.  
	
	\begin{table}[!ht]
		\begin{tabular}{llllll}
			\hline
			\textbf{Symbol}                      & \textbf{Parameter name}                          & \textbf{Test 1} & \textbf{Test 2} & \textbf{Unit}                    & \textbf{Eq.}				\\ \hline
			\textbf{elastic parameters} &                               &                             &                   & \multicolumn{1}{c}{}    \\  \hline
			$\Lambda$                   & Lam\'e constant                 & 25000                       & 101160            & MPa		& \eqref{eq5} \\
			$\mu$                       & shear modulus                 & 55000                       & 73255             & MPa 		& \eqref{eq5} \\ \hline
			\textbf{plastic parameters} &                               &                             &                   &              &           \\ \hline
			$\sigma _{y0}$              & initial yield stress          & 100                         & 340               & MPa           &  \eqref{eq17}          \\
			$a$                       & kinematic hardening parameter    & 62.5                        & 0                 & MPa          &  \eqref{eq6}           \\
			$P$                       & linear isotropic hardening parameter   & 0                           & 520               & MPa           &  \eqref{eq6}         \\
			$e_{p}$                   & $1^{\rm{st}}$ nonlinear isotropic hardening parameter          & 125                         & 296               & MPa           &  \eqref{eq6}         \\
			$f_{p}$                   & $2^{\rm{nd}}$ nonlinear isotropic hardening parameter         & 5                           & 18.9              & -             &  \eqref{eq6}       \\ \hline
			\textbf{damage parameters}  &                               &                             &                   &               &          \\ \hline
			$Y_0$                       & damage threshold              & 2.5                         & 750               & MPa           & \eqref{eq18}         \\
			$e_{d}$                     & $1^{\rm{st}}$ damage hardening parameter      & 5.0                         & 100               & MPa           & \eqref{eq7}          \\
			$f_{d}$                     & $2^{\rm{nd}}$ damage hardening parameter       & 100                         & 100               & -             & \eqref{eq7}        \\
			$A$                         & internal length scale         & 75                          & 50                & MPa $\rm{mm^2}$    & \eqref{eq8}          \\
			$H$                         & penalty parameter             & $10^6$                      & $10^4$            & MPa           & \eqref{eq8}          \\ \hline
			\textbf{thermal parameters} &                               &                             &                   &               &           \\ \hline
			$c$                         & volumetric heat capacity      & 3.59                        & 3.59              & $\rm{mJ/(mm^3 \, K}$)     & \eqref{eq31}      \\
			$\alpha$                    & thermal expansion coefficient & 1.1                         & 1.1               & $\rm{10^{-5}/K}$        & \eqref{eq2}    \\
			$K_0$                       & heat conductivity             & 50.2                        & 50.2              & $\rm{mW/(mm \, K}$)     &  \eqref{eq15}         \\
			$\theta _0$                 & reference temperature         & 273.15                      & 273.15            & $\rm{K}$                 & \eqref{eq2}      \\
			$\omega$                    & thermal softening parameter   & 0.002                       & 0.002             & $\rm{1/K}$               & \eqref{eq16}      \\ \hline
		\end{tabular}
		\label{T1}
		\caption{Material parameters list for Test 1 and Test 2 with relevant equations.}
	\end{table}
	
	The specimen is meshed by 8-node hexahedral elements (see \hyperref[fig1]{Fig.2}), leading to discretizations with 2192, 4663, 6718, and 8979 elements. Each node possesses 5 DOFs (three displacement components, non-local damage, and temperature). The coordinates of selected displacement and damage points are defined as ${P_u} \, (55, \, 36, \, 0)$ and ${P_d} \, (40, \, 5, \, 1)$, respectively. The temperature of the specimen is prescribed with a constant temperature $\tilde{\theta} = \SI{5}{\celsius}$. More details are given in \hyperref[fig2]{Fig.3}. The mechanical and thermal parameters for the asymmetrically notched specimen are taken from \citet{brepols2020gradient} and \citet{dittmann2020phase}, respectively, see \hyperref[T1]{Table 1}. The objective is to investigate the reduction performance of the multi-field decomposed MOR for different meshes.
	
	\subsubsection{Comparison of singular values of different fields}
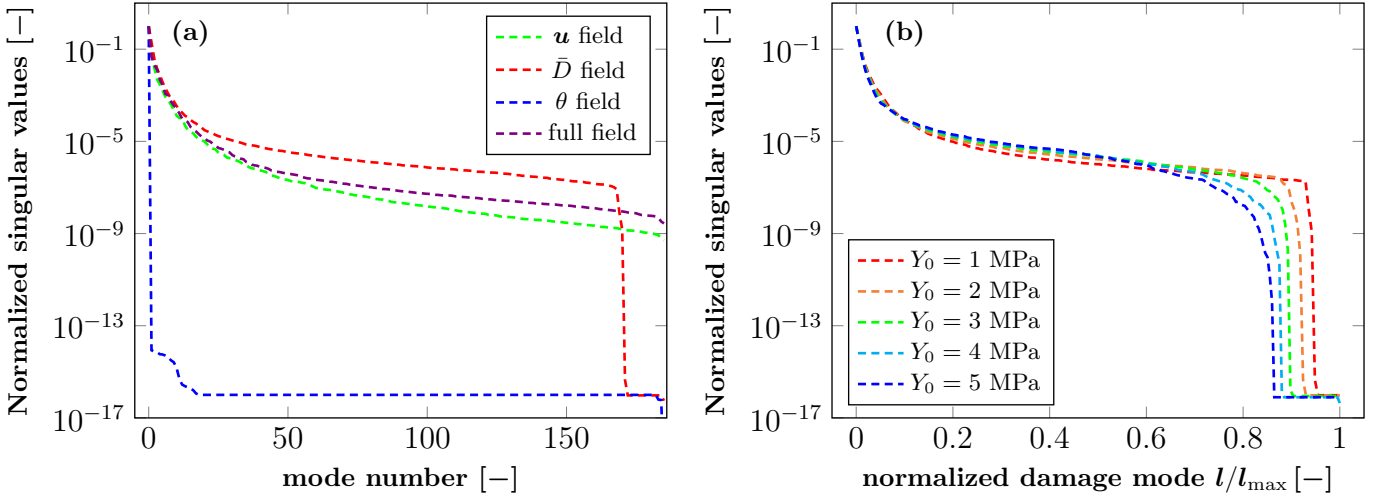
\begin{figure*}[!ht]
		\begin{minipage}[b]{0.5\linewidth}
			\begin{tikzpicture}
				\begin{axis}[name=plotA,
					xmin=-5,xmax=186,
					ymin=1e-17,ymax=10,
					scale only axis,
					ticklabel style={font=\large},
					xlabel={mode number $\bm{[-]}$}, 
					ylabel={Normalized singular values $\bm{[-]}$},
					line width=0.5pt,
					ymode=log,
					log basis y={10},
					axis line style={line width=0.5pt},
					font=\bfseries,
					legend pos=north east,
					legend style={fill=none, font=\small}
					]
					\addplot[
					line width=1.0pt,
					color=green, densely dashed,
					] table[x=modes,y=U]{./Asy/modes-SVD.txt};
					
					\addplot[
					line width=1.0pt,
					color=red, densely dashed,
					] table[x=modes,y=D]{./Asy/modes-SVD.txt};
					
					\addplot[
					line width=1.0pt,
					color=blue, densely dashed,
					] table[x=modes,y=T]{./Asy/modes-SVD.txt};
					
					\addplot[
					line width=1.0pt,
					color=violet, densely dashed,
					] table[x=modes,y=full]{./Asy/modes-SVD.txt};
					
					\legend{
						{\text{$\bm{u}$ field   }},
						{\text{$\bar{D}$ field   }},
						{\text{$\theta$ field   }},
						{\rm{full field}}
					}
					\node[align=center] at (15,0.5) {(a)};
				\end{axis}
			\end{tikzpicture}
		\end{minipage}
		\begin{minipage}[b]{0.5\linewidth}
			\centering
			\begin{tikzpicture}
				\begin{axis}[
					clip=true,
					xmin=-0.05,xmax=1.05,
					ymin=1e-17,ymax=10,
					scale only axis,
					ticklabel style={font=\large},
					xlabel={normalized damage mode $\bm{l}/\bm{l}_{\rm{max}} \, \bm{[-]}$}, 
					ylabel={Normalized singular values $\bm{[-]}$},
					ymode=log,
					log basis y={10},
					line width=0.5pt,
					axis line style={line width=0.5pt},
					font=\bfseries,
					legend pos=south west,
					legend style={fill=none, font=\small}
					]
					\addplot[
					line width=1.0pt,
					color=red, densely dashed] table[x=modes,y=sigma]{./Asy/modes-SVD-th1MPa.txt};
					
					\addplot[
					line width=1.0pt,
					color=Orange, densely dashed] table[x=modes,y=sigma]{./Asy/modes-SVD-th2MPa.txt};
					
					\addplot[
					line width=1.0pt,
					color=green, densely dashed] table[x=modes,y=sigma]{./Asy/modes-SVD-th3MPa.txt};
					
					\addplot[
					line width=1.0pt,
					color=cyan, densely dashed] table[x=modes,y=sigma]{./Asy/modes-SVD-th4MPa.txt};
					
					\addplot[
					line width=1.0pt,
					color=blue, densely dashed] table[x=modes,y=sigma]{./Asy/modes-SVD-th5MPa.txt};
					\legend{
						{$Y_{0}=1$ \rm{MPa}},
						{$Y_{0}=2$ \rm{MPa}},
						{$Y_{0}=3$ \rm{MPa}},
						{$Y_{0}=4$ \rm{MPa}},
						{$Y_{0}=5$ \rm{MPa}}
					}
					\node[align=center] at (0.1,0.5) {(b)};
				\end{axis}
			\end{tikzpicture}
		\end{minipage}
		\caption{(a) Normalized singular values ($\sigma /\sigma_{\rm{max}}$) in case of the multi-field decomposed MOR approach (separate singular values for the displacement field $\bm{u}$, non-local damage field $\bar{D}$, and temperature field $\theta$) and the classical non-decomposed MOR approach for a mesh with 8979 elements; (b) normalized singular values of the non-local damage field ($\sigma_{\bar{D}}/\sigma_{\bar{D}, \, \rm{max}}$) using different damage thresholds $Y_0$ in case of the multi-field decomposed MOR approach.}
		\label{fig3}
	\end{figure*}

	The normalized singular values of converged nodal snapshots are first investigated and compared by using the novel multi-field decomposed MOR approach and the classical non-decomposed MOR approaches, see \hyperref[fig3]{Fig.4 (a)}. The results of \hyperref[fig3]{Fig.4} are based on the offline pre-computations using the full-order model. It can be observed from \hyperref[fig3]{Fig.4} that the normalized singular values monotonically decrease as the mode number increases. This phenomenon can be attributed to the automatic ranking of singular values by the singular value decomposition. Interestingly enough, one can observe that the singular value of the temperature field quickly reaches approximately zero after considering only a few modes. This can be attributed to the constant temperature in the specimen. Consequently, only a few modes are required to accurately represent the temperature field. However, the micromorphic damage variable $\bar{D}$ significantly influences the mechanical response throughout the loading process of the asymmetrically notched specimen, see \hyperref[fig3]{Fig.4 (a)}. Only if the number of damage modes reaches $l=183$, the normalized $\sigma$ decreases to $8.5 \times 10^{-17}$. It is also quite instructive to observe that the curve of the normalized singular values of the displacement field $\bm{u}$ (using the multi-field decomposed MOR approach) and the curve of the singular values obtained using the standard non-decomposed MOR approach resemble each other quite well. This indicates that the displacements in this example mainly dominate the size of the singular values.
	
	Therefore, it is not surprising that a higher number of modes is required for the displacement and damage fields than for the temperature field when truncating the projection matrix. The main benefit of \hyperref[fig3]{Fig.4 (a)} is that it can be used to concretely determine how many modes $k$ for the displacement field, $l$ for the non-local damage field, and $m$ for the temperature field are needed. \hyperref[fig3]{Fig.4 (b)} indicates that the number of required damage modes decreases as the damage threshold increases. 
	
	\subsubsection{Influence of the number of displacement and damage modes}
	The difficulty of MOR lies in balancing the order of reduction with simulation accuracy. The reduction order is determined by the total number of modes $k+l+m$, as it controls the size of the reduced stiffness matrix, see \hyperref[eq45]{Eq.\eqref{eq45}}. Generally, an increase in the number of modes incorporated into the reduced-order model enriches the informational content of the projection matrix, thereby enhancing the accuracy of the reduced system. However, this also leads to higher computational costs. Therefore, it is essential to select the appropriate mode numbers for displacement, non-local damage, and temperature fields in a reasonable manner. Specifically, the number of modes for each field is chosen by considering the following tolerance value $\sigma_{\rm{tol}}$:
	\begin{equation}
		\frac{\sigma _j}{\sigma _{\rm{max}}} \le  \sigma_{\rm{tol}} = 10^{-7}.
		\label{eq51}
	\end{equation}
	The index $j$ represents the first mode number where ${\sigma _{j}}/\sigma_{\text{max}} \le 10^{-7}$. According to the normalized singular value development of \hyperref[fig3]{Fig.4 (a)}, the reference values for the number of modes for displacement field, non-local damage field, and temperature field are given as $k=59$, $l=167$, and $m=0$, respectively. However, \hyperref[eq51]{Eq.\eqref{eq51}} provides the considered mode numbers for various fields, it does not imply that the MOR model can reproduce the strong nonlinear damage evolution by employing these mode numbers. Furthermore, when employing tolerant mode $m=0$ in temperature, it is impossible to calculate the temperature development in the temperature field. With the consideration of the reference modes number over three fields and higher accuracy results, the numbers of modes for the displacement and temperature are finally chosen as $k=120$ and $m=5$, respectively.
	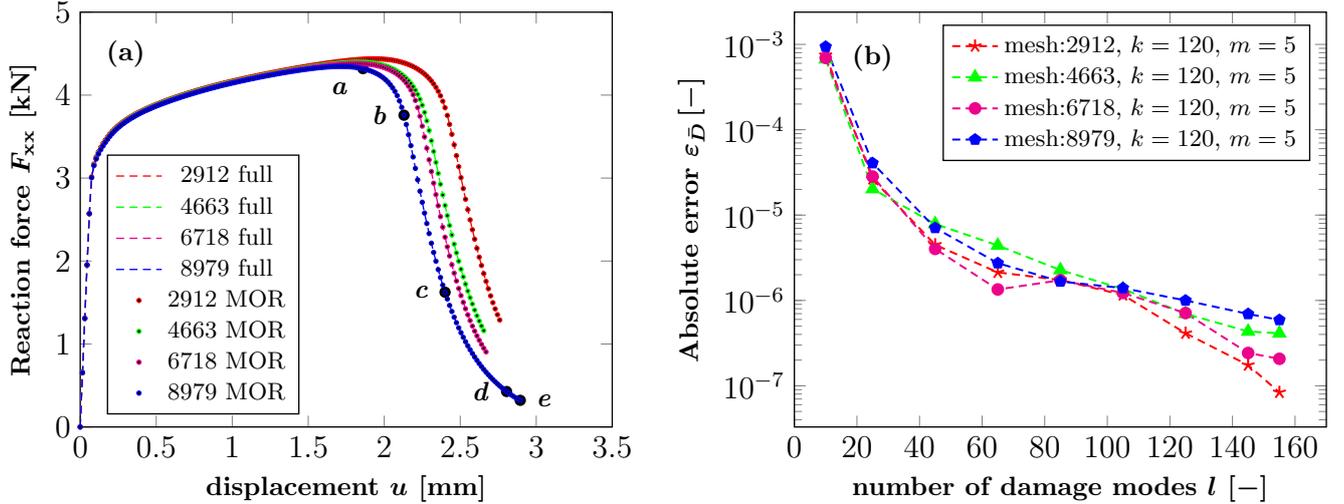
\begin{figure*}[!ht]
		\begin{minipage}[b]{0.5\linewidth}
			\centering
			\begin{tikzpicture}
				\begin{axis}[
					name=plotA,
					clip=true,
					xmin=0,xmax=3.5,
					ymin=0,ymax=5,
					scale only axis,
					log ticks with fixed point,
					ticklabel style={font=\large},
					xlabel={displacement $\bm{u}$ $\textbf{[mm]}$}, 
					ylabel={Reaction force $\bm{F_{\rm{xx}}}$ \textbf{[kN]}},
					line width=0.5pt,
					font=\bfseries,
					legend pos=south west,
					legend style={fill=none, font=\small, at={(0.05, 0.03)}},
					scatter/classes={%
						a={mark=*,draw=red},
						b={mark=*,draw=green},
						c={mark=*,draw=magenta},
						d={mark=*,draw=blue}
					}
					]
					\addplot[
					color=red, densely dashed,
					] table[x=u2912pre,y=f2912pre]{./Asy/FU2912pre.txt};
					
					\addplot[
					color=green, densely dashed,
					] table[x=u4663pre,y=f4663pre]{./Asy/FU4663pre.txt};
					
					\addplot[
					color=magenta, densely dashed,
					] table[x=u6718pre,y=f6718pre]{./Asy/FU6718pre.txt};
					
					\addplot[
					color=blue, densely dashed,
					] table[x=u8979pre,y=f8979pre]{./Asy/FU8979pre.txt};
					\addplot[
					scatter, only marks,
					mark=*,mark options={solid},
					scatter src=explicit symbolic,
					mark size=0.75pt,
					] table[x=u2912mor,y=f2912mor,meta=label]{./Asy/FU2912mor.txt};
					
					\addplot[
					scatter, only marks,
					mark=*,mark options={solid},
					scatter src=explicit symbolic,
					mark size=0.75pt,
					] table[x=u4663mor,y=f4663mor,meta=label]{./Asy/FU4663mor.txt};
					
					\addplot[
					scatter, only marks,
					mark=*,mark options={solid},
					scatter src=explicit symbolic,
					mark size=0.75pt,
					] table[x=u6718mor,y=f6718mor,meta=label]{./Asy/FU6718mor.txt};
					
					\addplot[
					scatter, only marks,
					mark=*,mark options={solid},
					scatter src=explicit symbolic,
					mark size=0.75pt,
					] table[x=u8979mor,y=f8979mor,meta=label]{./Asy/FU8979mor.txt};
					\legend{\rm{2912 full},\rm{4663 full},\rm{6718 full},\rm{8979 full},\rm{2912 MOR},\rm{4663 MOR},\rm{6718 MOR},\rm{8979 MOR}}
					\node[align=center] at (0.3,4.5) {(a)};
					\node[circle, fill, inner sep=1.5pt, label=below left:\textit{a}] (a) at (1.860000E+00 , 4.318769E+00) {};
					\node[circle, fill, inner sep=1.5pt, label=left:\textit{b}] (b) at (2.130000E+00 , 3.757736E+00) {};
					\node[circle, fill, inner sep=1.5pt, label=left:\textit{c}] (c) at (2.400000E+00 , 1.624238E+00) {};
					\node[circle, fill, inner sep=1.5pt, label=left:\textit{d}] (d) at (2.805000E+00, 4.271763E-01) {};
					\node[circle, fill, inner sep=1.5pt, label=right:\textit{e}] (e) at (2.895000E+00, 3.210785E-01) {};
				\end{axis}
			\end{tikzpicture}
		\end{minipage}
		\begin{minipage}[b]{0.5\linewidth}
			\begin{tikzpicture}
				\begin{axis}[
					xmin=0,xmax=170,
					scale only axis,
					ticklabel style={font=\large},
					xlabel={number of damage modes $\bm{l}$ $\bm{[-]}$}, 
					ylabel={Absolute error $\bm{\varepsilon _{\bar{D}}} \, \bm{[-]}$},
					line width=0.5pt,
					ymode=log,
					log basis y={10},
					font=\bfseries,
					legend pos=north east,
					legend style={fill=none, font=\small}
					]
					\addplot[
					line width=0.75pt,
					color=red, densely dashed,
					mark=star,
					mark size=2.5pt,
					mark options={solid},
					] table[x=modesD,y=D2912]{./Asy/DErrors.txt};
					
					\addplot[
					line width=0.75pt,
					color=green, densely dashed,
					mark=triangle*,
					mark size=2.5pt,
					mark options={solid},
					] table[x=modesD,y=D4663]{./Asy/DErrors.txt};
					
					\addplot[
					line width=0.75pt,
					color=magenta, densely dashed,
					mark=*,mark size=2pt,
					mark options={solid},
					] table[x=modesD,y=D6718]{./Asy/DErrors.txt};
					
					\addplot[
					line width=0.75pt,
					color=blue, densely dashed,
					mark=pentagon*,
					mark size=2pt,
					mark options={solid},
					] table[x=modesD,y=D8979]{./Asy/DErrors.txt};
					
					\legend{
						{\rm{mesh:2912}, $k=120$, $m=5$},
						{\rm{mesh:4663}, $k=120$, $m=5$},
						{\rm{mesh:6718}, $k=120$, $m=5$},
						{\rm{mesh:8979}, $k=120$, $m=5$}
					}
					\node[align=center] at (25,7e-4) {(b)};
				\end{axis}
			\end{tikzpicture}
		\end{minipage}
		\caption{(a) Force-displacement curve comparison between the full-order and the reduced-order models. The employed mode numbers of the reduced-order model are selected as $k=120$, $l=25$, and $m=5$ for the displacement field, non-local damage field, and temperature field, respectively; points $\bm{a}$ to $\bm{e}$ represent the states at $u=1.860, \, 2.130, \, 2.400, \, 2.850, \, 2.895$ mm, see \hyperref[fig5-1]{Fig.6}; (b) absolute errors in non-local damage for reduced simulations with different numbers of damage modes $l \in \left[10, \, 155\right]$.}
		\label{fig4}
	\end{figure*}
	
	Next, a comparison between the full-order model and reduced-order model is conducted by varying the number of elements in the mesh from 2912 to 8979, while keeping the number of modes fixed ($k=120$, $l=25$, and $m=5$), see \hyperref[fig4]{Fig.5 (a)}. It can be observed from \hyperref[fig4]{Fig.5 (a)} that the force-displacement curves of the reduced-order model accurately match each of the full-order models, even in regions where severe damage-induced softening occurs.  The error of the reduced-order solution lies between $2.022 \times 10^{-5} \sim 2.063 \times 10^{-5}$. In the second investigation, the number of non-local damage modes is varied from $l=10$ to $l=155$ to investigate the absolute error $\varepsilon _{\bar{D}}$, see \hyperref[fig4]{Fig.5 (b)}. It can be observed that $\varepsilon _{\bar{D}}$ decreases as the number of non-local damage modes increases. The errors lie between $2.071 \times 10^{-7} \sim 9.393 \times 10^{-4}$. The study's notable accuracy stems from $\mathbf{SVD}$, but also from the novel multi-field decomposed approach in which the individual projections are applied for the different fields.
	
	\begin{figure}[!ht]
		\centering
		\includegraphics[width=14cm]{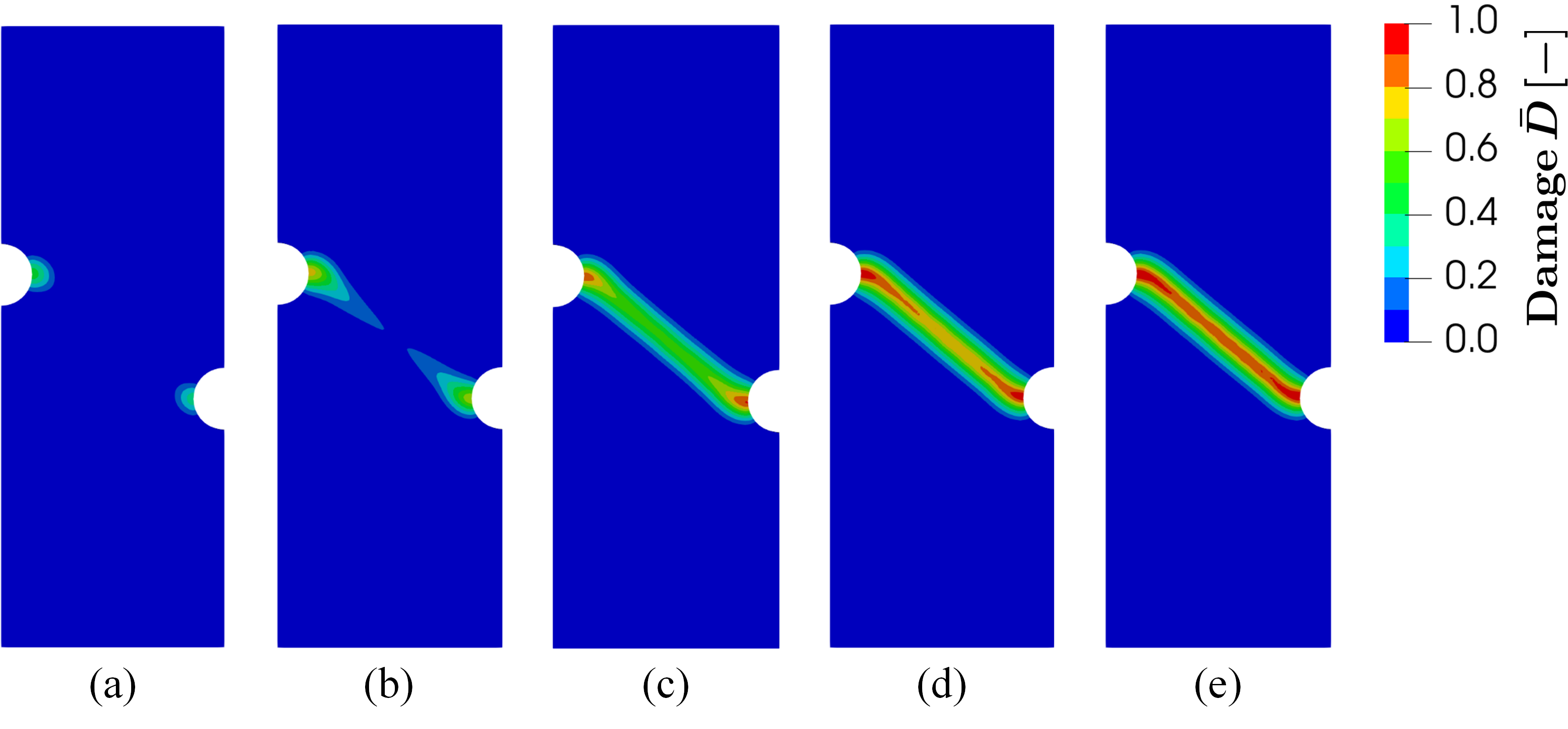}
		\caption{Non-local damage contour plots of a simulation carried out using the reduced-order model ($k=120$, $l=25$, and $m=5$) for the states $\bm{a}$ to $\bm{e}$ (see \hyperref[fig4]{Fig.5 (a)}) with corresponding displacements $u=1.860, \, 2.130, \, 2.400, \, 2.850, \, 2.895$ mm. The specimen is meshed by 8979 elements.} 
		\label{fig5-1}
	\end{figure}
	
	In summary, the multi-field decomposed MOR method effectively overcomes challenges in the damage-induced material softening region, ensuring high accuracy in multi-physical damage simulations with minimal errors. For detailed non-local damage contour plots obtained by using the reduced-order model and a comparison with results obtained by using the full-order model, please refer to \hyperref[fig5-1]{Fig.6} and \hyperref[fig9-1]{Fig.12 (a)}, respectively. 
	
	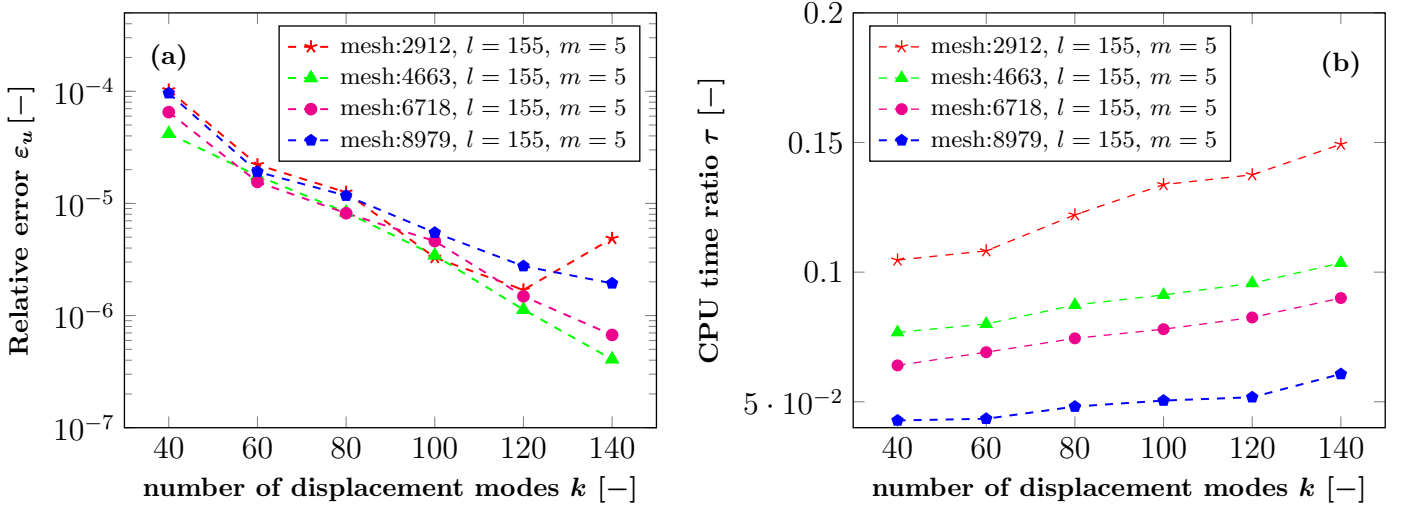
\begin{figure*}[!ht]
		\begin{minipage}[b]{0.5\linewidth}
			\begin{tikzpicture}
				\begin{axis}[
					xmin=30,xmax=150,
					ymin=1e-7,ymax=5e-4,
					scale only axis,
					ticklabel style={font=\large},
					xlabel={number of displacement modes $\bm{k}$ $\bm{[-]}$}, 
					ylabel={Relative error $\bm{\varepsilon _{u}} \, \bm{[-]}$},
					line width=0.5pt,
					ymode=log,
					log basis y={10},
					font=\bfseries,
					legend pos=north east,
					legend style={fill=none, font=\small}
					]
					\addplot[
					line width=0.75pt,
					color=red, dashed,
					mark=star,
					mark size=2.5pt,
					mark options={solid},
					] table[x=modesU,y=U2912]{./Asy/UErrors.txt};
					
					\addplot[
					line width=0.75pt,
					color=green, dashed,
					mark=triangle*,
					mark size=2.5pt,
					mark options={solid},
					] table[x=modesU,y=U4663]{./Asy/UErrors.txt};
					
					\addplot[
					line width=0.75pt,
					color=magenta, dashed,
					mark=*,mark size=2pt,
					mark options={solid},
					] table[x=modesU,y=U6718]{./Asy/UErrors.txt};
					
					\addplot[
					line width=0.75pt,
					color=blue, dashed,
					mark=pentagon*,
					mark size=2pt,
					mark options={solid},
					] table[x=modesU,y=U8979]{./Asy/UErrors.txt};
					
					\legend{
						{\rm{mesh:2912, $l=155$, $m=5$}},
						{\rm{mesh:4663, $l=155$, $m=5$}},
						{\rm{mesh:6718, $l=155$, $m=5$}},
						{\rm{mesh:8979, $l=155$, $m=5$}}
					}
					\node[align=center] at (40,2e-4) {(a)};
				\end{axis}
			\end{tikzpicture}
		\end{minipage}
		\begin{minipage}[b]{0.5\linewidth}
			\centering
			\begin{tikzpicture}
				\begin{axis}[
					clip=true,
					xmin=30,xmax=150,
					ymin=0.04,ymax=0.2,
					scale only axis,
					ticklabel style={font=\large},
					xlabel={number of displacement modes $\bm{k}$ $\bm{[-]}$}, 
					ylabel={CPU time ratio $\bm{\tau}$ $\bm{[-]}$},
					line width=0.5pt,
					font=\bfseries,
					legend pos=north west,
					legend style={fill=none, font=\small}
					]
					\addplot[
					color=red, dashed,
					mark=star,
					mark size=2.5pt,
					mark options={solid},
					] table[x=modesT,y=T2912]{./Asy/TCPU.txt};
					
					\addplot[
					color=green, dashed,
					mark=triangle*,
					mark size=2.5pt,
					mark options={solid},
					] table[x=modesT,y=T4663]{./Asy/TCPU.txt};
					
					\addplot[
					color=magenta, dashed,
					mark=*,mark size=2pt,
					mark options={solid},
					] table[x=modesT,y=T6718]{./Asy/TCPU.txt};
					
					\addplot[
					line width=0.75pt,
					color=blue, dashed,
					mark=pentagon*,
					mark size=2pt,
					mark options={solid},
					] table[x=modesT,y=T8979]{./Asy/TCPU.txt};
					
					\legend{
						{\rm{mesh:2912, $l=155$, $m=5$}},
						{\rm{mesh:4663, $l=155$, $m=5$}},
						{\rm{mesh:6718, $l=155$, $m=5$}},
						{\rm{mesh:8979, $l=155$, $m=5$}}
					}
					\node[align=center] at (140,0.18) {(b)};
				\end{axis}
			\end{tikzpicture}
		\end{minipage}
		\caption{(a) Relative errors of the displacement for different numbers of displacement modes $k$; (b) CPU time ratio comparison for different numbers of displacement modes $k$.}
		\label{fig5}
	\end{figure*}
	
	After investigating the influence of the number of damage modes, the influence of the number of displacement modes is investigated by comparing the relative error $\varepsilon _{u}$ of the displacement $u$. As illustrated in \hyperref[fig5]{Fig.7 (a)}, the general trend indicates a decrease in the relative error as the number of displacement modes increases. The relative error $\varepsilon _u$ lies between $4.077 \times 10^{-7} \sim 9.620 \times 10^{-5}$. A higher number of displacement modes is not considered here, because for $k>140$, $\varepsilon_{u}$ starts to fluctuate, indicating an instability that is induced by noise from unnecessary modes. Furthermore, as the mode number exceeds a certain threshold, $\sigma$ approaches zero. This phenomenon leads to multiple columns of the projection matrix $\bm{\Phi}$ becoming nearly zero, owing to the inclusion of superfluous modes. Consequently, the reduced stiffness matrix $\bm{K}_{\rm{red}}$ (\hyperref[eq45]{Eq.\eqref{eq45}}) could become singular during the Newton-Raphson iteration (\hyperref[eq47]{Eq.\eqref{eq47}}). In other words, selecting the mode number appropriately requires a delicate balance between desired accuracy and reduction order, which can be achieved by referring to the decomposed singular values from \hyperref[fig3]{Fig.4}.	
	
	Another important criterion for evaluating the performance of the novel MOR approach is the CPU time ratio $\tau$. As shown in \hyperref[fig5]{Fig.7 (b)}, the CPU time ratio increases as the number of modes increases. This is clear because an increased number of modes results in a higher-dimensional reduced-order problem, subsequently leading to higher computational costs. When comparing the results of different meshes, a constant number of modes will result in a constant $T_{\rm{red}}$, while $T_{\rm{full}}$ is approximately linear with the number of elements. Therefore, the results of \hyperref[fig5]{Fig.7 (b)} are understandable in the sense that  $\tau _{\rm{8979}}<\tau _{\rm{6718}}<\tau _{\rm{4663}}<\tau _{\rm{2912}}$ for each selected number of displacement modes $k$. The significant reduction in solution time for the mesh with 8979 elements is $95.72\%$ ($\tau _{8979} = 4.28\times10^{-2}$) and the corresponding relative error $\varepsilon _{u} $ reaches $9.62\times10^{-5}$. Compared to the full-order simulation, the reduced solution time for simulation is 23.37\footnote{The reduced performance is calculated by $\frac{1}{\tau}$.} times faster.

	\subsection{Test 2: Flat I-shaped specimen}
	
	\begin{figure*}[!ht]
		\centering
		\includegraphics[width=16cm]{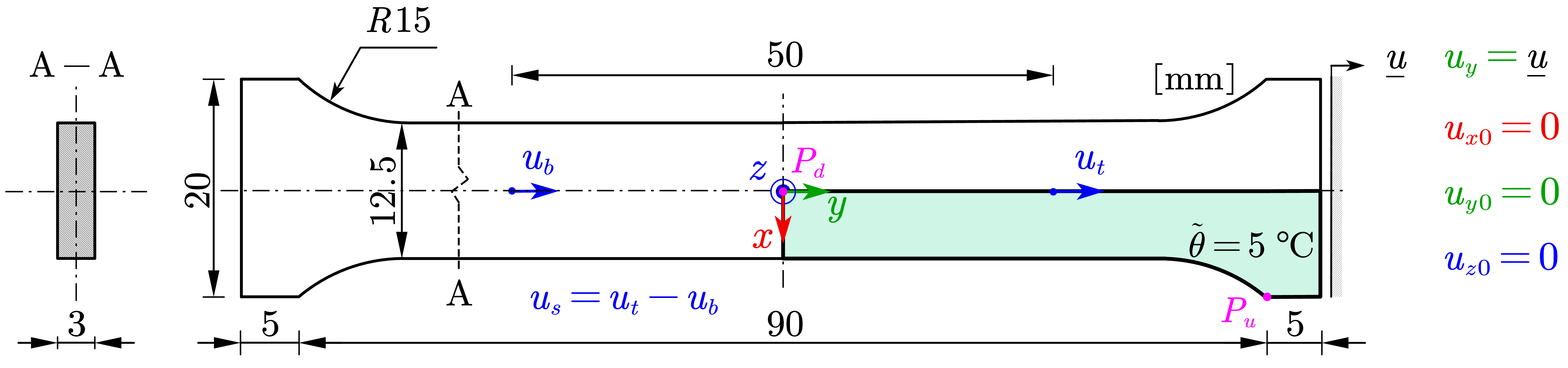}
		\caption{Schematic illustration of the I-shaped specimen geometry and the considered boundary value problem. The left picture shows the cross-section A-A. The displacement on the right edge is prescribed as $\underline{u}=9 $ mm. The symmetric part of the structure, highlighted in green, comprises one-eighth of the geometry utilized in the simulations. The specimen's thickness is 3 mm. A constant temperature $\tilde{\theta} = \SI{5}{\celsius}$ is prescribed.
		}
		\label{fig6}
	\end{figure*}
	
	Test 2 is concerned with an I-shaped specimen with isothermal expansion, see \hyperref[fig6]{Fig.8}. On its right surface, a displacement in the $y$ direction of $\underline{u}=9 \, mm$ is prescribed. Due to symmetry, only one-eighth of the structure highlighted in green needs to be considered in the simulations. The mechanical and thermal parameters are taken from \citet{felder2022thermo} and \citet{dittmann2020phase}, respectively, see \hyperref[T1]{Table 1}. A constant temperature of $\tilde{\theta} = \SI{5}{\celsius} $ is assumed in order to stay consistent with \hyperref[ANS]{Section. \textit{4.1}}. It is noteworthy that the value of $u_s$ is set to $\pm 25$ mm from the midpoint (see \hyperref[fig6]{Fig. 8}) in concordance with the prior research conducted by \citet{ambati2016phase}. To compare the reduced-order model with the full-order model, points ${P_u} \, (10, \, 45, \, 1.5)$ and ${P_d} \, (0, \, 0, \, 0)$ are selected to measure the nodal values of $u_y$ and $\bar{D}$, respectively.
	
	\subsubsection{Comparison of singular values of different fields}
	
	\begin{figure*}[!ht]
		\begin{minipage}[b]{0.5\linewidth}
			\begin{tikzpicture}
				\begin{axis}[
					xmin=-5,xmax=185,
					ymin=1e-17,ymax=10,
					scale only axis,
					ticklabel style={font=\large},
					xlabel={mode number $\bm{[-]}$}, 
					ylabel={Normalized singular values $\bm{[-]}$},
					line width=0.5pt,
					ymode=log,
					log basis y={10},
					font=\bfseries,
					legend pos=north east,
					legend style={fill=none, font=\small}
					]
					\addplot[
					line width=1.0pt,
					color=green, densely dashed,
					] table[x=modes,y=U]{./I-shape/modes-SVD.txt};
					
					\addplot[
					line width=1.0pt,
					color=red, densely dashed,
					] table[x=modes,y=D]{./I-shape/modes-SVD.txt};
					
					\addplot[
					line width=1.0pt,
					color=blue, densely dashed,
					] table[x=modes,y=T]{./I-shape/modes-SVD.txt};
					
					\addplot[
					line width=1.0pt,
					color=violet, densely dashed,
					] table[x=modes,y=full]{./I-shape/modes-SVD.txt};
					
					\legend{
						{\text{$\bm{u}$ field  }},
						{\text{$\bar{D}$ field  }},
						{\text{$\theta$ field  }},
						{\text{full field}}
					}
					\node[align=center] at (15,0.5) {(a)};
				\end{axis}
			\end{tikzpicture}
		\end{minipage}
		\begin{minipage}[b]{0.5\linewidth}
			\centering
			\begin{tikzpicture}
				\begin{axis}[
					clip=true,
					xmin=-0.05,xmax=1.05,
					ymin=1e-17,ymax=10,
					scale only axis,
					ticklabel style={font=\large},
					xlabel={normalized damage mode $\bm{l}/\bm{l}_{\rm{max}} \, \bm{[-]}$}, 
					ylabel={Normalized singular values $\bm{[-]}$},
					ymode=log,
					log basis y={10},
					line width=0.5pt,
					font=\bfseries,
					legend pos=north east,
					legend style={fill=none, font=\small}
					]
					\addplot[
					line width=1.0pt,
					color=red, densely dashed] table[x=modes,y=sigma]{./I-shape/modes-SVD-580.txt};
					
					\addplot[
					line width=1.0pt,
					color=Orange, densely dashed] table[x=modes,y=sigma]{./I-shape/modes-SVD-4113.txt};
					
					\addplot[
					line width=1.0pt,
					color=green, densely dashed] table[x=modes,y=sigma]{./I-shape/modes-SVD-8552.txt};
					
					\addplot[
					line width=1.0pt,
					color=cyan, densely dashed] table[x=modes,y=sigma]{./I-shape/modes-SVD-13660.txt};
					
					\addplot[
					line width=1.0pt,
					color=blue, densely dashed] table[x=modes,y=sigma]{./I-shape/modes-SVD-18510.txt};
					\legend{
						{$\text{mesh: 580    }$},
						{$\text{mesh: 4113  }$},
						{$\text{mesh: 8552  }$},
						{$\text{mesh: 13660}$},
						{$\text{mesh: 18510}$}
					}
					\node[align=center] at (0.1,0.5) {(b)};
				\end{axis}
			\end{tikzpicture}
		\end{minipage}
		\caption{(a) Normalized singular values ($\sigma / \sigma_{\text{max}}$) in case of the multi-field decomposed MOR approach and classical non-decomposed MOR approach for a mesh with 4113 elements; (b) comparison of normalized singular values of the non-local damage field ($\sigma_{\bar{D}}/\sigma_{\bar{D}, \, \rm{max}}$) for different meshes ranging from $580$ to $18510$ elements.}
		\label{fig7}
	\end{figure*}
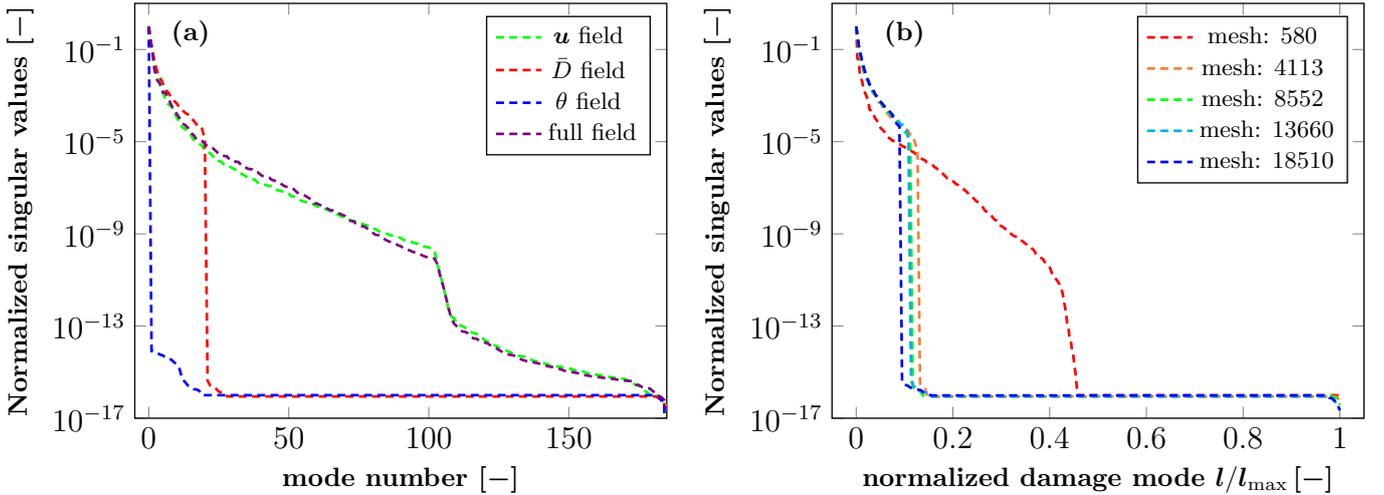
	
	This section compares the normalized singular values obtained using the novel multi-field decomposed MOR and standard MOR with each other. Afterwards, the normalized singular values of the non-local damage field are investigated for different meshes with elements ranging from $580$ to $18510$. It is interesting to see that the normalized singular values of the non-local damage approach zero as the number of damage modes $l$ exceeds 30, see \hyperref[fig7]{Fig.9 (a)}. Compared to the asymmetrically notched specimen (Test 1) in \hyperref[ANS]{Section. \textit{4.1}}, the influence of the damage field in the case of the I-shaped geometry (Test 2) is smaller. This is because the damage threshold chosen in Test 2 ($Y_{0}=750$ MPa) is much higher than in the case of Test 1 ($Y_{0}=2.5$ MPa), as listed in \hyperref[T1]{Table 1}. Utilizing such a high damage threshold indicates that damage initialization will occur only when the conjugate driving force $Y$ reaches a very high level, as described in \hyperref[eq18]{Eq.\eqref{eq18}}.  Therefore, the damage-induced softening behavior becomes significant only later in the process and makes a minor contribution to the overall simulation. Moreover, the singular value development of the displacement field closely aligns with the singular development obtained using the classical MOR approach, as depicted in \hyperref[fig7]{Fig.9 (a)}. In other words, the material behavior is dominated by plasticity in Test 2. To achieve mesh convergence in the damage simulations, \hyperref[fig7]{Fig.9 (b)} illustrates that the influence of the element numbers on the singular value development becomes negligible once the element number exceeds 4113. In summary, \hyperref[fig7]{Fig.9 (b)} provides the minimum number of damage modes required to obtain stable errors for various mesh discretizations when applying the MOR approach. 
	
	\subsubsection{Influence of the number of displacement and damage modes}
\begin{figure*}[!ht]
		\begin{minipage}[b]{0.5\linewidth}
			\centering
			\begin{tikzpicture}
				\begin{axis}[
					clip=true,
					xmin=0,xmax=18,
					ymin=0,ymax=25,
					scale only axis,
					log ticks with fixed point,
					ticklabel style={font=\large},
					xlabel={displacement $\bm{u_s}$ \textbf{[mm]}}, 
					ylabel={Reaction force $\bm{F_{\rm{yy}}}$ \textbf{[kN]}},
					line width=0.5pt,
					font=\bfseries,
					legend pos=south west,
					legend style={fill=none, font=\small, at={(0.08, 0.03)}},
					scatter/classes={%
						a={mark=*,draw=red},
						b={mark=*,draw=green},
						c={mark=*,draw=magenta},
						d={mark=*,draw=blue}
					}
					]
					\addplot[
					color=red, densely dashed,
					] table[x=u4113pre,y=f4113pre]{./I-shape/FU4113pre.txt};
					
					\addplot[
					color=green, densely dashed,
					] table[x=u8552pre,y=f8552pre]{./I-shape/FU8552pre.txt};
					
					\addplot[
					color=purple, densely dashed,
					] table[x=u13660pre,y=f13660pre]{./I-shape/FU13660pre.txt};
					
					\addplot[
					color=blue, densely dashed,
					] table[x=u18510pre,y=f18510pre]{./I-shape/FU18510pre.txt};
					\addplot[
					scatter, only marks,
					mark=*,mark options={solid},
					scatter src = explicit symbolic,
					mark size=0.75pt,
					] table[x=u4113mor,y=f4113mor,meta=label]{./I-shape/FU4113mor.txt};
					
					\addplot[
					scatter, only marks,
					mark=*,mark options={solid},
					scatter src = explicit symbolic,
					mark size=0.75pt,
					] table[x=u8552mor,y=f8552mor,meta=label]{./I-shape/FU8552mor.txt};
					
					\addplot[
					scatter, only marks,
					mark=*,mark options={solid},
					scatter src = explicit symbolic,
					mark size=0.75pt,
					] table[x=u13660mor,y=f13660mor,meta=label]{./I-shape/FU13660mor.txt};
					
					\addplot[
					scatter, only marks,
					mark=*,mark options={solid},
					scatter src = explicit symbolic,
					mark size=0.75pt,
					] table[x=u18510mor,y=f18510mor,meta=label]{./I-shape/FU18510mor.txt};
					\legend{\rm{4113 full},\rm{8552 full},\rm{13660 full},\rm{18510 full},\rm{4113 MOR},\rm{8552 MOR},\rm{13660 MOR}, \rm{18510 MOR}}
					\node[align=center] at (16,23) {(a)};
				\end{axis}
			\end{tikzpicture}
		\end{minipage}
		\begin{minipage}[b]{0.5\linewidth}
			\begin{tikzpicture}
				\begin{axis}[
					xmin=0,xmax=32,
					ymin=1e-13,ymax=5e-3,
					scale only axis,
					ticklabel style={font=\large},
					xlabel={number of damage modes $\bm{l}$ $\bm{[-]}$}, 
					ylabel={Absolute error $\bm{\varepsilon _{\bar{D}}} \, \bm{[-]}$},
					line width=0.5pt,
					ymode=log,
					log basis y={10},
					font=\bfseries,
					legend pos=south west,
					legend style={fill =none,font=\tiny}
					]
					\addplot[
					line width=0.75pt,
					color=red, densely dashed,
					mark=star,
					mark size=2.5pt,
					mark options={solid},
					] table[x=modesD,y=D4113]{./I-shape/DErrors.txt};
					
					\addplot[
					line width=0.75pt,
					color=green, densely dashed,
					mark=triangle*,
					mark size=2.5pt,
					mark options={solid},
					] table[x=modesD,y=D8552]{./I-shape/DErrors.txt};
					
					\addplot[
					line width=0.75pt,
					color=violet, densely dashed,
					mark=*,mark size=2pt,
					mark options={solid},
					] table[x=modesD,y=D13660]{./I-shape/DErrors_13660.txt};
					
					\addplot[
					line width=0.75pt,
					color=blue, densely dashed,
					mark=pentagon*,
					mark size=2pt,
					mark options={solid},
					] table[x=modesD,y=D18510]{./I-shape/DErrors_18510.txt};
					
					\legend{
						{\rm{4113, $k=120$, $m=5$}},
						{\rm{8552, $k=120$, $m=5$}},
						{\rm{13660, $k=120$, $m=5$}},
						{\rm{18510, $k=120$, $m=5$}}
					}
					\node[align=center] at (29,5e-4) {(b)};
				\end{axis}
			\end{tikzpicture}
		\end{minipage}
		\caption{(a) Force-displacement curve comparison between simulations carried out using the full-order and the reduced-order models for the I-shaped specimen. The number of modes for the displacement, damage, and temperature fields are defined as $k=50$, $l=25$, and $m=5$, respectively; (b) the absolute errors in non-local damage for different numbers of damage modes $l$.}
		\label{fig8}
	\end{figure*}
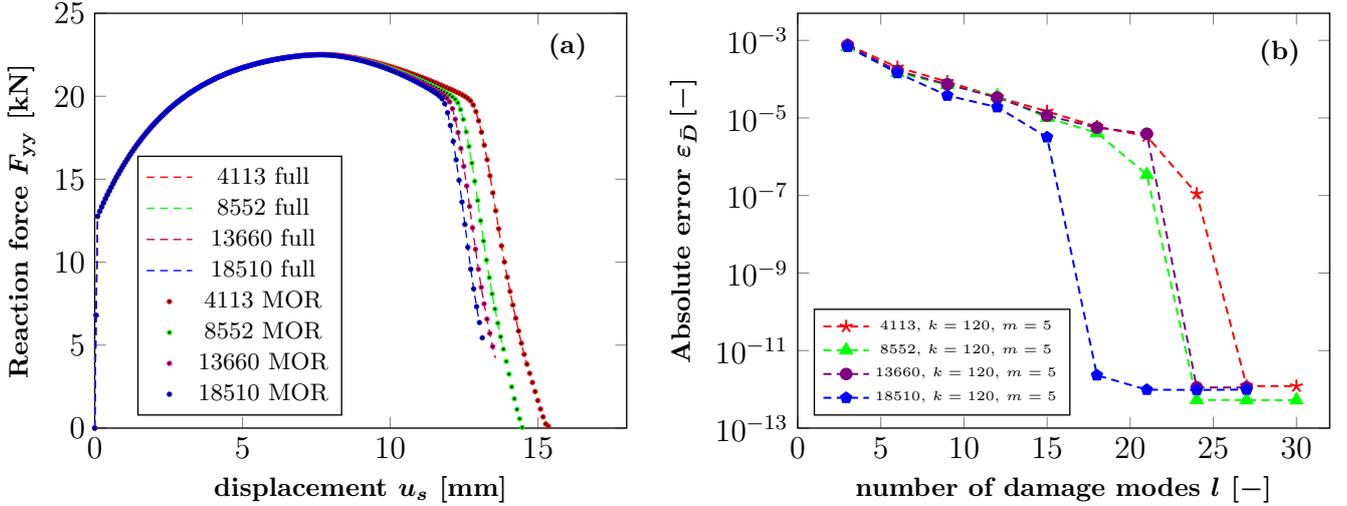

	When comparing the results obtained using the full-order model and the reduced-order model  (see, \hyperref[fig8]{Fig.10 (a)}), one can observe that the two match very well. To further study the influence of the error in the non-local damage variable, \hyperref[fig8]{Fig.10 (b)} shows $\varepsilon_{\bar{D}}$ for several damage modes $l$ ranging from 3 to 30. The values of $\varepsilon_{\bar{D}}$ are approximately the same, as the number of damage modes exceeds $25$, except for the mesh with 4113 elements.
	
	As the number of elements increases, there is a noticeable decrease in the number of the required damage modes compared to scenarios utilizing coarser meshes. This trend is illustrated in \hyperref[fig8]{Fig.10 (b)}, which mirrors the findings of \hyperref[fig7]{Fig.9 (b)}, wherein the absolute damage error in damage experiences a pronounced decrease and becomes stable beyond a specific number of damage modes. The absolute error of the non-local damage variable lies in the range of $5.34 \times 10^{-13} \sim 7.50 \times 10^{-4}$. Even for the lowest accuracy of $7.50 \times 10^{-4}$, the difference can hardly be identified in both the force-displacement curves and the contour plots. More detailed contour plots are given in \hyperref[fig9-1]{Fig.12 (b)} and \hyperref[fig9-1]{Fig.12 (c)}.
	
	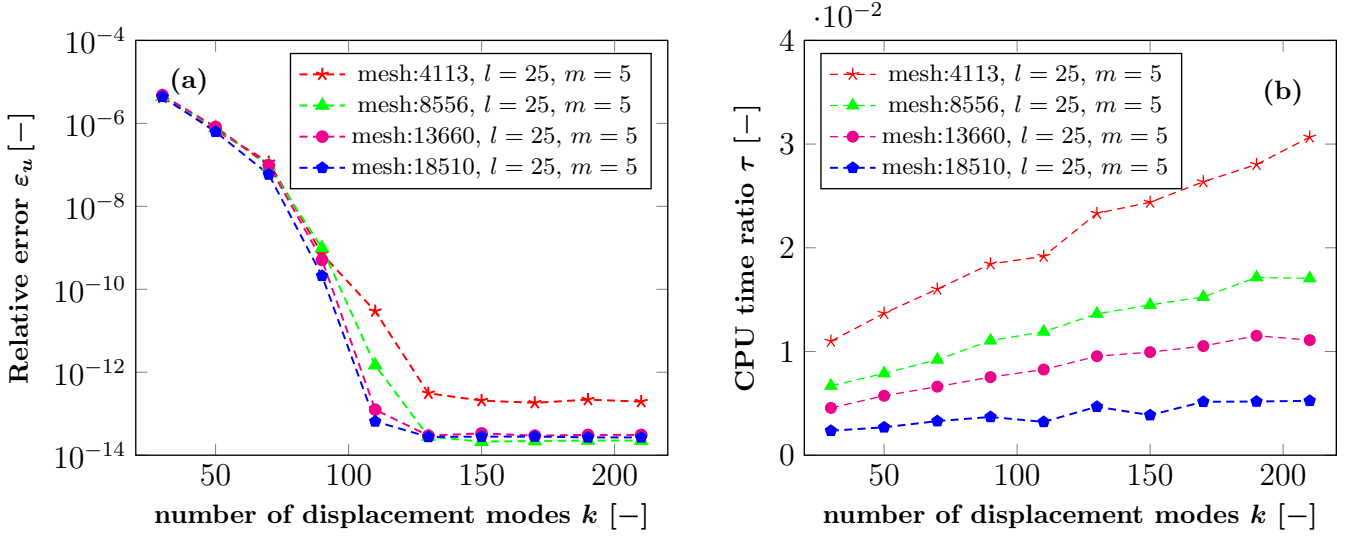
\begin{figure*}[!ht]
		\begin{minipage}[b]{0.5\linewidth}
			\begin{tikzpicture}
				\begin{axis}[
					xmin=20,xmax=220,
					ymin=1e-14,ymax=1e-4,
					scale only axis,
					ticklabel style={font=\large},
					xlabel={number of displacement modes $\bm{k}$ $\bm{[-]}$}, 
					ylabel={Relative error $\bm{\varepsilon _{u}} \, \bm{[-]}$},
					line width=0.5pt,
					ymode=log,
					log basis y={10},
					font=\bfseries,
					legend pos=north east,
					legend style={fill=none, font=\small}
					]
					\addplot[
					line width=0.75pt,
					color=red, densely dashed,
					mark=star,
					mark size=2.5pt,
					mark options={solid},
					] table[x=modesU,y=U4113]{./I-shape/UErrors.txt};
					
					\addplot[
					line width=0.75pt,
					color=green, densely dashed,
					mark=triangle*,
					mark size=2.5pt,
					mark options={solid},
					] table[x=modesU,y=U8552]{./I-shape/UErrors.txt};
					
					\addplot[
					line width=0.75pt,
					color=magenta, densely dashed,
					mark=*,mark size=2pt,
					mark options={solid},
					] table[x=modesU,y=U13660]{./I-shape/UErrors.txt};
					
					\addplot[
					line width=0.75pt,
					color=blue, densely dashed,
					mark=pentagon*,
					mark size=2pt,
					mark options={solid},
					] table[x=modesU,y=U18510]{./I-shape/UErrors.txt};
					
					\legend{
						{\rm{mesh:4113, $l=25$, $m=5$}},
						{\rm{mesh:8556, $l=25$, $m=5$}},
						{\rm{mesh:13660, $l=25$, $m=5$}},
						{\rm{mesh:18510, $l=25$, $m=5$}}
					}
					\node[align=center] at (40,1e-5) {(a)};
				\end{axis}
			\end{tikzpicture}
		\end{minipage}
		\begin{minipage}[b]{0.5\linewidth}
			\centering
			\begin{tikzpicture}
				\begin{axis}[
					clip=true,
					xmin=20,xmax=220,
					ymin=0,ymax=0.04,
					scale only axis,
					ticklabel style={font=\large},
					xlabel={number of displacement modes $\bm{k}$ $\bm{[-]}$}, 
					ylabel={CPU time ratio $\bm{\tau}$ $\bm{[-]}$},
					line width=0.5pt,
					font=\bfseries,
					legend pos=north west,
					legend style={fill=none, font=\small}
					]
					\addplot[
					color=red, densely dashed,
					mark=star,
					mark size=2.5pt,
					mark options={solid},
					] table[x=modesT,y=T4113]{./I-shape/TCPU.txt};
					
					\addplot[
					color=green, densely dashed,
					mark=triangle*,
					mark size=2.5pt,
					mark options={solid},
					] table[x=modesT,y=T8552]{./I-shape/TCPU.txt};
					
					\addplot[
					color=magenta, densely dashed,
					mark=*,mark size=2pt,
					mark options={solid},
					] table[x=modesT,y=T13660]{./I-shape/TCPU.txt};
					
					\addplot[
					line width=0.75pt,
					color=blue, densely dashed,
					mark=pentagon*,
					mark size=2pt,
					mark options={solid},
					] table[x=modesT,y=T18510]{./I-shape/TCPU.txt};
					
					\legend{
						{\rm{mesh:4113, $l=25$, $m=5$}},
						{\rm{mesh:8556, $l=25$, $m=5$}},
						{\rm{mesh:13660, $l=25$, $m=5$}},
						{\rm{mesh:18510, $l=25$, $m=5$}}
					}
					\node[align=center] at (200,3.5e-2) {(b)};
				\end{axis}
			\end{tikzpicture}
		\end{minipage}
		\caption{(a) Comparison of relative error of the displacement $\varepsilon_{u}$ over different displacement modes $k$; (b) CPU time ratio $\tau$ comparison for $k$ from $30$ to $210$.}
		\label{fig9}
	\end{figure*}
	
	Next, by comparing the relative error of the displacements $\varepsilon _{u}$ for a varying number of modes ranging from $k=30$ to $k=210$, it can be observed that $\varepsilon_{u}$ becomes stable around $10^{-14}$ in \hyperref[fig9]{Fig.11 (a)}, except for the mesh with 4113 elements ($\varepsilon_{u} \approx 10^{-13}$). The higher relative error observed for the mesh with 4113 elements, as compared to the other meshes, can be attributed to the use of a comparatively lower reduced number of damage modes ($l=25$) rather than employing the stable number ($l=27$). This lower mode number fails to adequately represent the coarse-meshed specimen, as shown in \hyperref[fig8]{Fig.10 (b)}. It is interesting to observe that the CPU time ratio of Test 2 (\hyperref[fig9]{Fig.11 (b)}) is significantly smaller than that of Test 1 (\hyperref[fig5]{Fig.7 (b)}). Furthermore, the CPU time ratio for Test 2 follows the same trends as in Test 1 across various mesh discretizations, i.e. $\tau _{\rm{18510}} < \tau _{\rm{13660}} < \tau _{\rm{8552}} < \tau _{\rm{4113}}$. The reason is that the number of stable damage modes for Test 2 ($l=25$) is much smaller than for Test 1 ($l=155$), see \hyperref[fig3]{Fig.4 (a)} and \hyperref[fig7]{Fig.9 (a)}. Although the required stable damage mode number for Test 2 is smaller, this difference does not affect achieving a low and stable error for both $\varepsilon_{u}$ and  $\varepsilon_{\bar{D}}$. Furthermore, \hyperref[fig9]{Fig.11} confirms that finer meshes lead to higher reductions in the computational cost achieved by the novel multi-field decomposed MOR approach.	
	
	As the number of elements in the mesh increases, the mode numbers for the I-shaped specimen (Test 2) stabilize, leading to a stabilization in the relative error of the reduced-order model, as depicted in \hyperref[fig9]{Fig.11 (a)}. This phenomenon is also observed in \hyperref[fig7]{Fig.9 (b)}. Consequently, the combined process of finding converged meshes and corresponding stable mode numbers is essential for minimizing computational costs.
	
	\begin{figure*}[!ht]
		\centering
		\includegraphics[width=16cm]{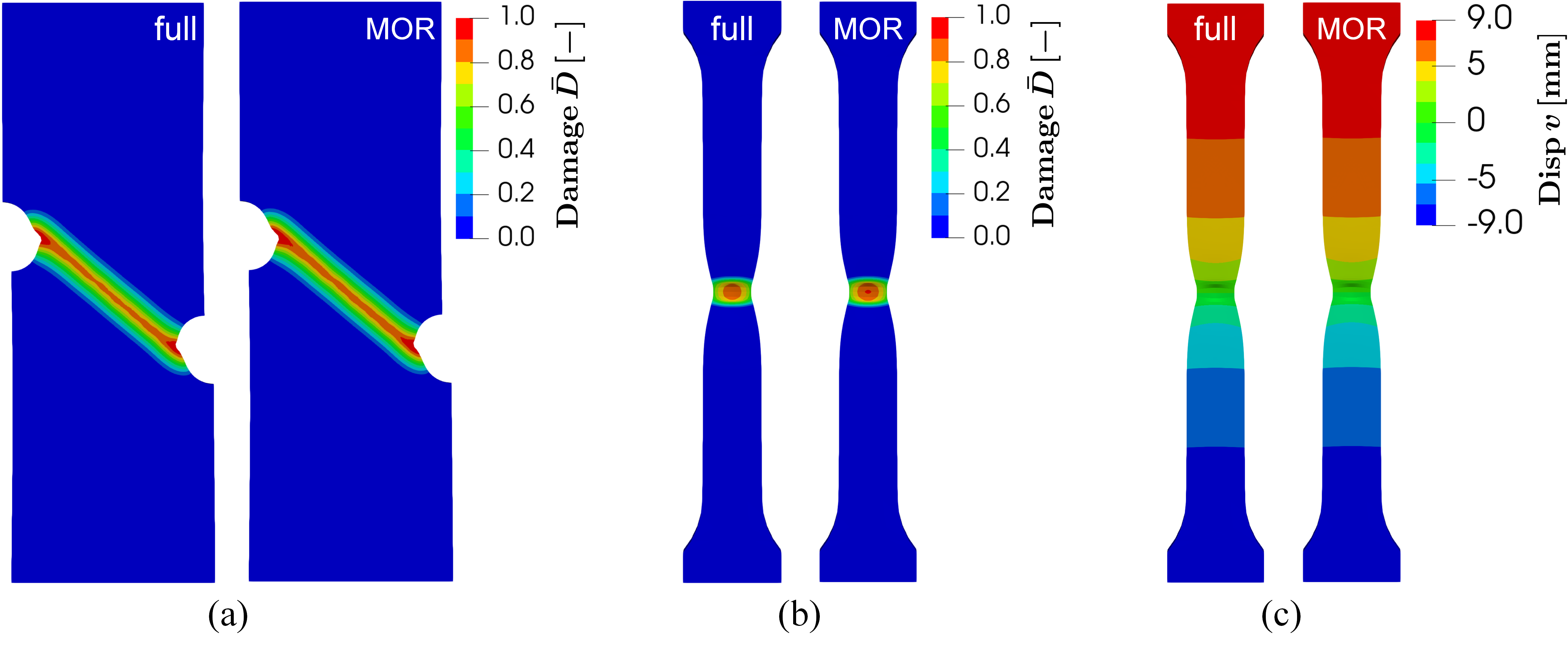}
		\caption{(a) Contour plots of the non-local damage $\bar{D}$ for the full-order and the reduced-order model (MOR, $k=120$, $l=25$, and $m=5$) simulations for a mesh with 8979 elements (Test 1); (b) contour plots of the non-local damage $\bar{D}$ for full-order model and reduced-order model (MOR, $k=50$, $l=25$, and $m=5$) simulations for a mesh with 4113 elements (Test 2); (c) contour plots of the vertical displacement $v$ for full-order model and reduced-order model (MOR, $k=50$, $l=25$, and $m=5$) simulations in Test 2.}
		\label{fig9-1}
	\end{figure*}
	
	\section{Results and discussion}
	\label{Results and discussion}
	\subsection{Converged versus non-converged nodal snapshots}
	\label{RAD}
	
	As stated earlier, the conventional MOR method has limitations in cases where material damage is considered. In order to investigate this further, three different MOR variants are additionally compared in the present study to emphasize the importance of using the novel multi-field decomposed MOR for simulations involving damage.
	\begin{itemize}
		\item Method-A: converged nodal solutions as snapshots in combination with the multi-field decomposition method ($k=120$, $l=25$, $m=5$); 
		\item Method-B: non-converged and converged nodal solutions as snapshots in combination with the multi-field decomposition method ($k=85$, $l=280$, $m=5$); 
		\item Method-C: non-converged and converged nodal solutions as snapshots without the field decomposition.
	\end{itemize}
	
	First, the standard MOR approach was applied in order to reduce the simulation cost for Test 1, specifically utilizing only converged nodal snapshots without field decomposition. However, this approach fails to accurately reconstruct the mechanical response, even in the elastoplastic region. For this reason, Method-C is now employed to refine the reduced-order model simulation and investigate it for varying additional mode numbers from 50 to 500.
	
	\begin{figure*}[!ht]
		\centering
		\begin{tikzpicture} [spy using outlines={rectangle, magnification=3, size=0.5cm, connect spies}]
			\begin{axis}[
				clip=true,
				xmin=0,xmax=4,
				ymin=0,ymax=5,
				scale only axis,
				log ticks with fixed point,
				ticklabel style={font=\large},
				xlabel={displacement $\bm{u}$ \textbf{[mm]}}, 
				ylabel={Reaction force $\bm{F_{\rm{xx}}}$ \textbf{[kN]}},
				line width=0.5pt,
				font=\bfseries,
				legend pos=south east,
				legend style={font=\tiny},
				scatter/classes={%
					a={mark=*,draw=blue}
				}
				]
				\addplot[
				line width=1.5pt,
				color=blue, 
				] table[x=ufull,y=ffull]{./Asy/snap_130step_compare/full-F-U.txt};
				
				\addplot[
				line width=0.75pt,
				color=red, densely dashed,
				] table[x=uMOR450,y=fMOR450]{./Asy/snap_130step_compare/mor-F-U-M450.txt};
				
				\addplot[
				line width=0.75pt,
				color=green, densely dashed,
				] table[x=uMOR500,y=fMOR500]{./Asy/snap_130step_compare/mor-F-U-M500.txt};
				
				\addplot[
				line width=0.75pt,
				color=black,
				densely dashed,
				] table[x=uMOR400,y=fMOR400]{./Asy/snap_130step_compare/mor-F-U-M400.txt};
				
				\addplot[
				line width=0.75pt,
				color=yellow, densely dashed,
				] table[x=uMOR350,y=fMOR350]{./Asy/snap_130step_compare/mor-F-U-M350.txt};
				
				\addplot[
				line width=0.75pt,
				color=magenta, densely dashed,
				] table[x=uMOR300,y=fMOR300]{./Asy/snap_130step_compare/mor-F-U-M300.txt};
				
				\addplot[
				line width=0.75pt,
				color=pink, densely dashed,
				] table[x=uMOR250,y=fMOR250]{./Asy/snap_130step_compare/mor-F-U-M250.txt};
				
				\addplot[
				line width=0.75pt,
				color=orange, densely dashed,
				] table[x=uMOR200,y=fMOR200]{./Asy/snap_130step_compare/mor-F-U-M200.txt};
				
				\addplot[
				line width=0.75pt,
				color=violet, densely dashed,
				] table[x=uMOR100,y=fMOR100]{./Asy/snap_130step_compare/mor-F-U-M100.txt};
				
				\addplot[
				line width=0.75pt,
				color=cyan, densely dashed,
				] table[x=uMOR150,y=fMOR150]{./Asy/snap_130step_compare/mor-F-U-M150.txt};
				
				\addplot[
				line width=0.75pt,
				color=purple, densely dashed,
				] table[x=uMOR50,y=fMOR50]{./Asy/snap_130step_compare/mor-F-U-M50.txt};

				\legend{full model,\rm{C, 450 modes},\rm{C, 500 modes},\rm{C, 400 modes},\rm{C, 350 modes},\rm{C, 300 modes},\rm{C, 250 modes},\rm{C, 200 modes},\rm{C, 100 modes},\rm{C, 150 modes},\rm{C, 50 modes}}
				
				\coordinate (spypoint) at (axis cs:0.48,3.8772); 
				\coordinate (spyviewer) at (axis cs:1.5,2.5);
				\spy[width=1.5cm,height=1.5cm] on (spypoint) in node [fill=none] at (spyviewer);
				
			\end{axis}
		\end{tikzpicture}
		\caption{Comparison of force-displacement curves across a range of mode numbers (50 to 500) using non-decomposed and additionally non-converged nodal snapshots (Method-C) in the MOR approach.}
		\label{fig10-1}
	\end{figure*}
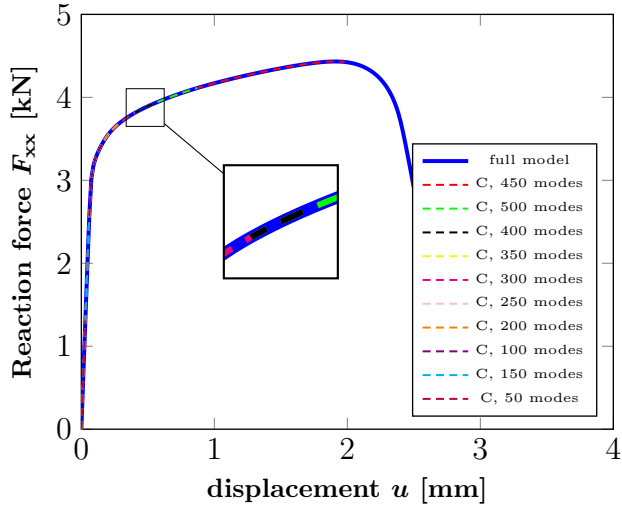
	
	\hyperref[fig10-1]{Fig.13} indicates that Method-C can accurately reproduce the elastoplastic and damage behavior in the structure before global softening occurs when increasing the number of modes. However, it is noteworthy that Method-C is limited to this region in the plot, it fails for the part leading to global softening. This observation is in accordance with previous research findings \cite{kerfriden2011bridging,kerfriden2013partitioned,selvaraj2024adaptive,kastian2023discrete}. To address this issue, this study proposes the above-explained novel multi-field decomposed MOR methodology. In the following, two different variants of the latter are investigated in more detail, namely Method-A and Method-B.
	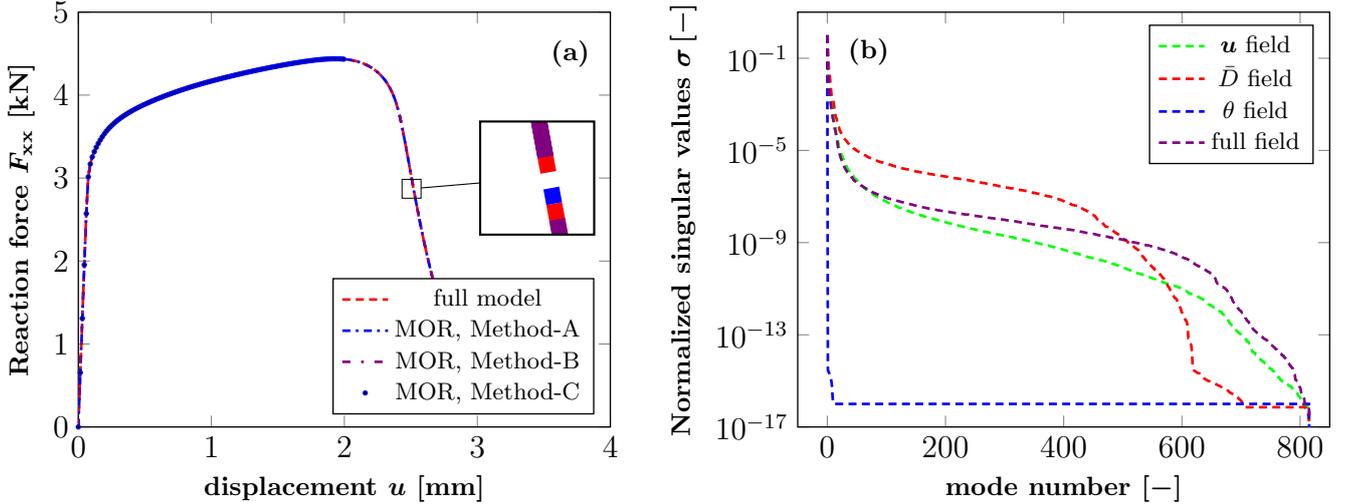
\begin{figure*}[!ht]
		\begin{minipage}[b]{0.5\linewidth}
			\centering
			\begin{tikzpicture} [spy using outlines={rectangle, magnification=6, size=0.5cm, connect spies}]
				\begin{axis}[
					clip=true,
					xmin=0,xmax=4,
					ymin=0,ymax=5,
					scale only axis,
					log ticks with fixed point,
					ticklabel style={font=\large},
					xlabel={displacement $\bm{u}$ \textbf{[mm]}}, 
					ylabel={Reaction force $\bm{F_{\rm{xx}}}$ \textbf{[kN]}},
					line width=0.5pt,
					font=\bfseries,
					legend pos=south east,
					legend style={font=\small},
					scatter/classes={%
						a={mark=*,draw=blue}
					}
					]
					\addplot[
					line width=1pt,
					color=red, densely dashed,
					] table[x=u2912pre,y=f2912pre]{./Asy/FU2912pre.txt};
					
					\addplot[
					line width=1pt,
					color=blue, dashdotted,
					] table[x=u2912mor,y=f2912mor]{./Asy/con-vs-incon.txt};
					
					\addplot[
					line width=1pt,
					color=violet, loosely dashdotted,
					] table[x=u2912nrmor,y=f2912nrmor]{./Asy/con-vs-incon.txt};
					
					\addplot[
					scatter, only marks,
					mark=*,mark options={solid},
					scatter src=explicit symbolic,
					mark size=0.75pt,
					] table[x=u2912nr130,y=f2912nr130,meta=label]{./Asy/NR130-2912.txt};
					
					\legend{\rm{full model},\rm{MOR, Method-A},\rm{MOR, Method-B},\rm{MOR, Method-C}}
					\node[align=center] at (3.7,4.5) {(a)};
					
					\coordinate (spypoint) at (axis cs:2.505,2.871); 
					\coordinate (spyviewer) at (axis cs:3.45,3);
					\begin{scope}
						\spy[width=1.5cm,height=1.5cm] on (spypoint) in node [fill=none] at (spyviewer);
					\end{scope}
					
				\end{axis}
							
							
							
			\end{tikzpicture}
		\end{minipage}
		\begin{minipage}[b]{0.5\linewidth}
			\begin{tikzpicture}
				\begin{axis}[
					xmin=-50,xmax=850,
					ymin=1e-17,ymax=10,
					scale only axis,
					ticklabel style={font=\large},
					xlabel={mode number $\bm{[-]}$}, 
					ylabel={Normalized singular values $\bm{\sigma} \, \bm{[-]}$},
					line width=0.5pt,
					ymode=log,
					log basis y={10},
					font=\bfseries,
					legend pos=north east,
					legend style={font=\small}
					]
					\addplot[
					line width=1pt,
					color=green, densely dashed,
					] table[x=modes,y=U]{./Asy/modes-SVD-NR.txt};
					
					\addplot[
					line width=1pt,
					color=red, densely dashed,
					] table[x=modes,y=D]{./Asy/modes-SVD-NR.txt};
					
					\addplot[
					line width=1pt,
					color=blue, densely dashed,
					] table[x=modes,y=T]{./Asy/modes-SVD-NR.txt};
					
					\addplot[
					line width=1pt,
					color=violet, densely dashed,
					] table[x=modes,y=full]{./Asy/modes-SVD-NR.txt};
					
					\legend{
						{\rm{$\bm{u}$ field}},
						{\rm{$\bar{D}$ field}},
						{\rm{$\theta$ field}},
						{\rm{full field}}
					}
					\node[align=center] at (70,0.2) {(b)};
				\end{axis}
			\end{tikzpicture}
		\end{minipage}
		\caption{(a) Force-displacement curve comparison between different MOR variants and the full-order model for Test 1. Method-A: converged nodal snapshots with multi-field decomposition method ($k=120$, $l=25$, $m=5$); Method-B: non-converged and converged nodal snapshots with multi-field decomposition method ($k=85$, $l=280$, $m=5$); Method-C: non-convergent and converged nodal snapshots without field decomposition; (b) singular value decomposition for the non-converged and converged nodal snapshots with multi-field decomposition for Test 1 with  2912 elements.}
		\label{fig10}
	\end{figure*}
	
	\hyperref[fig10]{Fig.14 (a)} illustrates that both Method-A and Method-B are capable of reconstructing mechanical responses within the global softening region induced by damage, owing to their multi-field decomposed character. However, due to the higher number of modes, the computational cost of Method-B ($k + l + m = 370$) is 2.5 times higher than that of Method-A ($k + l + m = 150$). The comparison of the normalized singular values between Method-B (see \hyperref[fig10]{Fig.14 (b)}) and Method-A (see \hyperref[fig3]{Fig.4 (a)}) indicates an approximately four times difference in mode numbers, accompanied by different developments for the individual field quantities. These findings underscore the necessity of employing the novel multi-field decomposition to effectively address damage-induced global softening in reduced-order simulations, see \hyperref[T2]{Table 2}.
	
	\begin{table}[!ht]
		\fontsize{8pt}{8pt}
		\centering
		\begin{tabular}{lccccc}
			\hline
			\multirow{2}{*}{} & \multicolumn{1}{l}{\multirow{2}{*}{\textbf{\begin{tabular}[c]{@{}l@{}}only converged \\ snapshots\end{tabular}}}} & \multicolumn{1}{l}{\multirow{2}{*}{\textbf{\begin{tabular}[c]{@{}l@{}}decompose \\ field\end{tabular}}}} & \multicolumn{1}{l}{\multirow{2}{*}{\textbf{\begin{tabular}[c]{@{}l@{}}elastoplastic \\ region\end{tabular}}}} & \multicolumn{1}{l}{\multirow{2}{*}{\textbf{\begin{tabular}[c]{@{}l@{}}global softening \\ region\end{tabular}}}} & \multicolumn{1}{l}{\multirow{2}{*}{\textbf{\begin{tabular}[c]{@{}l@{}}computational \\ cost\end{tabular}}}} \\
			& \multicolumn{1}{l}{}                                                                                              & \multicolumn{1}{l}{}                                                                                     & \multicolumn{1}{l}{}                                                                                         & \multicolumn{1}{l}{}                                                                                             & \multicolumn{1}{l}{}                                                                                        \\  \hline
			\textbf{Method-A} & $\surd$                             & $\surd$                          & $\surd$                            & $\surd$                         &  ++                        \\
			\textbf{Method-B} & $\times$                             & $\surd$                          & $\surd$                            & $\surd$                         & +                          \\
			\textbf{Method-C} & $\times$                             & $\times$                          & $\surd$                            & $\times$                         & -                           \\
			\textbf{Classic}  & $\surd$                             & $\times$                          & $\times$                            & $\times$                         & -                           \\ \hline
		\end{tabular}
		\label{T2}
		\caption{A comparison of the computational performance among different MOR variants in Test 1 simulation. The markers $\surd$ and $\times$ denote the \textit{yes} and \textit{no}, respectively.}
	\end{table}
	
	The primary advantage and methodological innovation for the multi-field decomposed MOR lies in its capacity to individually project the coupled nonlinear fields, such as displacement, damage, and temperature onto independent POD-subspaces. In contrast to conventional methods employing single-field projection, the detailed influential factors and contributions of displacement, damage, and temperature fields in multi-physics problems are distinctly revealed through comparative error analysis across displacement and damage modes, as well as normalized singular value developments over different fields.
	
	However, the limitation of the multi-field decomposed MOR becomes apparent, analogous to the shortcomings observed in conventional MOR approaches, as articulated by \citet{chaturantabut2010nonlinear}: \textit{For a given system, the solution space often converges towards a low-dimensional manifold}. This method is limited to scenarios with the same mesh which is a characteristic inherent to global Galerkin-projection-based MOR techniques. When using different meshes, the Galerkin projection process becomes infeasible due to the dimensional inconsistent nature of the projection matrix. Finally, the authors believe that multi-field decomposed MOR is one of the missing keys for accurately and effectively addressing damage-induced global softening.
	
	\subsection{Damage modes}
	
\begin{figure*}[!ht]
		\centering
		\includegraphics[width=18cm]{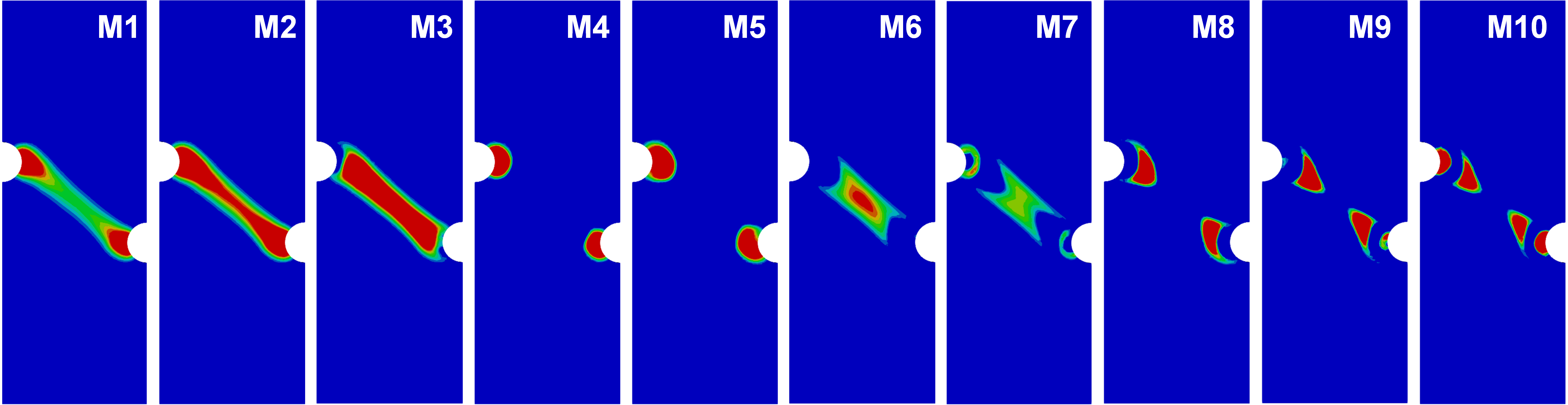}
		\caption{Schematic illustration of the dominant non-local damage modes ($\bar{D}$) of the asymmetrically notched specimen, ranging from mode 1 (M1) to mode 10 (M10).}
		\label{fig12}
	\end{figure*}

\begin{figure*}[!ht]
		\centering
		\includegraphics[width=18cm]{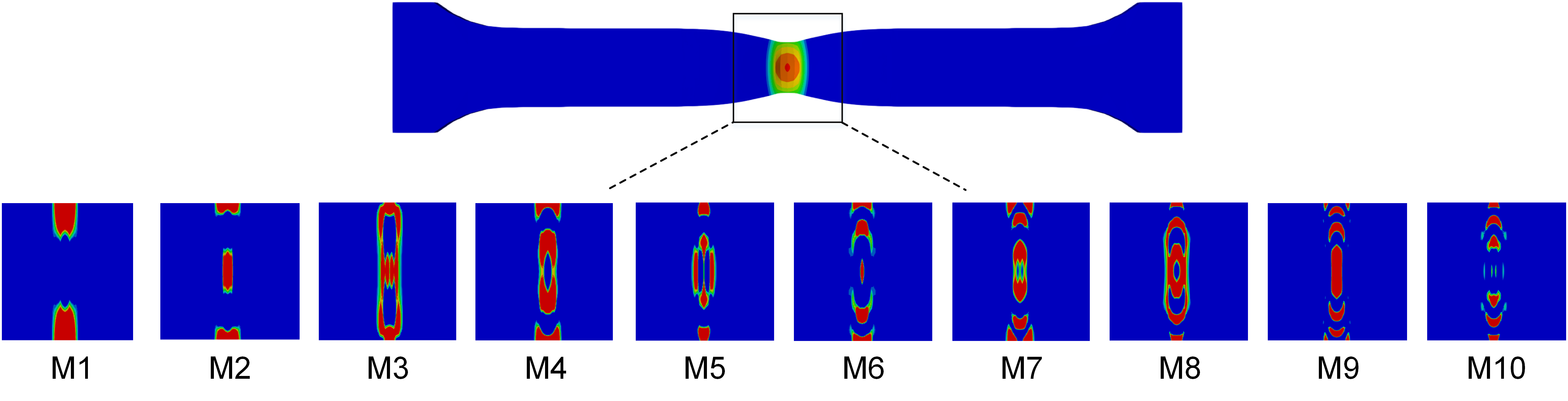}
		\caption{Schematic illustration of the dominant non-local damage modes ($\bar{D}$) of the I-shaped specimen ranging from mode 1 (M1) to mode 10 (M10).}
		\label{fig13}
	\end{figure*}

	This study investigates and compares the error in the non-local damage fields obtained by using the full-order and reduced-order model in the simulations, employing additionally various numbers of damage modes as illustrated in \hyperref[fig4]{Fig.5 (b)} and \hyperref[fig8]{Fig.10 (b)}. While these modes are crucial for reconstructing the damage evolution in MOR simulations, their physical meaning remains somewhat unclear. To overcome this, the dominant damage modes are shown in \hyperref[fig12]{Fig.15} and \hyperref[fig13]{Fig.16}.
	
	From the structural perspective, each mode of the non-local damage projection matrix $\bm{\Phi}_{\bar{D}, \, l}$ (see \hyperref[eq43]{Eq.\eqref{eq43}} and \hyperref[eq44]{Eq.\eqref{eq44}}) represents a specific aspect of damage states. From modes 1 to 10, non-zero values of $\bar{D}$ cluster predominantly around the notched areas. Understanding this pattern requires insight into the stages of damage evolution: initialization, propagation, branching or merging, and fracture. \hyperref[fig12]{Fig.15} illustrates how different modes correspond to these stages. Modes M1, M4, and M5 depict damage initiation at the notch tip, while M2 and M3 represent orientation and propagation. Modes M6 to M7 signify the tiny potential branching, and M8 to M10 indicate the merging of damaged regions.
	
	Compared to \hyperref[fig12]{Fig.15}, the dominant damage modes of the I-shaped specimen (\hyperref[fig13]{Fig.16}) are simpler. Considering the dominant factor in the evolution of damage within the specimen, shear stress plays an essential role in the loading of the asymmetrically notched specimen, due to the slippage of the specimen in vertical direction. The dominant factor for damage evolution in the I-shaped specimen test is normal stress rather than shear stress. Consequently, the final diffuse damage zone forms as a vertical band in the center of the specimen, see \hyperref[fig13]{Fig.16}. To precisely replicate this damage band using the projection matrix of damage in MOR, it is reasonable that the non-zero data should be predominantly located within this damage band. It can be observed in \hyperref[fig12]{Fig.15} and \hyperref[fig13]{Fig.16} that the distinct patterns of non-local damage modes are similar to the corresponding states of damage evolution.
	\section{Conclusion}
	\label{Conclusion}
	In this study, a novel multi-field decomposed model order reduction (MOR) approach was first proposed to address thermo-mechanically coupled gradient-extended damage simulations. By utilizing multi-field decomposition, complex and coupled nonlinear coupling problems were systematically converted into linear projection tasks. The integration of multi-field decomposition into the standard POD-G approach can lead to much higher accuracy in addressing complex multi-physics and damage-induced global softening simulations.
	This was shown in this study via two numerical examples, confirming that the MOR technique based on multi-field decomposition outperforms the conventional POD-based MOR approach, particularly in cases where damage-induced global softening dominates the structural response. To quantify the efficiency of the proposed method, comparative analyses of the relative error in displacement and the absolute error in damage were conducted by comparing the results obtained by using the full-order and reduced-order models. The findings of this study indicate that the novel multi-field decomposed model order reduction approach exhibits a noticeable computational efficiency as well as a high level of accuracy for thermo-mechanically coupled gradient-extended damage simulations.

	\section{Acknowledgement}
	\label{Acknowledgement}
	Q. Zhang, J. Kehls, T. Brepols, and S. Reese acknowledge the financial support of subproject B05 within SFB/TRR 339 (project number: 453596084). Q. Zhang acknowledges the valuable discussions as well as suggestions from Erik Prume and Hooman Danesh. J. Zhang acknowledges the financial support from the China Scholarship Council (CSC). S. Ritzert and S. Reese acknowledge the financial support of the subproject A01 within SFB/TRR 280 (project number: 417002380).
	
	\section{Conflict of interest}
	The authors declare that they have no known competing financial interests or personal relationships that could have appeared to influence the work reported in this paper.
	
	\section{Contributions by the authors}
	\textbf{Qinghua Zhang}: Conceptualization, Methodology, Theoretical description, Implementation, Data Curation, Formal analysis, Investigation, Visualization, Writing – original draft, Writing – review \& editing.
	\textbf{Stephan Ritzert}: Theoretical description, Implementation, Visualization, Review \& editing.
	\textbf{Jian Zhang}: Theoretical description, Writing – review \& editing.
	\textbf{Jannick Kehls}: Review \& editing.
	\textbf{Stefanie Reese}: Conceptualization, Supervision, Funding acquisition.
	\textbf{Tim Brepols}: Conceptualization, Theoretical description, Supervision, Funding acquisition, Writing – review \& editing.
	\newpage
	\appendix 
	
	\section{1}
	\label{AppendixA}
	The total internal heat generation is expressed as $r_{\rm{int}}={r_e}+{r_p}+{r_d}$ with elastic $\left(r_{e}\right)$, plastic $\left(r_{p}\right)$, and damaged $\left(r_{d}\right)$ parts. They are given as:
	\begin{equation}
		{r_e} = \theta \, \left( {\frac{1}{2} \, \frac{{\partial {\bm{S}}}}{{\partial \theta }} - \alpha \, \frac{{\partial {\bm{S}}}}{{\partial {\bm{C}}}} \, {{:}} \, {\bm{C}} - \alpha \, {\bm{S}}} \right) \, {:} \, {\dot {\bm{C}}},
		\label{eq23}
	\end{equation}
	
	\begin{equation}
		{r_p} = \frac{1}{2} \, {\bm{C}} \, \left( {{\bm{S}} - \frac{{\partial {\bm{S}}}}{{\partial \theta }} \, \theta } \right) \, {:} \, {\bm{C}}_p^{ - 1} \, {{\dot{\bm{C}}}_p} - \frac{1}{2} \, \left( {{\bm{X}} - \frac{{\partial {\bm{X}}}}{{\partial \theta }} \, \theta } \right) \, {:} \, {{\dot{\bm{C}}}_p} - \left( {{q_p} - \frac{{\partial {q_p}}}{{\partial \theta }} \, \theta } \right) \, {{\dot \xi }_p},
		\label{eq24}
	\end{equation}
	
	\begin{equation}
		{r_d} = \left( {Y - \frac{{\partial Y}}{{\partial \theta }} \, \theta } \right) \, \dot D - \left( {{q_d} - \frac{{\partial {q_d}}}{{\partial \theta }} \, \theta } \right) \, {{\dot \xi }_d} - \theta \, \left( {\frac{{\partial {a_i}}}{{\partial \theta }} - \alpha \, \frac{{\partial {\bm{S}}}}{{\partial \bar D}} \, {:} \, {\bm{C}}} \right) \, \dot {\bar D} + \theta \, \frac{{\partial {\bm{b}_i}}}{{\partial \theta }} \cdot \nabla \dot {\bar D}.
		\label{eq25}
	\end{equation}
	
	\section{1}
	\label{AppendixB}
	The Claus-Duhem inequality can be derived with a time derivative of the elastic and plastic right Cauchy-Green tensors $\dot{\bm{C}}_e$ and $\dot{\bm{C}}_p$, respectively:
	\begin{equation}
		\begin{gathered}
			{{{\dot{\bm C}}}_e} =  - {\bm{F}}_p^{{-\text{T}}} \, {\dot{\bm F}}_p^{\text{T}} \, {{\bm{C}}_e} + \frac{1}{{{\vartheta ^2}}} \, {\bm{F}}_p^{{-\text{T}}} \, {\dot{\bm C}} \, {{\bm F}}_p^{ - 1} - {{\bm{C}}_e} \, {{{\dot{\bm F}}}_p} \, {\bm{F}}_p^{ - 1} - \frac{2}{{{\vartheta ^2}}} \, {\bm{F}}_p^{{-\text{T}}} \, {\bm{C\, F}}_p^{ - 1} \, \alpha \, \dot \theta ,  \hfill \\
			{{{\dot{\bm C}}}_p} = 2 \, {\bm{F}}_p^{\text{T}} \, {{\bm{D}}_p}{{\bm{F}}_p} \quad {\text{with}} \quad {{\bm{D}}_p} = \frac{1}{2} \, {\bm{F}}_p^{{-\text{T}}} \, {{{\dot{\bm C}}}_p} \, {\bm{F}}_p^{ - 1}. \hfill \\ 
		\end{gathered}
	\end{equation}
	Then the \hyperref[eq10]{Eq.\eqref{eq10}} can be expressed as:
	\begin{equation}
		\begin{gathered}
			\left( {{\bm{S}} - 2\, {f_d}(D) \, \frac{1}{{{\vartheta ^2}}} \, {\bm{F}}_p^{ - {1}} \, \frac{{\partial {\psi _e}}}{{\partial {{\bm{C}}_e}}} \, {\bm{F}}_p^{ - {\text{T}}}} \right):\frac{1}{2} \,  \dot {\bm{C}} - {f_d}(D) \, \left[ {2 \, {\bm{F}}_p^{\text{T}} \, \frac{{\partial {\psi _p}}}{{\partial {{\bm{C}}_p}}} \, {{\bm{F}}_p} - 2 \, \frac{{\partial {\psi _e}}}{{\partial {{\bm{C}}_e}}} \, {{\bm{C}}_e}} \right]:{{\bm{D}}_p} - {f_d}(D) \, \frac{{\partial {\psi _p}}}{{\partial {\xi _p}}} \, {{\dot \xi }_p} \hfill \\
			- \left( {\eta  + \frac{{\partial \psi }}{{\partial \theta }} - 2 \, {f_d}(D) \,  \bm{C} \, \alpha \, \frac{1}{{{\vartheta ^2}}}{\bm{F}}_p^{ - {1}} \, \frac{{\partial {\psi _e}}}{{\partial {{\bm{C}}_e}}} \, {\bm{F}}_p^{ - {\text{T}}}:{\bm{I}}} \right) \, \dot \theta  - \left[ {\frac{{\partial {f_d}(D)}}{{\partial D}} \, \left( {{\psi _e} + {\psi _p}} \right) + \frac{{\partial {\psi _{\bar d}}}}{{\partial D}}} \right]\dot D - \frac{{\partial {\psi _d}}}{{\partial {\xi _d}}} \, {{\dot \xi }_d} \hfill \\
			\left( {{a_i} - \frac{{\partial {\psi _{\bar d}}}}{{\partial \bar D}}} \right) \, \dot {\bar D} + \left( {{\bm{b}_i} - \frac{{\partial {\psi _{\bar d}}}}{{\partial \nabla \bar D}}} \right) \, \nabla \dot {\bar D} - \frac{1}{\theta } \, q \, \nabla \theta  \geqslant 0. \hfill \\ 
		\end{gathered}
	\end{equation}
	Therefore, the Mandel stress $\bm{M}$, the back-stress $\bm{\chi}$, and the corresponding stress tensor $\bm{X}$ are expressed as:
	\begin{equation}
		\begin{gathered}
			{\bm{M}} = 2 \, {f_d}(D) \, {{\bm{C}}_e} \, \frac{{\partial {\psi _e}}}{{\partial {{\bm{C}}_e}}} = 2{f_d}(D) \, {{\bm{C}}_e} \, \frac{1}{{{\vartheta ^2}}} \, {\bm{F}}_p^{ - {\text{T}}} \, {\bm{CF}}_p^{ - 1} \, \frac{{\partial {\psi _e}}}{{\partial {{\bm{C}}_e}}} = {\bm{F}}_p^{ - {\text{T}}} \, {\bm{SCF}}_p^{\text{T}}, \hfill \\
			{\bm{\chi }} \, {\text{: = }} \, 2{f_d}(D) \, {{\bm{F}}_p} \, \frac{{\partial {\psi _p}}}{{\partial {{\bm{C}}_p}}} \, {\bm{F}}_p^{\text{T}}, \quad {\bm{X}} \, {\text{:}} = 2{f_d}(D) \, \frac{{\partial {\psi _p}}}{{\partial {{\bm{C}}_p}}}. \hfill \\ 
		\end{gathered}
	\end{equation}
	The detailed operation of \hyperref[eq14]{Eq.\eqref{eq14}} can be expressed with $\bm{C}$, $\bm{C}_p$, $\bm{S}$, and $\bm{X}$ as: 
	\begin{equation}
		{\bm{M}} \, {\text{:}} \, {{\bm{D}}_p} = \frac{1}{2} \, {\text{tr}}\left( {{\bm{F}}_p^{ - {\text{T}}} \, {\bm{S \, C \, F}}_p^{\text{T}} \, {\bm{F}}_p^{ - {\text{T}}} \, {{{\dot{\bm C}}}_p}{\bm{F}}_p^{ - 1}} \right) = \frac{1}{2} \, {\text{tr}}\left( {{\bm{S\, C \, F}}_p^{\text{T}} \, {\bm{F}}_p^{ - {\text{T}}} \, {{{\dot{\bm C}}}_p} \, {\bm{F}}_p^{ - 1} \, {\bm{F}}_p^{ - {\text{T}}}} \right) = \frac{1}{2} \, {\bm{S \, C}} \, {\text{:}} \, {{{\dot{\bm C}}}_p} \, {\bm{C}}_p^{ - 1},
	\end{equation}
	\begin{equation}
		\begin{gathered}
			{\bm{\chi }} \, {\text{:}} \, {{\bm{D}}_p} = {\text{tr}}\left( {2{f_d}(D) \, {{\bm{F}}_p} \, \frac{{\partial {\psi _p}}}{{\partial {{\bm{C}}_p}}} \, {\bm{F}}_p^{\text{T}} \, \frac{1}{2} \, {\bm{F}}_p^{ - {\text{T}}} \, {{{\dot{\bm C}}}_p} \, {\bm{F}}_p^{ - 1}} \right) = \frac{1}{2} \, {\text{tr}}\left( {2{f_d}(D) \, \frac{{\partial {\psi _p}}}{{\partial {{\bm{C}}_p}}} \, {\bm{F}}_p^{\text{T}} \, {\bm{F}}_p^{ - {\text{T}}} \, {{{\dot{\bm C}}}_p} \, {\bm{F}}_p^{ - 1} \, {{\bm{F}}_p}} \right), \hfill \\
			{\bm{\chi }} \, {\text{:}} \, {{\bm{D}}_p} = \frac{1}{2} \, {\bm{X}} \, {\text{:}} \, {{{\dot{\bm C}}}_p}. \hfill \\ 
		\end{gathered}
	\end{equation}
	
	\newpage	
	\section*{References}
	\bibliography{Article}

\begin{thebibliography}{67}
\providecommand{\natexlab}[1]{#1}
\providecommand{\url}[1]{\texttt{#1}}
\expandafter\ifx\csname urlstyle\endcsname\relax
  \providecommand{\doi}[1]{doi: #1}\else
  \providecommand{\doi}{doi: \begingroup \urlstyle{rm}\Url}\fi

\bibitem[Ahrens et~al.(2005)Ahrens, Geveci, Law, Hansen, and
  Johnson]{ahrens200536}
J.~Ahrens, B.~Geveci, C.~Law, C.~Hansen, and C.~Johnson.
\newblock 36-paraview: An end-user tool for large-data visualization.
\newblock \emph{The visualization handbook}, 717:\penalty0 717 – 731, 2005.

\bibitem[Alter et~al.(2000)Alter, Brown, and Botstein]{alter2000singular}
O.~Alter, P.~O. Brown, and D.~Botstein.
\newblock Singular value decomposition for genome-wide expression data
  processing and modeling.
\newblock \emph{Proceedings of the National Academy of Sciences of the United
  States of America}, 97\penalty0 (18):\penalty0 10101 – 10106, 2000.

\bibitem[Ambati et~al.(2016)Ambati, Kruse, and De~Lorenzis]{ambati2016phase}
M.~Ambati, R.~Kruse, and L.~De~Lorenzis.
\newblock A phase-field model for ductile fracture at finite strains and its
  experimental verification.
\newblock \emph{Computational Mechanics}, 57:\penalty0 149--167, 2016.

\bibitem[Ames et~al.(2009)Ames, Srivastava, Chester, and Anand]{ames2009thermo}
N.~M. Ames, V.~Srivastava, S.~A. Chester, and L.~Anand.
\newblock A thermo-mechanically coupled theory for large deformations of
  amorphous polymers. part ii: Applications.
\newblock \emph{International Journal of Plasticity}, 25\penalty0 (8):\penalty0
  1495--1539, 2009.

\bibitem[Amor et~al.(2009)Amor, Marigo, and Maurini]{amor2009regularized}
H.~Amor, J.-J. Marigo, and C.~Maurini.
\newblock Regularized formulation of the variational brittle fracture with
  unilateral contact: Numerical experiments.
\newblock \emph{Journal of the Mechanics and Physics of Solids}, 57\penalty0
  (8):\penalty0 1209 – 1229, 2009.

\bibitem[Ballani et~al.(2013)Ballani, Grasedyck, and Kluge]{ballani2013black}
J.~Ballani, L.~Grasedyck, and M.~Kluge.
\newblock Black box approximation of tensors in hierarchical tucker format.
\newblock \emph{Linear Algebra and its Applications}, 438\penalty0
  (2):\penalty0 639--657, 2013.

\bibitem[Bhattacharjee and Matou{\v{s}}(2016)]{bhattacharjee2016nonlinear}
S.~Bhattacharjee and K.~Matou{\v{s}}.
\newblock A nonlinear manifold-based reduced order model for multiscale
  analysis of heterogeneous hyperelastic materials.
\newblock \emph{Journal of Computational Physics}, 313:\penalty0 635--653,
  2016.

\bibitem[Brepols et~al.(2017)Brepols, Wulfinghoff, and
  Reese]{brepols2017gradient}
T.~Brepols, S.~Wulfinghoff, and S.~Reese.
\newblock Gradient-extended two-surface damage-plasticity: micromorphic
  formulation and numerical aspects.
\newblock \emph{International Journal of Plasticity}, 97:\penalty0 64--106,
  2017.

\bibitem[Brepols et~al.(2018)Brepols, Wulfinghoff, and
  Reese]{brepols2018micromorphic}
T.~Brepols, S.~Wulfinghoff, and S.~Reese.
\newblock A micromorphic damage-plasticity model to counteract mesh dependence
  in finite element simulations involving material softening.
\newblock \emph{Lecture Notes in Applied and Computational Mechanics},
  86:\penalty0 235 – 255, 2018.

\bibitem[Brepols et~al.(2020)Brepols, Wulfinghoff, and
  Reese]{brepols2020gradient}
T.~Brepols, S.~Wulfinghoff, and S.~Reese.
\newblock A gradient-extended two-surface damage-plasticity model for large
  deformations.
\newblock \emph{International Journal of Plasticity}, 129:\penalty0 102635,
  2020.

\bibitem[Breuer and Sirovich(1991)]{breuer1991use}
K.~S. Breuer and L.~Sirovich.
\newblock The use of the karhunen-loeve procedure for the calculation of linear
  eigenfunctions.
\newblock \emph{Journal of Computational Physics}, 96\penalty0 (2):\penalty0
  277--296, 1991.

\bibitem[Chatterjee(2000)]{chatterjee2000introduction}
A.~Chatterjee.
\newblock An introduction to the proper orthogonal decomposition.
\newblock \emph{Current Science}, 78\penalty0 (7):\penalty0 808--817, 2000.

\bibitem[Chaturantabut and Sorensen(2010)]{chaturantabut2010nonlinear}
S.~Chaturantabut and D.~C. Sorensen.
\newblock Nonlinear model reduction via discrete empirical interpolation.
\newblock \emph{SIAM Journal on Scientific Computing}, 32\penalty0
  (5):\penalty0 2737--2764, 2010.

\bibitem[Chen et~al.(2021)Chen, Wang, Hesthaven, and Zhang]{chen2021physics}
W.~Chen, Q.~Wang, J.~S. Hesthaven, and C.~Zhang.
\newblock Physics-informed machine learning for reduced-order modeling of
  nonlinear problems.
\newblock \emph{Journal of Computational Physics}, 446:\penalty0 110666, 2021.

\bibitem[Cordebois and Sidoroff(1982)]{cordebois1982damage}
J.~Cordebois and F.~Sidoroff.
\newblock Damage induced elastic anisotropy.
\newblock In \emph{Mechanical Behavior of Anisotropic Solids/Comportment
  M{\'e}chanique des Solides Anisotropes: Proceedings of the Euromech
  Colloquium 115 Villard-de-Lans, June 19--22, 1979/Colloque Euromech 115
  Villard-de-Lans, 19--22 juin 1979}, pages 761--774. Springer, 1982.

\bibitem[Dinachandra and Alankar(2022)]{dinachandra2022adaptive}
M.~Dinachandra and A.~Alankar.
\newblock Adaptive finite element modeling of phase-field fracture driven by
  hydrogen embrittlement.
\newblock \emph{Computer Methods in Applied Mechanics and Engineering},
  391:\penalty0 114509, 2022.

\bibitem[Dittmann et~al.(2020)Dittmann, Aldakheel, Schulte, Schmidt,
  Kr{\"u}ger, Wriggers, and Hesch]{dittmann2020phase}
M.~Dittmann, F.~Aldakheel, J.~Schulte, F.~Schmidt, M.~Kr{\"u}ger, P.~Wriggers,
  and C.~Hesch.
\newblock Phase-field modeling of porous-ductile fracture in non-linear
  thermo-elasto-plastic solids.
\newblock \emph{Computer Methods in Applied Mechanics and Engineering},
  361:\penalty0 112730, 2020.

\bibitem[Farhat et~al.(2015)Farhat, Chapman, and Avery]{farhat2015structure}
C.~Farhat, T.~Chapman, and P.~Avery.
\newblock Structure-preserving, stability, and accuracy properties of the
  energy-conserving sampling and weighting method for the hyper reduction of
  nonlinear finite element dynamic models.
\newblock \emph{International Journal for Numerical Methods in Engineering},
  102\penalty0 (5):\penalty0 1077--1110, 2015.

\bibitem[Felder et~al.(2022)Felder, Kopic-Osmanovic, Holthusen, Brepols, and
  Reese]{felder2022thermo}
S.~Felder, N.~Kopic-Osmanovic, H.~Holthusen, T.~Brepols, and S.~Reese.
\newblock Thermo-mechanically coupled gradient-extended damage-plasticity
  modeling of metallic materials at finite strains.
\newblock \emph{International Journal of Plasticity}, 148:\penalty0 103142,
  2022.

\bibitem[Forest(2009)]{forest2009micromorphic}
S.~Forest.
\newblock Micromorphic approach for gradient elasticity, viscoplasticity, and
  damage.
\newblock \emph{Journal of Engineering Mechanics}, 135\penalty0 (3):\penalty0
  117--131, 2009.

\bibitem[Freddi and Iurlano(2017)]{fredddi2017numerical}
F.~Freddi and F.~Iurlano.
\newblock Numerical insight of a variational smeared approach to cohesive
  fracture.
\newblock \emph{Journal of the Mechanics and Physics of Solids}, 98:\penalty0
  156--171, 2017.
\newblock ISSN 0022-5096.

\bibitem[Geelen et~al.(2023)Geelen, Wright, and Willcox]{geelen2023operator}
R.~Geelen, S.~Wright, and K.~Willcox.
\newblock Operator inference for non-intrusive model reduction with quadratic
  manifolds.
\newblock \emph{Computer Methods in Applied Mechanics and Engineering},
  403:\penalty0 115717, 2023.

\bibitem[Georgaka et~al.(2018)Georgaka, Stabile, Rozza, and
  Bluck]{georgaka2018parametric}
S.~Georgaka, G.~Stabile, G.~Rozza, and M.~J. Bluck.
\newblock Parametric pod-galerkin model order reduction for unsteady-state heat
  transfer problems.
\newblock \emph{arXiv preprint arXiv:1808.05175}, 2018.

\bibitem[Georgaka et~al.(2020)Georgaka, Stabile, Star, Rozza, and
  Bluck]{georgaka2020hybrid}
S.~Georgaka, G.~Stabile, K.~Star, G.~Rozza, and M.~J. Bluck.
\newblock A hybrid reduced order method for modelling turbulent heat transfer
  problems.
\newblock \emph{Computers \& Fluids}, 208:\penalty0 104615, 2020.

\bibitem[Ghavamian et~al.(2017)Ghavamian, Tiso, and Simone]{ghavamian2017pod}
F.~Ghavamian, P.~Tiso, and A.~Simone.
\newblock Pod--deim model order reduction for strain-softening viscoplasticity.
\newblock \emph{Computer Methods in Applied Mechanics and Engineering},
  317:\penalty0 458--479, 2017.

\bibitem[Golub and Van~Loan(2013)]{golub2013matrix}
G.~H. Golub and C.~F. Van~Loan.
\newblock \emph{Matrix computations}.
\newblock JHU press, 2013.

\bibitem[Gurtin et~al.(2010)Gurtin, Fried, and Anand]{gurtin2010mechanics}
M.~E. Gurtin, E.~Fried, and L.~Anand.
\newblock \emph{The mechanics and thermodynamics of continua}.
\newblock Cambridge university press, 2010.

\bibitem[Hernandez et~al.(2017)Hernandez, Caicedo, and
  Ferrer]{hernandez2017dimensional}
J.~A. Hernandez, M.~A. Caicedo, and A.~Ferrer.
\newblock Dimensional hyper-reduction of nonlinear finite element models via
  empirical cubature.
\newblock \emph{Computer Methods in Applied Mechanics and Engineering},
  313:\penalty0 687--722, 2017.

\bibitem[Holmes(2012)]{holmes2012turbulence}
P.~Holmes.
\newblock \emph{Turbulence, coherent structures, dynamical systems and
  symmetry}.
\newblock Cambridge university press, 2012.

\bibitem[Junker et~al.(2022)Junker, Riesselmann, and
  Balzani]{junker2022efficient}
P.~Junker, J.~Riesselmann, and D.~Balzani.
\newblock Efficient and robust numerical treatment of a gradient-enhanced
  damage model at large deformations.
\newblock \emph{International Journal for Numerical Methods in Engineering},
  123\penalty0 (3):\penalty0 774--793, 2022.

\bibitem[Kastian et~al.(2020)Kastian, Moser, Grasedyck, and
  Reese]{kastian2020two}
S.~Kastian, D.~Moser, L.~Grasedyck, and S.~Reese.
\newblock A two-stage surrogate model for neo-hookean problems based on
  adaptive proper orthogonal decomposition and hierarchical tensor
  approximation.
\newblock \emph{Computer Methods in Applied Mechanics and Engineering},
  372:\penalty0 113368, 2020.

\bibitem[Kastian et~al.(2023)Kastian, Kehls, Brepols, and
  Reese]{kastian2023discrete}
S.~Kastian, J.~Kehls, T.~Brepols, and S.~Reese.
\newblock Discrete empirical interpolation method for nonlinear softening
  problems involving damage and plasticity.
\newblock \emph{arXiv preprint arXiv:2311.17485}, 2023.

\bibitem[Kerfriden et~al.(2011)Kerfriden, Gosselet, Adhikari, and
  Bordas]{kerfriden2011bridging}
P.~Kerfriden, P.~Gosselet, S.~Adhikari, and S.~P.-A. Bordas.
\newblock Bridging proper orthogonal decomposition methods and augmented
  newton--krylov algorithms: an adaptive model order reduction for highly
  nonlinear mechanical problems.
\newblock \emph{Computer Methods in Applied Mechanics and Engineering},
  200\penalty0 (5-8):\penalty0 850--866, 2011.

\bibitem[Kerfriden et~al.(2012)Kerfriden, Passieux, and
  Bordas]{kerfriden2012local}
P.~Kerfriden, J.-C. Passieux, and S.~P.-A. Bordas.
\newblock Local/global model order reduction strategy for the simulation of
  quasi-brittle fracture.
\newblock \emph{International Journal for Numerical Methods in Engineering},
  89\penalty0 (2):\penalty0 154--179, 2012.

\bibitem[Kerfriden et~al.(2013)Kerfriden, Goury, Rabczuk, and
  Bordas]{kerfriden2013partitioned}
P.~Kerfriden, O.~Goury, T.~Rabczuk, and S.~P.-A. Bordas.
\newblock A partitioned model order reduction approach to rationalise
  computational expenses in nonlinear fracture mechanics.
\newblock \emph{Computer Methods in Applied Mechanics and Engineering},
  256:\penalty0 169--188, 2013.

\bibitem[Korelc(2002)]{korelc2002multi}
J.~Korelc.
\newblock Multi-language and multi-environment generation of nonlinear finite
  element codes.
\newblock \emph{Engineering with Computers}, 18:\penalty0 312--327, 2002.

\bibitem[Korelc and Wriggers(2016)]{korelc2016automation}
J.~Korelc and P.~Wriggers.
\newblock \emph{Automation of Finite Element Methods}.
\newblock Springer, 2016.

\bibitem[Lamm et~al.(2024)Lamm, Awad, Jan, Holthusen, Felder, Stefanie, and
  Brepols]{lamm2024gradient}
L.~Lamm, A.~S. Awad, M.~Jan, H.~Holthusen, S.~Felder, R.~Stefanie, and
  T.~Brepols.
\newblock A gradient-extended thermomechanical model for rate-dependent damage
  and failure within rubberlike polymeric materials at finite strains.
\newblock \emph{International Journal of Plasticity}, 173:\penalty0 103883,
  2024.

\bibitem[Lange et~al.(2024)Lange, H{\"u}tter, and Kiefer]{lange2024monolithic}
N.~Lange, G.~H{\"u}tter, and B.~Kiefer.
\newblock A monolithic hyper rom fe2 method with clustered training at finite
  deformations.
\newblock \emph{Computer Methods in Applied Mechanics and Engineering},
  418:\penalty0 116522, 2024.

\bibitem[Liang et~al.(2002)Liang, Lee, Lim, Lin, Lee, and Wu]{liang2002proper}
Y.~Liang, H.~Lee, S.~Lim, W.~Lin, K.~Lee, and C.~Wu.
\newblock Proper orthogonal decomposition and its applications—part i:
  Theory.
\newblock \emph{Journal of Sound and Vibration}, 252\penalty0 (3):\penalty0
  527--544, 2002.

\bibitem[Lu and Pister(1975)]{lu1975decomposition}
S.~Lu and K.~S. Pister.
\newblock Decomposition of deformation and representation of the free energy
  function for isotropic thermoelastic solids.
\newblock \emph{International Journal of Solids and Structures}, 11\penalty0
  (7-8):\penalty0 927--934, 1975.

\bibitem[Miehe et~al.(2010)Miehe, Hofacker, and Welschinger]{miehe2010phase}
C.~Miehe, M.~Hofacker, and F.~Welschinger.
\newblock A phase field model for rate-independent crack propagation: Robust
  algorithmic implementation based on operator splits.
\newblock \emph{Computer Methods in Applied Mechanics and Engineering},
  199\penalty0 (45-48):\penalty0 2765 – 2778, 2010.

\bibitem[Mishra et~al.(2024)Mishra, Carrara, Marfia, Sacco, and
  De~Lorenzis]{mishra2024enhanced}
A.~Mishra, P.~Carrara, S.~Marfia, E.~Sacco, and L.~De~Lorenzis.
\newblock Enhanced transformation field analysis for reduced-order modeling of
  problems with cohesive interfaces.
\newblock \emph{Computer Methods in Applied Mechanics and Engineering},
  421:\penalty0 116755, 2024.

\bibitem[Nguyen et~al.(2020{\natexlab{a}})Nguyen, Zhuang, Chamoin, Zhao,
  Nguyen-Xuan, and Rabczuk]{nguyen2020three}
C.~Nguyen, X.~Zhuang, L.~Chamoin, X.~Zhao, H.~Nguyen-Xuan, and T.~Rabczuk.
\newblock Three-dimensional topology optimization of auxetic metamaterial using
  isogeometric analysis and model order reduction.
\newblock \emph{Computer Methods in Applied Mechanics and Engineering},
  371:\penalty0 113306, 2020{\natexlab{a}}.

\bibitem[Nguyen et~al.(2020{\natexlab{b}})Nguyen, Fassin, Eggersmann, Reese,
  and Wulfinghoff]{nguyen2020gradient}
T.~T. Nguyen, M.~Fassin, R.~Eggersmann, S.~Reese, and S.~Wulfinghoff.
\newblock Gradient-extended brittle damage modeling.
\newblock \emph{Technische Mechanik-European Journal of Engineering Mechanics},
  40\penalty0 (1):\penalty0 53--58, 2020{\natexlab{b}}.

\bibitem[Oliver et~al.(2017)Oliver, Caicedo, Huespe, Hern{\'a}ndez, and
  Roubin]{oliver2017reduced}
J.~Oliver, M.~Caicedo, A.~E. Huespe, J.~Hern{\'a}ndez, and E.~Roubin.
\newblock Reduced order modeling strategies for computational multiscale
  fracture.
\newblock \emph{Computer Methods in Applied Mechanics and Engineering},
  313:\penalty0 560--595, 2017.

\bibitem[Peerlings et~al.(1996)Peerlings, de~Borst, Brekelmans, and
  de~Vree]{peerlings1996gradient}
R.~H. Peerlings, R.~de~Borst, W.~M. Brekelmans, and J.~de~Vree.
\newblock Gradient enhanced damage for quasi-brittle materials.
\newblock \emph{International Journal for Numerical Methods in Engineering},
  39\penalty0 (19):\penalty0 3391--3403, 1996.

\bibitem[Radermacher and Reese(2013)]{radermacher2013comparison}
A.~Radermacher and S.~Reese.
\newblock A comparison of projection-based model reduction concepts in the
  context of nonlinear biomechanics.
\newblock \emph{Archive of Applied Mechanics}, 83\penalty0 (8):\penalty0
  1193--1213, 2013.

\bibitem[Radermacher and Reese(2016)]{radermacher2016pod}
A.~Radermacher and S.~Reese.
\newblock Pod-based model reduction with empirical interpolation applied to
  nonlinear elasticity.
\newblock \emph{International Journal for Numerical Methods in Engineering},
  107\penalty0 (6):\penalty0 477--495, 2016.

\bibitem[Rezaei et~al.(2017)Rezaei, Wulfinghoff, and
  Reese]{rezaei2017prediction}
S.~Rezaei, S.~Wulfinghoff, and S.~Reese.
\newblock Prediction of fracture and damage in micro/nano coating systems using
  cohesive zone elements.
\newblock \emph{International Journal of Solids and Structures}, 121:\penalty0
  62--74, 2017.

\bibitem[Rocha et~al.(2020)Rocha, van~der Meer, Moror{\'o}, and
  Sluys]{rocha2020accelerating}
I.~B. Rocha, F.~P. van~der Meer, L.~A. Moror{\'o}, and L.~J. Sluys.
\newblock Accelerating crack growth simulations through adaptive model order
  reduction.
\newblock \emph{International Journal for Numerical Methods in Engineering},
  121\penalty0 (10):\penalty0 2147--2173, 2020.

\bibitem[Ruan et~al.(2023)Ruan, Rezaei, Yang, Gross, and Xu]{ruan2023thermo}
H.~Ruan, S.~Rezaei, Y.~Yang, D.~Gross, and B.-X. Xu.
\newblock A thermo-mechanical phase-field fracture model: Application to hot
  cracking simulations in additive manufacturing.
\newblock \emph{Journal of the Mechanics and Physics of Solids}, 172:\penalty0
  105169, 2023.

\bibitem[Saji et~al.(2024)Saji, Pantidis, and Mobasher]{saji2024new}
R.~P. Saji, P.~Pantidis, and M.~E. Mobasher.
\newblock A new unified arc-length method for damage mechanics problems.
\newblock \emph{Computational Mechanics}, pages 1--32, 2024.

\bibitem[Selvaraj and Hallett(2024)]{selvaraj2024adaptive}
J.~Selvaraj and S.~R. Hallett.
\newblock Adaptive and variable model order reduction method for fracture
  modelling using explicit time integration.
\newblock \emph{Computer Methods in Applied Mechanics and Engineering},
  418:\penalty0 116506, 2024.

\bibitem[Simo and Miehe(1992)]{simo1992associative}
J.~Simo and C.~Miehe.
\newblock Associative coupled thermoplasticity at finite strains: Formulation,
  numerical analysis and implementation.
\newblock \emph{Computer Methods in Applied Mechanics and Engineering},
  98\penalty0 (1):\penalty0 41--104, 1992.

\bibitem[Stojanovic et~al.(1964)Stojanovic, Djuric, and
  Vujosevic]{stojanovic1964finite}
R.~Stojanovic, S.~Djuric, and L.~Vujosevic.
\newblock On finite thermal deformations.
\newblock \emph{Archiwum Mechaniki Stosowanej}, 16:\penalty0 103--108, 1964.

\bibitem[Taylor(2014)]{taylor2014feap}
R.~L. Taylor.
\newblock Feap-a finite element analysis program, 2014.

\bibitem[Tiso and Rixen(2013)]{tiso2013discrete}
P.~Tiso and D.~J. Rixen.
\newblock Discrete empirical interpolation method for finite element structural
  dynamics.
\newblock In \emph{Topics in Nonlinear Dynamics, Volume 1: Proceedings of the
  31st IMAC, A Conference on Structural Dynamics, 2013}, pages 203--212.
  Springer, 2013.

\bibitem[Volkwein(2013)]{volkwein2013proper}
S.~Volkwein.
\newblock Proper orthogonal decomposition: Theory and reduced-order modelling.
\newblock \emph{Lecture Notes, University of Konstanz}, 4\penalty0
  (4):\penalty0 1--29, 2013.

\bibitem[Wriggers et~al.(2021{\natexlab{a}})Wriggers, {De Bellis}, and
  Hudobivnik]{wriggers2021taylor}
P.~Wriggers, M.~{De Bellis}, and B.~Hudobivnik.
\newblock A taylor–hood type virtual element formulations for large
  incompressible strains.
\newblock \emph{Computer Methods in Applied Mechanics and Engineering},
  385:\penalty0 114021, 2021{\natexlab{a}}.

\bibitem[Wriggers et~al.(2021{\natexlab{b}})Wriggers, Hudobivnik, and
  Aldakheel]{wriggers2021nurbs}
P.~Wriggers, B.~Hudobivnik, and F.~Aldakheel.
\newblock Nurbs-based geometries: A mapping approach for virtual serendipity
  elements.
\newblock \emph{Computer Methods in Applied Mechanics and Engineering},
  378:\penalty0 113732, 2021{\natexlab{b}}.

\bibitem[Wu(2017)]{wu2017unified}
J.-Y. Wu.
\newblock A unified phase-field theory for the mechanics of damage and
  quasi-brittle failure.
\newblock \emph{Journal of the Mechanics and Physics of Solids}, 103:\penalty0
  72--99, 2017.

\bibitem[Zhang et~al.(2022{\natexlab{a}})Zhang, Lu, Tang, and
  Liu]{zhang2022hidenn}
L.~Zhang, Y.~Lu, S.~Tang, and W.~K. Liu.
\newblock Hidenn-td: Reduced-order hierarchical deep learning neural networks.
\newblock \emph{Computer Methods in Applied Mechanics and Engineering},
  389:\penalty0 114414, 2022{\natexlab{a}}.

\bibitem[Zhang et~al.(2021)Zhang, Mortazavi, and Aldakheel]{zhang2021molecular}
Q.~Zhang, B.~Mortazavi, and F.~Aldakheel.
\newblock Molecular dynamics modeling of mechanical properties of polymer
  nanocomposites reinforced by c7n6 nanosheet.
\newblock \emph{Surfaces}, 4\penalty0 (3):\penalty0 240--254, 2021.

\bibitem[Zhang et~al.(2022{\natexlab{b}})Zhang, Mortazavi, Zhuang, and
  Aldakheel]{zhang2022exploring}
Q.~Zhang, B.~Mortazavi, X.~Zhuang, and F.~Aldakheel.
\newblock Exploring the mechanical properties of two-dimensional carbon-nitride
  polymer nanocomposites by molecular dynamics simulations.
\newblock \emph{Composite Structures}, 281:\penalty0 115004,
  2022{\natexlab{b}}.

\bibitem[Zhuang et~al.(2022)Zhuang, Zhou, Huynh, Areias, and
  Rabczuk]{zhuang2022phase}
X.~Zhuang, S.~Zhou, G.~D.~H. Huynh, P.~Areias, and T.~Rabczuk.
\newblock Phase field modeling and computer implementation: A review.
\newblock \emph{Engineering Fracture Mechanics}, 262:\penalty0 108234, 2022.

\bibitem[Zhuang et~al.(2023)Zhuang, Li, and Zhou]{zhuang2023transverse}
X.~Zhuang, X.~Li, and S.~Zhou.
\newblock Transverse penny-shaped hydraulic fracture propagation in
  naturally-layered rocks under stress boundaries: A 3d phase field modeling.
\newblock \emph{Computers and Geotechnics}, 155:\penalty0 105205, 2023.

\end{thebibliography}


@article{clausius1870xvi,
  title={XVI. On a mechanical theorem applicable to heat},
  author={Clausius, Rudolf},
  journal={The London, Edinburgh, and Dublin Philosophical Magazine and Journal of Science},
  volume={40},
  number={265},
  pages={122--127},
  year={1870},
  publisher={Taylor \& Francis}
}

@article{maxwell1870reciprocal,
  title={I.—on reciprocal figures, frames, and diagrams of forces},
  author={Maxwell, J Clerk},
  journal={Earth and Environmental Science Transactions of the Royal Society of Edinburgh},
  volume={26},
  number={1},
  pages={1--40},
  year={1870},
  publisher={Royal Society of Edinburgh Scotland Foundation}
}

@article{novoselov2004electric,
  title={Electric field effect in atomically thin carbon films},
  author={Novoselov, Kostya S and Geim, Andre K and Morozov, Sergei V and Jiang, De-eng and Zhang, Yanshui and Dubonos, Sergey V and Grigorieva, Irina V and Firsov, Alexandr A},
  journal={science},
  volume={306},
  number={5696},
  pages={666--669},
  year={2004},
  publisher={American Association for the Advancement of Science}
}

@article{geim2007rise,
  title={The rise of graphene},
  author={Geim, Andre K and Novoselov, Konstantin S},
  journal={Nature materials},
  volume={6},
  number={3},
  pages={183--191},
  year={2007},
  publisher={Nature Publishing Group UK London}
}

@article{schedin2007detection,
  title={Detection of individual gas molecules adsorbed on graphene},
  author={Schedin, Fredrik and Geim, Andrei Konstantinovich and Morozov, Sergei Vladimirovich and Hill, Ew W and Blake, Peter and Katsnelson, Mi I and Novoselov, Kostya Sergeevich},
  journal={Nature materials},
  volume={6},
  number={9},
  pages={652--655},
  year={2007},
  publisher={Nature Publishing Group UK London}
}

@article{neto2009electronic,
  title={The electronic properties of graphene},
  author={Neto, AH Castro and Guinea, Francisco and Peres, Nuno MR and Novoselov, Kostya S and Geim, Andre K},
  journal={Reviews of modern physics},
  volume={81},
  number={1},
  pages={109},
  year={2009},
  publisher={APS}
}
\end{document}